\documentclass[reqno]{amsart}

\usepackage{xypic}
\input xy
\xyoption{all}
\usepackage{epsfig}
\usepackage{color}
\usepackage{amsthm}
\usepackage{changepage}
\usepackage{amssymb}
\usepackage{mathdots}
\usepackage{amsmath}
\usepackage{amscd}
\usepackage{amsopn}
\usepackage{enumerate}
\usepackage{mathtools}
\usepackage{url}
\usepackage{graphicx}
\usepackage{hyperref}

\usepackage{hyperref}\hypersetup{colorlinks}

\usepackage{color} 

\definecolor{darkred}{rgb}{1,0,0} 
\definecolor{darkgreen}{rgb}{0,0.8,0}
\definecolor{darkblue}{rgb}{0,0,1}

\hypersetup{colorlinks,
linkcolor=darkblue,
filecolor=darkgreen,
urlcolor=darkred,
citecolor=darkgreen}

\numberwithin{equation}{section}

\newtheorem {Theorem}{Theorem}

\numberwithin{Theorem}{section}

\newtheorem {Lemma}[Theorem]    {Lemma}

\newtheorem {Proposition}[Theorem]{Proposition}
\newtheorem {Corollary}[Theorem]{Corollary}
\theoremstyle{definition}
\newtheorem{Definition}[Theorem]{Definition}

\newtheorem{Remark}[Theorem]{Remark}
\newtheorem{Example}[Theorem]{Example}
\theoremstyle{remark}

\newtheoremstyle{TheoremForIntro} 
        {.6em}{.6em}              
        {\itshape}                      
        {}                              
        {\bfseries}                     
        {. }                             
        { }                             
        {\thmname{#1}\thmnote{ \bfseries #3}}
    \theoremstyle{TheoremForIntro}
    \newtheorem{TheoremIntro}[Theorem]{Theorem}
    \newtheorem{PropositionIntro}[Theorem]{Proposition}

\expandafter\chardef\csname pre amssym.def at\endcsname=\the\catcode`\@
\catcode`\@=11
\def\undefine#1{\let#1\undefined}
\def\newsymbol#1#2#3#4#5{\let\next@\relax
 \ifnum#2=\@ne\let\next@\msafam@\else
 \ifnum#2=\tw@\let\next@\msbfam@\fi\fi
 \mathchardef#1="#3\next@#4#5}
\def\mathhexbox@#1#2#3{\relax
 \ifmmode\mathpalette{}{\m@th\mathchar"#1#2#3}%
 \else\leavevmode\hbox{$\m@th\mathchar"#1#2#3$}\fi}
\def\hexnumber@#1{\ifcase#1 0\or 1\or 2\or 3\or 4\or 5\or 6\or 7\or 8\or
 9\or A\or B\or C\or D\or E\or F\fi}

\font\teneufm=eufm10
\font\seveneufm=eufm7
\font\fiveeufm=eufm5
\newfam\eufmfam
\textfont\eufmfam=\teneufm
\scriptfont\eufmfam=\seveneufm
\scriptscriptfont\eufmfam=\fiveeufm

\catcode`\@=\csname pre amssym.def at\endcsname

\newcommand{\fsu}{{\mathfrak {su}}}
\newcommand{\fsl}{{\mathfrak {sl}}}
\newcommand{\fgl}{{\mathfrak {gl}}}

\newcommand{\fo}{{\mathfrak {o}}}
\newcommand{\fsp}{{\mathfrak {sp}}}
\newcommand{\fa}{{\mathfrak a}}
\newcommand{\fg}{{\mathfrak g}}
\newcommand{\fh}{{\mathfrak h}}
\newcommand{\fn}{{\mathfrak n}}
\newcommand{\fs}{{\mathfrak s}}
\newcommand{\fc}{{\mathfrak c}}

\newcommand{\ft}{{\mathfrak t}}

\newcommand{\fm}{{\mathfrak m}}
\newcommand{\fk}{{\mathfrak k}}

\newcommand{\sSU}{{\mathsf {SU}}}
\newcommand{\sSL}{{\mathsf {SL}}}
\newcommand{\sSO}{{\mathsf {SO}}}
\newcommand{\sO}{{\mathsf {O}}}
\newcommand{\sSp}{{\mathsf {Sp}}}
\newcommand{\sGL}{{\mathsf {GL}}}

\newcommand{\sG}{{\mathsf G}}
\newcommand{\sE}{{\mathsf E}}
\newcommand{\sH}{{\mathsf H}}

\newcommand{\sC}{{\mathsf C}}

\newcommand{\sT}{{\mathsf T}}

\newcommand{\sK}{{\mathsf K}}
\newcommand{\sU}{{\mathsf U}}

\newcommand{\Xx}{{\mathcal X}}

\newcommand{\Ll}{{\mathcal L}}

\newcommand{\Ee}{{\mathcal E}}
\newcommand{\Ff}{{\mathcal F}}

\newcommand{\Gg}{{\mathcal G}}

\newcommand{\Mm}{{\mathcal M}}

\newcommand{\Oo}{{\mathcal O}}
\newcommand{\Pp}{{\mathcal P}}
\newcommand{\Rr}{{\mathcal R}}

\newcommand{\Tt}{{\mathcal T}}

\newcommand{\gen}[1]{\left\langle  #1 \right\rangle}

\newcommand{\dwn}[1]{_{_{#1}}}

\newcommand{\mtrx}[1]{\left (\begin{matrix}#1\end{matrix}\right)}
\newcommand{\eqtns}[2]{\begin{equation}
\begin{dcases}#1\end{dcases}\label{#2}\end{equation}}

\def    \C      {{\mathbb C}}
\def    \R      {{\mathbb R}}
\def    \Z      {{\mathbb Z}}
\def    \N      {{\mathbb N}}

\def    \P    {{\mathbb P}}

\def    \ra     {{\rightarrow}}
\def    \lra     {{\longrightarrow}}
\def    \haf    {{\frac{1}{2}}}

\def    \p      {\partial}

\def    \H  {\operatorname{\scriptscriptstyle{H}}}

\newcommand{\An}{\xymatrix{ *{\circ} \ar@{-}[r]|*\dir{ } & *{\circ}\ar@{-}[r]&{\cdots}&*{\circ}\ar@{-}[l]|*\dir{ }\ar@{-}[r]|*\dir{ }&*{\circ} }}
\newcommand{\Anlabel}{\xymatrix@R=.25em{ *{\circ}\ar@<-1ex>@{}[d]^{\alpha_{1}} \ar@{-}[r]|*\dir{ } & *{\circ}\ar@<-1ex>@{}[d]^{\alpha_{2}}\ar@{-}[r]&{\cdots}&*{\circ}\ar@{-}[l]|*\dir{ }\ar@{-}[r]|*\dir{ }\ar@<-2ex>@{}[d]^{\alpha_{n-1}}&*{\circ}\ar@<-1ex>@{}[d]^{\alpha_{n}}\\&&&& }}
\newcommand{\Bn}{\xymatrix{ *{\circ} \ar@{-}[r]|*\dir{ } & *{\circ}\ar@{-}[r]&{\cdots}&*{\circ}\ar@{-}[l]|*\dir{ }\ar@{=}[r]|*\dir{>}&*{\circ} }}
\newcommand{\Bnlabel}{\xymatrix@R=.25em{ *{\circ}\ar@<-1ex>@{}[d]^{\alpha_{1}} \ar@{-}[r]|*\dir{ } & *{\circ}\ar@<-1ex>@{}[d]^{\alpha_{2}}\ar@{-}[r]&{\cdots}&*{\circ}\ar@{-}[l]|*\dir{ }\ar@{=}[r]|*\dir{>}\ar@<-2ex>@{}[d]^{\alpha_{n-1}}&*{\circ}\ar@<-1ex>@{}[d]^{\alpha_{n}}\\&&&& }}
\newcommand{\Cn}{\xymatrix{ *{\circ} \ar@{-}[r]|*\dir{ } & *{\circ}\ar@{-}[r]&{\cdots}&*{\circ}\ar@{-}[l]|*\dir{ }\ar@{=}[r]|*\dir{<}&*{\circ} }}
\newcommand{\Cnlabel}{\xymatrix@R=.25em{ *{\circ}\ar@<-1ex>@{}[d]^{\alpha_{1}} \ar@{-}[r]|*\dir{ } & *{\circ}\ar@<-1ex>@{}[d]^{\alpha_{2}}\ar@{-}[r]&{\cdots}&*{\circ}\ar@{-}[l]|*\dir{ }\ar@{=}[r]|*\dir{<}\ar@<-2ex>@{}[d]^{\alpha_{n-1}}&*{\circ}\ar@<-1ex>@{}[d]^{\alpha_{n}}\\&&&& }}
\newcommand{\Dn}{\xymatrix@R=.25em{&&&&*{\circ} \\ *{\circ} \ar@{-}[r]|*\dir{ } & *{\circ}\ar@{-}[r]&{\cdots}&*{\circ}\ar@{-}[l]|*\dir{ }\ar@{-}[ur]|*\dir{ }\ar@{-}[dr]|*\dir{ }& \\ &&&&*{\circ} }}
\newcommand{\Dnlabel}{\xymatrix@R=.25em{&&&&*{\circ}\ar@<-1ex>@{}[d]^{\alpha_{n-1}} \\ *{\circ}\ar@<-1ex>@{}[d]^{\alpha_{1}} \ar@{-}[r]|*\dir{ } & *{\circ}\ar@<-1ex>@{}[d]^{\alpha_{2}}\ar@{-}[r]&{\cdots}&*{\circ}\ar@{-}[l]|*\dir{ }\ar@{-}[ur]|*\dir{ }\ar@{-}[dr]|*\dir{ }\ar@<-4ex>@{}[d]^{\alpha_{n-2}}& \\ &&&&*{\circ}\ar@<-2ex>@{}[u]^(-.5){\alpha_{n}} }}

\usepackage[margin=1in]{geometry}
\author{Brian Collier}

\title{Maximal $\sSp(4,\R)$ surface group representations, minimal immersions and cyclic surfaces}

\begin{document}

\setlength{\smallskipamount}{6pt}
\setlength{\medskipamount}{10pt}
\setlength{\bigskipamount}{16pt}

\maketitle
\begin{abstract}
	Let $S$ be a closed surface of genus at least $2$. For each maximal representation $\rho: \pi_1(S)\ra\sSp(4,\R)$ in one of the $2g-3$ exceptional connected components, we prove there is a unique conformal structure on the surface in which the corresponding equivariant harmonic map to the symmetric space $\sSp(4,\R)/\sU(2)$ is a minimal immersion. Using a Higgs bundle parameterization of these components, we give a mapping class group invariant parameterization of such components as fiber bundles over Teichm\"uller space.
	Unlike Labourie's recent results on Hitchin components, these bundles are not vector bundles.  
\end{abstract}
\tableofcontents

\section{Introduction}
\smallskip
The nonabelian Hodge correspondence between Higgs bundles and surface group representations provides powerful tools for studying surface group character varieties, but these tools come at a price. The correspondence depends on a choice of conformal (or, equivalently, complex) structure on the surface, and, as a result, breaks the mapping class group symmetry.
 In this paper, we restore the symmetry by associating a unique ``preferred'' conformal structure (Theorem \ref{uniqueminThm}) to each point in certain connected components of the surface group character variety for $\sSp(4,\R).$  Using a Higgs bundle parameterization (Proposition \ref{ParamMdProposition}), we give a mapping class group invariant parameterization of these components (Theorem \ref{MCGinvparamTheorem}). 
 Recently, Labourie \cite{cyclicSurfacesRank2} provided a mapping class group invariant parameterization of the Hitchin component for $\sSp(4,\R).$ Our approach builds on his work and provides the first mapping class group parameterization of what we call the \em Gothen components.\em

For a a closed oriented surface $S$ of genus $g\geq2,$ denote the associated character variety by $\Rr(\sSp(4,\R)),$
 \[\Rr(\sSp(4,\R))\ =\ Hom(\pi_1(S),\sSp(4,\R))\ //\ \sSp(4,\R).\]
Since $\sSp(4,\R)$ is a group of Hermitian type, there is a Toledo invariant $\tau\in\Z$ which helps to distinguish connected coponents of $\Rr(\sSp(4,\R)).$ 
Given two representations $\rho$ and $\rho'$ in $\Rr(\sSp(4,\R)),$ if $\tau(\rho)\neq\tau(\rho'),$ then $\rho$ and $\rho'$ are in different connected components. 
The Toledo invariant satisfies the Milnor-Wood inequality $|\tau|\leq -\haf\chi(S)\cdot rk(\sSp(4,\R))= 2g-2$ \cite{TuraevMilnorWoodInq}, giving a decomposition:
\begin{equation}\label{ToledoRepDecomp}
	\Rr(\sSp(4,\R))=\bigsqcup\limits_{-(2g-2)\leq\tau\leq2g-2}\Rr^\tau(\sSp(4,\R)),
\end{equation}
with $\Rr^\tau(\sSp(4,\R))\cong\Rr^{-\tau}(\sSp(4,\R)).$  

A conformal structure $J$ on the surface $S$ defines a Riemann surface structure $(S,J)=\Sigma.$ 
For each conformal structure, let $\Mm_J(\sSp(4,\R))$ be the moduli space of $\sSp(4,\R)$-Higgs bundles over $\Sigma.$ 
By the nonabelian Hodge correspondence \cite{canonicalmetrics,harmoicmetric,selfduality,localsystems}, there is a homeomorphism between the representation variety $\Rr(\sSp(4,\R))$ and the Higgs bundle moduli space $\Mm_J(\sSp(4,\R)).$
Using this correspondence, the Toledo invariant gives a decomposition:
\[\xymatrix{\Mm_J(\sSp(4,\R))=\bigsqcup\limits_{-(2g-2)\leq\tau\leq2g-2}\Mm_{J,\tau}(\sSp(4,\R))}.\]
Such a decomposition can be seen directly using Higgs bundle techniques \cite{sp4GothenConnComp,HiggsbundlesSP2nR}.

In \cite{sp4GothenConnComp}, Gothen showed that, for $|\tau|=2g-2$, $\Mm_{J,\tau}(\sSp(4,\R))$ has $3\cdot 2^{2g}+2g-4$ connected components. 
When $|\tau|=2g-2$, representations $\rho\in\Rr^\tau(\sSp(4,\R))$ are called \em maximal \em and have many special properties. 
For example, $2^{2g}$ of the $3\cdot 2^{2g}+2g-4$ connected components are comprised of Hitchin representations \cite{liegroupsteichmuller}. 
Maximal representations are examples of Anosov representations \cite{MaxRepsAnosov}, and thus, they are all discrete and faithful. 
Furthermore, Labourie has shown \cite{CrossRatioAnosoveProperEnergy} that the mapping class group acts properly on the space of maximal $\sSp(4,\R)$ representations. 
From a detailed study of the geometry of the symmetric space $\sSp(4,\R)/\sU(2),$ Burger, Iozzi, and Wienhard  \cite{BIWmaximalToledoAnnals} classified the possible subgroups $\sG\subset\sSp(4,\R)$ which can be the Zariski closure of a maximal $\sSp(4,\R)$ representation. 
Using Higgs bundles, Bradlow, Garcia-Prada, and Gothen \cite{MaximalSP4} showed that $\Rr^{2g-2}(\sSp(4,\R))$ has $2^{2g}+2g-3$ smooth connected components which are distinguished by an integer $g-1<d\leq 3g-3:$ 
\[\xymatrix@=1.5em{\bigsqcup\limits_{g-1<d\leq 3g-3} \Mm_{J}^{d}(\sSp(4,\R))\subset \Mm_{J,{2g-2}}(\sSp(4,\R))&\text{and}&\bigsqcup\limits_{g-1<d\leq 3g-3} \Rr_{d}(\sSp(4,\R))\subset \Rr^{2g-2}(\sSp(4,\R))},\]
with $\Mm_J^d(\sSp(4,\R))\cong \Rr_d(\sSp(4,\R)).$

When $d=3g-3,$ $\Rr_d(\sSp(4,\R))$ has $2^{2g}$ connected components which are all isomorphic; these are the Hitchin components. They are contractible and contain representations with Zariski closure the irreducible $\sSL(2,\R)\subset\sSp(4,\R).$ 
For $g-1<d<3g-3$ the smooth components $\Rr_d(\sSp(4,\R))$ are connected, noncontractible, and contain only Zariski dense representations. As these $2g-3$ connected components were originally discovered by Gothen, we call them the \em Gothen components\em. Although all representations in the Gothen components are Zariski dense, Guichard and Wienhard \cite{TopInvariantsAnosov} have identified some special representations in these components, called hybrid representations, which are built by gluing together Fuchsian representation on subsurfaces of $S$.
In Theorem \ref{MCGinvparamTheorem}, we give a mapping class group invariant parameterization of these $2^{2g}+2g-3$ smooth connected components $\Rr_d(\sSp(4,\R))$ as fiber bundles over Teichm\"uller space. This is done by associating a preferred conformal structure to each representation $\rho\in\Rr_d(\sSp(4,\R))$ (Theorem \ref{uniqueminThm}). 

To make this more precise, recall that if $\sH\subset\sG$ is the maximal compact subgroup of a semisimple Lie group and $\rho:\pi_1(S)\ra\sG$ is a representation, then a metric on the flat $\sG$-bundle $\widetilde \Sigma\times_\rho \sG$ is equivalent to a $\rho$-equivariant map
$h_\rho:\widetilde S\lra\sG/\sH.$
If $\rho$ is completely reducible, then, by Corlette's theorem \cite{canonicalmetrics}, for each choice of conformal structure $\Sigma=(S,J)$ there exists a \em harmonic \em metric on the associated flat bundle. That is, the $\rho$-equivariant map $h_\rho:\widetilde \Sigma\ra \sG/\sH$ can be chosen to be a critical point of the energy functional
\[\Ee_J(h_\rho)=\haf\int\limits_\Sigma|dh_\rho|^2 dvol.\]
For any harmonic map $h:\Sigma\ra (N,g),$ with domain a Riemann surface, the $(2,0)$-part of the pullback metric $h^*g$ defines a holomorphic quadratic differential $q_2(h)\in H^0(\Sigma, K^2),$ called the Hopf differential:
\[q_2(h)=(h^*g)^{2,0}.\]
Furthermore, for a  harmonic map $h,$ the Hopf differential is zero if and only if $h$ is weakely conformal or equivalently, a branched minimal immersion \cite{MinImmofRiemannSurf,SchoenYauMinimalSurfEnergy}. 

\begin{TheoremIntro}[\ref{uniqueminThm}]
 	If $S$ be a closed surface of genus at least 2 and $\rho\in\Rr_d(\sSp(4,\R))$ for $g-1<d\leq 3g-3,$ then there exists a unique conformal structure $(S,J_\rho)=\Sigma$ so that the unique $\rho$-equivariant harmonic map 
 	\[h_\rho:\widetilde\Sigma\ra\sSp(4,\R)/\sU(2)\] is a minimal immersion. 
 \end{TheoremIntro}
%

\begin{Remark}
	Independent of the component, for all $\rho\in\Rr^{2g-2}(\sSp(4,\R))$ there exists a conformal structure so that the corresponding $\rho$-equivariant harmonic map $h_\rho$ is a branched minimal immersion. 
	This follows from the properness of the mapping class group action on $\Rr^{2g-2}(\sSp(4,\R))$ \cite{CrossRatioAnosoveProperEnergy}. For the components $\Rr_d(\sSp(4,\R))$ with $g-1<d\leq 3g-3,$ we show that the conformal structure is unique, and that the associated branched minimal immersion is branch point free. For the Hitchin components, Theorem \ref{uniqueminThm} was proven by Labourie \cite{cyclicSurfacesRank2} but is new for the Gothen components.
\end{Remark}

The proof of Theorem \ref{uniqueminThm} generalizes techniques of \cite{cyclicSurfacesRank2}, together with a Higgs bundle description of the $2g-3$ non-Hitchin smooth components $\Mm_J^d(\sSp(4,\R))\subset\Mm^{2g-2}_J(\sSp(4,\R))$. 
Fix a representation $\rho\in\Rr_d(\sSp(4,\R))$ and choose a conformal structure $\Sigma=(S,J)$ in which the equivariant harmonic map $h_\rho$ is minimal. 
The minimal map condition is interpreted as the vanishing of the holomorphic quadratic differential associated, via the Hitchin fibration, to the Higgs bundle corresponding to $\rho.$ 
In Proposition \ref{MetricSplittingProposition}, we show that these Higgs bundles have the special property that they are fixed by a $4^{th}$-roots of unity action on $\Mm(\sSp(4,\R)),$ and thus, the harmonic map $h_\rho$ lifts (see Remark \ref{metricliftingRemark}) equivariantly, to the homogeneous space $\sSp(4,\R)/(\sU(1)\times \sU(1))$
\[
\xymatrix{
&\sSp(4,\R)/(\sU(1)\times \sU(1))\ar[d]\\
\widetilde\Sigma\ar[r]_{h_\rho\ \ \ \ }\ar@{-->}[ur]&\sSp(4,\R)/\sU(2)
}
\]
For the Hitchin component, such Higgs bundles and their metric splitting properties were studied in detail by Baraglia \cite{g2geometry,cyclichiggsaffinetoda}. 
To show local uniqueness, the lifted harmonic maps are interpreted as a \em cyclic Pfaffian system \em (see Defintion \ref{GG0cyclicsurfaceDef}) which generalizes the cyclic surfaces of \cite{cyclicSurfacesRank2}, and shown to be rigid in Theorem \ref{RigidityOfCyclicSurf}. 
Then, using the smoothness of the connected components, we appeal to an argument in \cite{cyclicSurfacesRank2} to obtain global uniqueness. 
As a special case, we recover the analogous result for the $\sSp(4,\R)$-Hitchin component.
\begin{Remark}
The local uniqueness of Theorem \ref{RigidityOfCyclicSurf} holds for a slightly more general class of representations. For $d=g-1,$ there is a connected component $\Rr_{g-1}(\sSp(4,\R))\subset\Rr^{2g-2}(\sSp(4,\R))$ which shares many properties with $\Rr_d(\sSp(4,\R))$ for $g-1<d\leq 3g-3$ (see Remark \ref{StabilityofMdRemark}). 
\end{Remark}
To describe the mapping class group invariant parameterization of the components $\Rr_d(\sSp(4,\R))$ for $g-1<d\leq3g-3$, we need to describe the Higgs bundles in $\Mm_J^d(\sSp(4,\R))$. 
Let $K\ra\Sigma$ is the canonical bundle, Higgs bundles in $\Mm_J^d(\sSp(4,\R))$ are determined by quadruples $(N,\mu,\nu,q_2)$ \cite{MaximalSP4}, where:
\begin{itemize}
	\item $N\ra\Sigma$ is a holomorphic line bundle of degree $d,$
	\item $\mu\in H^0(N^{-2}K^3)$ is a \em nonzero \em holomorphic section,
	\item $\nu\in H^0(N^2K)$ is any holomorphic section,
	\item $q_2\in H^0(K^2)$ is a holomorphic quadratic differential.
\end{itemize} 
In the above description, the tuple $(N,\mu,\nu,q_2)$ does not uniquely determine the isomorphism class of the associated Higgs bundle; instead, there is a $1$-parameter gauge symmetry to account for. For $\lambda\in \C^*,$ there is an isomorphism 
\begin{equation}\label{fiberwiseactionIntro}
	\xymatrix{(N,\mu,\nu,q_2)\ar[r]^{g_\lambda\ \ \ \ \ }&(N,\lambda^{-2}\mu,\lambda^2\nu,q_2)}.
\end{equation}
Denote the variety of degree $d$ line bundles on $\Sigma$ by $Pic^d(\Sigma).$ 
For $g-1<d\leq3g-3,$ the moduli space $\Mm_J^d(\sSp(4,\R))$ is given by 
\[\left.\left\{(N,\mu,\nu,q_2)\ |\ N\in Pic^d(\Sigma),\ 0\neq\mu\in H^0(N^{-2}K^3),\ \nu\in H^0(N^2K),\ q_2\in H^0(K^2)\right\}\right/\C^*,\]
where $\C^*$ acts by \eqref{fiberwiseactionIntro}. 
Since the holomorphic quadratic differential is independent of $N$ and is acted on trivially by \eqref{fiberwiseactionIntro}, define
\begin{equation}
	\Ff_J^d=\left.\left\{(N,\mu,\nu)|N\in Pic^d(\Sigma)\ \text{with}\ h^0(N^{-2}K^3)>0,\ 0\neq\mu\in H^0(N^{-2}K^3),\ \nu\in H^0(N^2K)\right\}\right/\C^*\label{FJdeq}
\end{equation}
The moduli space $\Mm_J^d$ is then given by 
\[\Mm_J^d(\sSp(4,\R))\cong\Ff_J^d\times H^0(K^2).\]
Using Theorem \ref{uniqueminThm} we prove the following mapping class group invariant parameterization of $\Rr_d(\sSp(4,\R)).$
\begin{TheoremIntro}[\ref{MCGinvparamTheorem}]
 	For $g-1<d\leq3g-3,$ let $\Rr_d(\sSp(4,\R))$ be the component of the maximal $\sSp(4,\R)$ representation variety corresponding to the Higgs bundle component $\Mm_J^d(\sSp(4,\R)).$ There is a mapping class group equivariant diffeomorphism 
 	\[\Psi:\Ff^d\lra\Rr_d(\sSp(4,\R))\]
 	where $\pi:\Ff^d\ra\Tt(S)$ is a fiber bundle over Teichm\"uller space with $\pi^{-1}([J])= \Ff_J^d$ from \eqref{FJdeq}.
 \end{TheoremIntro}
\begin{Remark}
	Since Teichm\"uller space is contractible, the bundle $\Ff\ra\Tt(S)$ is trivial and $\Rr_d(\sSp(4,\R))$ is homotopic to the fiber $\Ff_J^d.$ 
	\end{Remark}

To describe the fiber $\Ff_J^d$ in more detail we follow ideas of \cite{MaximalSP4}. Let $\Pp_d\ra Pic^d(\Sigma)\times \Sigma$ be the Poincar\'e line bundle; $\Pp_d$ is constructed so that the line over a point $(N,p)\in Pic^d(\Sigma)\times\Sigma$ is the line $N_p.$ 
Consider the line bundles 
\[\xymatrix{\Ll_1=\Pp_d^{-2}\otimes\pi_\Sigma^* K^3&\text{and} &\Ll_2=\Pp_d^{2}\otimes\pi_\Sigma^* K}\]
on $Pic^d(\Sigma)\times\Sigma,$ here $\pi_{Pic}$ and $\pi_{\Sigma}$ are the projections onto the appropriate factors.
Pushing forward gives coherent sheaves on $Pic^d(\Sigma)$ 
\[\xymatrix{(\pi_{Pic})_*\Ll_1&\text{and}&(\pi_{Pic})_*\Ll_2}.\]
The stalks of $(\pi_{Pic})_*\Ll_1$ and $(\pi_{Pic})_*\Ll_2$ over a line bundle $N\in Pic^d(\Sigma)$ are $H^0(N^{-2}K^3)$ and $H^0(N^2K)$ respectively. 
By Riemann-Roch $h^0(N^2K)=2d+g-1,$ and thus $(\pi_{Pic})_*(\Ll_2)$ is a vector bundle.
For $g-1<d<2g-2,$ Riemann-Roch gives $h^0(N^{-2}K^3)=-2d+5g-5,$ so, in this range, $(\pi_{Pic})_*(\Ll_1)$ is also a vector bundle over $Pic^d(\Sigma).$ 
In the range $g-1<d<2g-2$, Bradlow, Garcia-Prada and Gothen \cite{MaximalSP4} proved that $\Mm_J^d(\sSp(4,\R))$ is given by the following fiber bundle 
\[\xymatrix{\left.\left(((\pi_{Pic})_*\Ll_1-\{zero-section\})\oplus(\pi_{Pic})_*\Ll_2\oplus(\pi_{Pic})_*\pi_\Sigma^* K^2\right)\right/\C^*\ar[d]_\pi\\Pic^d(\Sigma)}\]
where $\C^*$ acts by \eqref{fiberwiseactionIntro}. Furthermore, for $N\in\ Pic^d(\Sigma),$ the fiber is given by \[\pi^{-1}(N)\cong \Oo_{\P^{-2d+5g-6}}(1)^{\oplus 2d+g-1}\times\C^{3g-3}.\]
However, for $2g-2\leq d\leq 3g-3,$ $h^0(N^{-2}K^3)$ depends on $N,$ and not just $d.$ For all $d$ we have:
\begin{PropositionIntro}[\ref{ParamMdProposition}]
	There is a surjection from the smooth manifold $\Ff_J^d$ to the support of $(\pi_{Pic})_*(\Pp^{-2}_d\otimes\pi_\Sigma^*)$
	\[\pi:\xymatrix{\Ff_J^d\ar@{->>}[r]&supp((\pi_{Pic})_*(\Pp^{-2}_d\otimes\pi_\Sigma^*))} \subset\ Pic^d(\Sigma).\]
	If $N\in supp(\pi_{Pic})_*(\Pp^{-2}_d\otimes\pi_\Sigma^*))$ with $h^0(N^{-2}K^3)>1$ then the fiber is $\pi^{-1}(N)\cong\Oo_{\P^{a-1}}(1)^{\oplus b},$ where
  	$a=h^0(N^{-2}K^3) $ and  $b=h^0(N^2K)=2d+g-1$. 
  	If $N\in supp(\pi_{Pic})_*(\Pp^{-2}_d\otimes\pi_\Sigma^*))$ with $h^0(N^{-2}K^3)=1$  then $\pi^{-1}(N)\cong\C^{2d+g-1}.$
\end{PropositionIntro}
\smallskip
\textbf{The paper is oraganized as follows: }\ \
We start by recalling the necessary definitions and tools of the theory of Higgs bundles and harmonic maps. In section \ref{sp4Higgs}, we recall the work of \cite{MaximalSP4,CrossRatioAnosoveProperEnergy} on maximal $\sSp(4,\R)$ representations and their corresponding Higgs bundles and harmonic maps. 
With this set up, Theorem \ref{uniqueminThm} is introduced, and an important metric splitting property is proven (Proposition \ref{MetricSplittingProposition}). 
Using the parameterization of $\Mm_J^d$ and assuming Theorem \ref{uniqueminThm}, we prove the mapping class group parameterization of $\Rr_d(\sSp(4,\R))$ for $g-1<d\leq3g-3$ in Theorem \ref{MCGinvparamTheorem}.

The rest of the paper is devoted to proving Theorem \ref{uniqueminThm}. To do this, for any complex semisimple Lie group $\sG,$ we introduce a special class of maps (Definition \ref{GG0cyclicsurfaceDef}) called \em $\sG$-cylic surfaces \em from a Riemann surface to the homogeneous space $\sG/\sT_0,$ where $\sT_0$ is the maximal compact torus of the split real form $\sG_0\subset\sG.$ 
In section \ref{LieTheory} and \ref{Homspaces}, we introduce the necessary Lie theory and homogeneous space geometry. Then in section \ref{cyclicSurf}, we define cyclic surfaces in general and prove certain properties of their deformations. For representations in $\Rr_d(\sSp(4,\R)),$ we show in Proposition \ref{cyclicSurfacesforRd} that, in a suitable conformal structure, the associated Higgs bundles give rise to special equivariant $\sSp(4,\R)$-cyclic surfaces.
Finally, in section \ref{proofofTHM}, we
prove a rigidity result (Theorem \ref{RigidityOfCyclicSurf}) for special equivariant cyclic surfaces, then complete the proof of Theorem \ref{uniqueminThm}.

\vspace{.3cm}

\textbf{Acknowledgements:}\ \ I would like to thank Daniele Alessandrini, Steve Bradlow, and Fran\c{c}ois Labourie for many fruitful discussions. I am very grateful to Marco Spinaci for many enlightening email correspondences and useful comments. Also, I would like to thank Qiongling Li and Andy Sanders for countless stimulating conversations about representation varieties, harmonic maps and Higgs bundles. I acknowledge the support from U.S. National Science Foundation grants DMS 1107452, 1107263, 1107367 ``RNMS: GEometric structures And Representation varieties'' (the GEAR Network). I have benefited greatly from the opportunities the GEAR Network has provided me.

\bigskip
\section{Higgs bundles and harmonic maps}
\smallskip
Let $S$ be a closed orientable surface of genus at least two, fix a conformal structure $J$ and denote the corresponding Riemann surface $(S,J)$ by $\Sigma.$ Let $K\ra\Sigma$ be the canonical bundle, i.e. $K=T^{*,(10)}_\C\Sigma$ is the holomorphic cotangent bundle. As the main results of this paper exploit all four aspects of the nonabelian Hodge correspondence, we introduce the main objects at length.
\subsection{Higgs bundles in general}
Let $\sG$ be a real simple Lie group with Lie algebra $\fg$ and Killing form $B_\fg.$ 
Everything discussed below has an analogous statement for reductive groups, however we will only consider simple groups. Any involution $\sigma:\fg\ra\fg$ gives a decomposition $\fg=\fh\oplus\fm$ into $\pm1$ eigenspaces such that
\[\xymatrix{[\fh,\fh]\subset\fh\ ,&[\fh,\fm]\subset\fm\ ,&[\fm,\fm]\subset\fh}.\]
An involution $\sigma$ is called a \em Cartan involution \em if 
\[B_\sigma(X,Y)=-B_\fg(X,\sigma(Y))\]
is a symmetric positive definite bilinear form. 
For a Cartan involution, it follows that the splitting $\fg=\fh\oplus\fm$ is orthogonal and that $B$ is positive definite on $\fm$ and negative definite on $\fh.$ 
Thus, $\fh$ is the Lie algebra of a maximal compact subgroup $\sH\subset\sG.$ 
Furthermore, the splitting $\fg=\fh\oplus\fm$ is $Ad_H$-invariant, and thus, after complexifying,  the splitting $\fg_\C=\fh_\C\oplus\fm_\C$ is a splitting of $H_\C$-representations. 
When $\fg$ is a complex semisimple Lie algebra, we will denote a Cartan involution by $\theta.$ 
A Cartan involution on a complex semisimple Lie algebra is equivalent to a compact real form; in particular, a Cartan involution $\theta$ is conjugate linear. 

\begin{Definition}
	A \em $\sG$-Higgs bundle \em over a compact Riemann surface $\Sigma$ is a pair $(\Ee,\phi)$ where 
	\begin{itemize}
		\item $\Ee\ra\Sigma$ is a holomorphic principal $H_\C$ bundle,
		\item $\phi\in H^0(\Sigma,(\Ee\times_{H_\C}\fm_\C)\otimes K)$ is a holomorphic section of the associated $K$-twisted $\fm_\C$-bundle. 
	\end{itemize}
\end{Definition}

We will focus on $\sG=\sSL(n,\C)$ and $\sG=\sSp(2n,\R),$ in these cases, we can work with vector bundles with extra structure. For $\sSL(n,\C),$ the maximal compact is $\sSU(n)$ and the Cartan decomposition is 
\[\fsl(n,\C)=\fsu(n)\oplus i \fsu(n).\] 
Complexifying yields $H_\C=\sSL(n,\C)$ and $\fm_\C=\fsl(n,\C).$ Furthermore, the representation of $H_\C$ on $\fm_\C$ is the adjoint representation of $\sSL(n,\C)$ on $\fsl(n,\C).$ 
Using the standard representation of $\sSL(n,\C)$ on $\C^n,$ a $\sSL(n,\C)$-Higgs bundle is given by a pair $(E,\phi)$ where 
\begin{itemize}
	\item $E\ra\Sigma$ is a rank $n$ holomorphic vector bundle with $det(E)=\Oo$
	\item$\phi\in H^0(\Sigma, End(E)\otimes K)$ is a traceless holomorphic $K$-twisted endomorphism
\end{itemize} 
This is the original definition of Hitchin \cite{selfduality}. 

For $\sSp(2n,\R),$ the maximal compact is $\sU(n),$ and if $Sym^2(V)$ is the symmetric tensor product of the standard representation of $\sGL(n,\C)$ on $\C^n,$ then complexifying the Cartan decomposition of $\fsp(2n,\R)$ gives the decomposition 
\[\fsp(2n,\C)=\fgl(n,\C)\oplus \left(Sym^2(V)\oplus Sym^2(V^*)\right)\]
(see subsection \ref{sp4subsection} for more details).
Thus, an $\sSp(2n,\R)$-Higgs bundle is given by a triple $(V,\beta,\gamma)$ where: 
\begin{itemize}
	\item $V\ra \Sigma$ is a rank $n$ holomorphic vector bundle
	\item $\beta\in H^0(Sym^2(V)\otimes K)$ is a holomorphic symmetric section of $V\otimes V$ twisted by $K$
	\item $\gamma\in H^0(Sym^2(V^*)\otimes K)$ is a holomorphic symmetric section of $V^*\otimes V^*$ twisted by $K.$
\end{itemize}
Given an $\sSp(2n,\R)$-Higgs bundle $(V,\beta,\gamma)$, the associated $\sSL(2n,\C)$-Higgs bundle is $(E,\phi)=\left(V\oplus V^*,\mtrx{0&\beta\\\gamma&0}\right).$ 

To form a moduli space of Higgs bundles, we need a notion of stability.
\begin{Definition}
An $\sSL(n,\C)$-Higgs bundle $(E,\phi)$ is \em semistable\em, if for all proper holomorphic subbundles $F\subset E$ we have $deg(F)\leq 0.$ If the inequality is always strict, then $(E,\phi)$ is called \em stable, \em and if $(E,\phi)$ is a direct sum of stable $\sSL(n_i,\C)$-Higgs bundles, then it is called \em polystable\em.
\end{Definition}

\begin{Remark}\label{SLstabenough}
For an arbitrary real reductive Lie group $\sG,$ the definition of stability and semistability is significantly more involved, see \cite{HiggsPairsSTABILITY,HiggsbundlesSP2nR}. However, if $\sG\subset\sSL(n,\C)$ then a $\sG$-Higgs bundle is \em unstable \em if and only if the corresponding $\sSL(n,\C)$-Higgs bundle is unstable. 
The $\sSp(4,\R)$-Higgs bundles $(V,\beta,\gamma)$ studied in this paper (see Remark \ref{StabilityofMdRemark}) are stable or polystable as $\sSp(4,\R)$-Higgs bundles if and only if the associated the $\sSL(4,\C)$-Higgs bundle $\left(V\oplus V^*,\mtrx{0&\beta\\\gamma&0}\right)$ is stable or polystable, respectively. 
\end{Remark}

The \em moduli space of $\sG$-Higgs bundles \em over $\Sigma=(S,J)$ is defined as
\begin{equation}
	\Mm_J(\sG)=\{(\Ee,\phi)\ |\ (\Ee,\phi)\ \text{a\ polystable\ } \sG \text{-Higgs\ bundle}\}/\Gg_{\sH_\C}
\end{equation}
where $\Gg_{\sH_\C}$ is $\sH_\C$-gauge group (i.e. the group of smooth bundle automorphism of $\Ee$). 

\begin{Definition}
 	A $\sG$-Higgs bundle $(\Ee,\phi)$ is called \em simple \em if $Aut(\Ee,\phi)\subset \Gg_{\sH_\C}=Z(\sH_\C).$ For $\sG=\sSL(n,\C),$ stable Higgs bundles are automatically simple, this is not true in general.
 \end{Definition} 
 
To go from the moduli space of $\sG$-Higgs bundles to the space of representations of $\pi_1(\Sigma)$ we have the following theorem, it was originally proved by Hitchin \cite{selfduality} for $\sSL(2,\C)$ Simpson \cite{localsystems} for general complex reductive groups and Bradlow, Garcia-Prada, Mundet i Riera \cite{RelKHpairscorresp} and Garcia-Prada, Gothen, Mundet i Riera \cite{HiggsPairsSTABILITY} for real reductive groups. 

\begin{Theorem}\label{GKHcorr} Let $\sG$ be a real reductive Lie group and $\theta$ a fixed Cartan involution on $\fg_\C.$ If $(\Ee,\phi)$ is a stable and simple $\sG$-Higgs bundle then there exists a unique reduction of structure group $\sigma:\xymatrix{\Sigma\ra \Ee/H}$
	with associated `Chern connection' $A_\sigma$ solving the Higgs bundle equations
	\begin{eqnarray*}
		F_{A_\sigma}+[\phi,-\theta(\phi)]=0\\
	\nabla_{A_\sigma}^{01}\phi=0 
	\end{eqnarray*}
Conversely, if $(A_\sigma,\phi)$ is such a solution, then the corresponding Higgs bundle is polystable Higgs bundle.
\end{Theorem}
For the appropriate notion of the Chern connection of a reduction of structure of a holomorphic principal bundle to its maximal compact see the appendix of \cite{UniversalKHcorr}. 
A solution $(A,\phi)$ to the Higgs bundle equations gives rise to the flat $\sG$-connection $A+\phi-\theta(\phi).$ This gives a map to the representation variety
\[\Mm_{J}(\sG)\lra \Rr(\pi_1(S),\sG).\]

\begin{Remark}\label{simpleRepsDef}
	For reductive Lie group $\sG,$ a representation $\rho$ is called \em irreducible \em if the centralizer of $\rho$ is the center of $\sG$. Given a stable and simple $\sG$-Higgs bundle $(\Ee,\phi)$, the corresponding representation $\rho$ is irreducible (see Theorem 3.32 of \cite{HiggsPairsSTABILITY}). Irreducible representations and stable and simple Higgs bundles define the smooth points of $\Rr(\sG)$ and $\Mm_J(\sG).$
\end{Remark}

 In the case of $\sSL(n,\C),$ Theorem \ref{GKHcorr} simplifies to the following. 
\begin{Theorem}\label{SLKHcorr} Let $(E,\phi)$ be a stable $\sSL(n,\C)$-Higgs bundle, then there exists a unique hermitian metric $h,$ with Chern connection $A_h,$ solving the Hitchin equations
\begin{eqnarray}\label{Hit}F_{A_h}+[\phi,\phi^{*_h}]=0\\ 
\nabla_{A_h}^{01}\phi=0\nonumber
\end{eqnarray}
where $F_{A_h}$ is the curvature of $A_h$ and $\phi^{*_h}$ is the hermitian adjoint. Conversely,  if $(A_h,\phi)$ is a solution then the corresponding Higgs bundle is polystable. 
\end{Theorem}

\begin{Remark}\label{MetricRemark}
Denote the holomorphic structure on $E$ by $\bar \p_E.$ Given a stable Higgs bundle $(\bar \p_E,\phi),$ let $h$ be the unique metric which solves the Hitchin equations (\ref{Hit}). 
For any $\sSL(n,\C)$-gauge transformation $g,$ the pair $(g^{-1}\bar \p_E g,g^{-1}\phi g)$ also has a unique metric $h'$ solving (\ref{Hit}). The metrics $h$ and $h'$ are related by $h'=hg^{*_h}g.$ 
This follows from general gauge theoretic arguments, for example see \cite{VorticesInLineBundles}.
\end{Remark}

\begin{Remark}\label{U1Remark}
	Note that if $h$ is a solution metric for a Higgs bundle $(E,\phi),$ then for all $\lambda\in \sU(1),$ $h$ is also the solution metric for $(E,\lambda\phi).$ This gives a $\sU(1)$ action on the moduli space of Higgs bundles.
	This action and its restriction to the $k^{th}$ roots of unity $\langle \zeta\dwn{k}\rangle$ play a key role in the proof of Proposition \ref{MetricSplittingProposition}.
\end{Remark}

\begin{Example}\label{HitEx}
	Fix a square root $K^\haf$ of the canonical bundle $K,$ and let $q_2\in H^0(K^2)$ and $q_4\in H^0(K^4)$ be holomorphic differentials. Consider the $\sSp(4,\R)$-Higgs bundle
	\[(V,\beta,\gamma)=\left(K^\frac{3}{2}\oplus K^{-\haf}, \mtrx{q_4&q_2\\q_2&1},\mtrx{0&1\\1&0}\right).\]
	The associate $\sSL(4,\C)$-Higgs bundle is 
	\[\left(V\oplus V^*,\mtrx{0&\beta\\\gamma&0}\right)=\left( K^\frac{3}{2}\oplus K^{-\haf}\oplus K^{-\frac{3}{2}}\oplus K^\haf,\mtrx{&&q_4&q_2\\&&q_2&1\\0&1&&\\1&0&&}\right).\]
	This is the Higgs bundle description of the $\sSp(4,\R)$-Hitchin components \cite{liegroupsteichmuller}. 
\end{Example}	

\begin{Definition} Let $\sG$ be a complex simple Lie group of rank $\ell,$ and fix a homogeneous basis $(p_1,\cdots, p_{\ell})$ of the invariant polynomials $\C[\fg]^{\sG}$ with $deg(p_j)=m_j+1$. Let $(E,\phi)$ be a $\sG$-Higgs bundle, then $p_j(\phi)\in H^0(\Sigma, K^{m_j+1})$ and we have a map 
	\[\xymatrix@R=0em{
	\Mm_J(\sSL(n,\C))\ar[rr]^{(p_1,\cdots,p_{\ell})\ \ \ }&&\bigoplus\limits_{j=1}^n H^0(\Sigma,K^{m_j+1})\\
	[\Ee,\phi]\ar@{|->}[rr]&&(p_1(\phi),\cdots,p_{\ell}(\phi))
	}\]
	called the \em Hitchin Fibration.\em
\end{Definition}
For $\sSL(n,\C),$ $m_j=1, 2, \cdots, n-1$ and for $\sSp(2n,\C)$ $m_j=1,3,\cdots, 2n-1.$ Up to a constant, the holomorphic quadratic differential $p_1(\phi)$ is given by $Tr(\phi^2).$ 

\subsection{Harmonic maps and Corlette's theorem}
To go from the representation variety to the moduli space of Higgs bundles we need to develop some harmonic map theory. Let $(N_1,g_1)$ and $(N_2,{g_2})$ be two Riemannian manifolds with $N_1$ compact and orientable. A smooth map $f:N_1\ra N_2$ is harmonic if it is a critical point of the energy functional 
\[\Ee(f)=\haf\int_{N_1}|df|^2dvol_{g_1}.\]
Here $df$ is a section of $T^*N_1\otimes f^*TN_2$ and the norm is taken with the metric induced by $g_1$ and ${g_2}.$ 
If $(D^{f^*\nabla_{{g_2}}})^*$ is the adjoint (with respect to the metric $g_1$) of the covariant derivative induced by the Levi-Civitta connection for ${g_2},$ then the Euler-Lagrange equations for harmonic maps can be written as
\[(D^{\nabla_{f^*{g_2}}})^*(df)=0.\] 
When the domain $N_1$ is a surface, harmonic maps have many additional properties. For instance, the energy, and hence the Euler-Lagrange equations, are invariant under a conformal change of the domain metric. 
Thus, harmonicity of maps $f:\Sigma\ra(N,g)$ with domain a Riemann surface $(S,J)=\Sigma$ makes sense. 
Given a map $f:\Sigma\ra (N,{g})$, denote the decomposition of the exterior derivative $d$ and the covariant derivative $D^{f^*\nabla_h}:\Omega^1(\Sigma,f^*TN)\ra \Omega^2(\Sigma,f^*TN)$ into their $(1,0)$ and $(0,1)$ parts by
\[\xymatrix{d=d^{10}+d^{01}&\txt{ and }& D^{f^*\nabla_{g}}=(D^{f^*\nabla_{g}})^{10}+(D^{f^*\nabla_{g}})^{01}}. \]
For a Riemann surface, a map $f:\Sigma\ra (N,{g})$ is harmonic if and only if $df^{10}$ is holomorphic, that is 
\[(D^{\nabla_{g}})^{01}(df^{10})=0.\]

Any smooth map $f:\Sigma\ra (N,{g})$ gives a symmetric 2-tensor $f^*{g}\in\Omega^0(\Sigma,S^2(T^*\Sigma))$; the $(2,0)$ part $q_f=(f^*{g})^{20}\in \Omega^0(\Sigma,K^2)$ is a quadratic differential. 
\begin{Definition}
	If $f:\Sigma\ra (N,{g})$ is harmonic then $q_f=(f^*{g})^{20}\in H^0(\Sigma,K^2)$ is a holomorphic quadratic differential called the \em Hopf differential \em of $f$. 
\end{Definition}
If $z$ is a local complex coordinate, then the Hopf differential of a harmonic map $f$ is defined by 
\[q_f\left(\frac{\p}{\p z},\frac{\p}{\p z}\right)=g\left(df\left(\frac{\p}{\p z}\right),df\left(\frac{\p}{\p z}\right)\right)
.\]
Writing $q_f$ in terms of $z=x+iy,$ it is clear that if $q_f=0$ then the rank of $df$ is always $0$ or $2$. 
\begin{Remark}\label{BranchedMinSurfRem}
	A harmonic map $f$ is \em weakly conformal \em if and only if $q_f=0.$ Equivalently \cite{MinImmofRiemannSurf,SchoenYauMinimalSurfEnergy}, $q_f=0$ if and only if $f$ is a \em branched minimal immersion\em. If $f$ is a branched minimal immersion with $df$ nowhere zero, then $f$ is a \em minimal immersion.\em
\end{Remark}

Given a representation $\rho:\pi_1(S)\ra \sG$, consider the associated flat principal $\sG$-bundle $\widetilde S\times_\rho \sG.$ A metric on a principal $G$-bundle $E$ is a reduction of structure group to the maximal compact subgroup which can be interpreted as a section of the fiber bundle
$h:S\ra E/H.$
Since the bundle $\widetilde S\times_\rho G$ is flat, such a section $h$ is equivalent to an $\rho$-equivariant map to the symmetric space
\[h_\rho:\widetilde S\longrightarrow G/H.\]
To obtain a Higgs bundle from representation $\rho\in\Rr(\sG),$ to each choice of conformal structure $\Sigma=(S,J),$ Corlette's theorem \cite{canonicalmetrics} provides a special metric on the flat bundle $\widetilde \Sigma\times_\rho G.$
\begin{Theorem}\label{CorletteThm}
	Let $\rho\in\Rr(\sG)$ be an irreducible representation, and fix a conformal structure $(S,J)=\Sigma.$ Then there exists a unique metric 
	\[h_\rho:\widetilde\Sigma\longrightarrow G/H\]
	which is harmonic.
\end{Theorem}

Given $\rho:\pi_1(S)\ra\sG,$ denote the flat connection on $\widetilde\Sigma\times_\rho\sG$ by $D\in\Omega^1(\widetilde\Sigma\times_\rho\sG,\fg).$ 
Any metric $h_\rho$ gives a $\sH$-subbundle $E_\sH\subset \widetilde\Sigma\times_\rho\sG,$ and allows us to decompose the flat connection as $D=A+\psi$ where
 \[\xymatrix{A\in\Omega^1(E_\sH,\fh)&\text{and}& \psi\in\Omega^1(\Sigma,E_\sH\times_\sH\fm)}.\]
 Here $A$ is a connection $1$-form on $E_\sH$ and $\psi$ is identified with the derivative of the metric $dh_\rho.$
If the metric $h_\rho$ is harmonic, then $\nabla_A^*\psi=0,$ or equivalently $\nabla^{01}{\psi^{10}}=0.$ 
 With respect to this decomposition, the flatness for $D$ and harmonic equation for $h_\rho$ decompose as 
\begin{equation}\label{HarmonicEq}
	\xymatrix{F_A+\haf[\psi,\psi]=0\ ,&&\nabla_A\psi=0\ ,&&  \nabla_A^{01}\psi^{10}=0 }
\end{equation}
The Higgs bundle corresponding to $h_\rho$ is $(\Ee,\phi)=(E_H\times_{\sH_\C}\sH_\C,\psi^{10}).$ Since $\nabla_A\psi=\nabla_A^*\psi=0,$ we have $\psi^{01}=-\theta(\psi^{10}).$ 
Thus, the pair $(A,\psi^{10})$ solves the Higgs bundle equations 
\[\xymatrix{F_A+[\phi,-\theta(\phi)]=0&\text{and}&\nabla_A^{01}\psi^{10}=0}.\]
\begin{Remark}\label{HiggsMinimalSurfaceRem}
Under the identification of the representation variety $\Rr(\sG)$ with the Higgs bundle moduli space $\Mm_J(\sG),$ the $(1,0)$-part of the differential of the harmonic map from Corelette's Theorem is identified with Higgs field $\phi$ (see Remark \ref{Higgsfieldanddh}). 
For $\sG$ simple, the metric on $\sG/\sH$ is induced by the Killing form $B_\fg.$ 
 Hence, for $\rho\in\Rr(\pi_1(S),\sG),$ up to a constant, the Hopf differential of the corresponding harmonic map $h_\rho$ is given by 
\[q_{h_\rho}=Tr(dh_\rho^{10}\otimes dh_\rho^{10})=Tr(\phi^2).\]
Thus, the harmonic map $h_\rho$ is a branched minimal immersion if and only if the quadratic differential associated the corresponding Higgs bundle, via the Hitchin fibration, vanishes. This will be crucial in linking the Higgs bundle description of \cite{MaximalSP4} with the minimal surface theory of \cite{CrossRatioAnosoveProperEnergy}.
\end{Remark}

\bigskip
\section{$\sSp(4,\R)$-Higgs bundles and harmonic maps}\label{sp4Higgs}\smallskip
For $\sG=\sSp(4,\R),$ the complexification of the maximal compact subgroup is $\sH_\C=\sGL(2,\C).$ 
For a $\sSp(4,\R)$ Higgs bundle $(V,\beta,\gamma)$, $\tau=deg(V)\in\Z$ defines an integer invariant called the Toledo invariant.
 Given two $\sSp(4,\R)$-Higgs bundles $(V,\beta,\gamma)$ and $\ (V',\beta',\gamma'),$ if $deg(V)\neq deg(V')$ then $(V,\beta,\gamma)$ and $(V',\beta',\gamma')$ are in different connected components.
  This gives a decomposition
\[\Mm_J(\sSp(4,\R))=\bigsqcup\limits_{\tau\in\Z}\Mm_{J,\tau}(\sSp(4,\R)).\]
The map sending $(V,\beta,\gamma)$ to $(V^*,\gamma,\beta)$ gives an isomorphism $\Mm_{J,\tau}(\sSp(4,\R))\cong \Mm_{J,-\tau}(\sSp(4,\R)).$
Note that if $deg(V)>2g-2,$ then $deg(V^*\otimes V^*\otimes K)<0.$ Thus, $\gamma\in H^0(Sym^2(V^*)\otimes K)$ implies $\gamma=0.$ 
In this case, $V\subset V\oplus V^*$ is a $\phi$-invariant subbundle with positive degree, and hence $(V,\beta,\gamma)$ is unstable by Remark \ref{StabilityofMdRemark}. 
So, $|\tau|\leq 2g-2$ and 
\[\Mm_J(\sSp(4,\R))=\bigsqcup\limits_{-2g+2\leq\tau\leq2g-2}\Mm_{J,\tau}(\sSp(4,\R)).\] 
Gothen \cite{sp4GothenConnComp} showed that, for $\tau=0,$ $\Mm_{J,\tau}(\sSp(4,\R))$ is connected, and, for $|\tau|=2g-2,$ the moduli space $\Mm_{J,\tau}(\sSp(4,\R))$ has $3^{2g}+2g-4$ connected components. 
In \cite{nonmaxSp4}, it is shown that $\Mm_{J,\tau}(\sSp(4,\R))$ is connected for all other values of the Toledo invariant. This gives $1+2(2g-1)+2(3^{2g}+2g-4)$ total connected components for $\Mm_J(\sSp(4,\R)).$

\subsection{Maximal components for $\sSp(4,\R)$}
Higgs bundles $(V,\beta,\gamma)$ with $|deg(V)|=|\tau|=2g-2$ are called maximal. 
When $\tau=2g-2,$ polystablility forces the holomorphic map $\gamma:V\ra V^*\otimes K$ to be an isomorphism \cite{sp4GothenConnComp}. Using this fact, to a maximal $\sSp(4,\R)$ Higgs bundle $(V,\beta,\gamma)$ one associates a $\sGL(2,\R)$ $K^2$-twisted Higgs bundle $(W,\varphi)$ (i.e. a $\sGL(2,\R)$ Higgs bundle where the Higgs field is twisted by $K^2$ instead of $K$), called its \em Cayley partner\em .
 The Cartan decomposition of $\fgl(2,\R)$ is $\fo(2,\R)\oplus sym(\R^2),$ and, complexifying, we have
\[\fgl(2,\C)=\fo(2,\C)\oplus sym(\C^2).\]
Thus, a $K^2$-twisted $\sGL(2,\R)$-Higgs bundle is a triple $(W,Q_W,\varphi)$ where $(W,Q_W)$ is a $\sO(2,\C)$ bundle and $\varphi\in H^0(End(W)\otimes K^2)$ satisfying $\varphi^TQ_W=Q_W\varphi.$

The characteristic classes of the Cayley partner help to distinguish connected components of $\Mm^{2g-2}_J(\sSp(4,\R))$. 
We will recall how this works for $\sSp(4,\R)$ \cite{sp4GothenConnComp,MaximalSP4}, for a general development of the theory of Cayley partners see \cite{CayleyPartnerThesis}.
Fix a square root of the canonical bundle $K^\haf$ and set $W=V^*\otimes K^\haf.$ 
Using the fact that $\gamma:V\ra V^*\otimes K$ is an isomorphism, define an orthogonal structure $Q_W:W^*\ra W$ by
\[Q_W=\gamma\otimes Id_{K^{-\haf}}:V\otimes K^{-\haf}\longrightarrow V^*\otimes K\otimes K^{-\haf}. \]
For the Cayley partner, the Higgs field $\varphi:W\ra W\otimes K^2$ is given by $\varphi=(\gamma\otimes Id_{K\otimes K^\haf})\circ(\beta\otimes K^\haf),$ i.e.
\[W=\xymatrix@=4em{V^*\otimes K^{\haf}\ar[r]^{\beta\otimes Id\ \ \ \ }&V\otimes K\otimes K^{\haf}\ar[r]^{\gamma\otimes Id\ \ \ \ }& V^*\otimes K\otimes K\otimes K^\haf}=W\otimes K^2.\]
The map $\varphi$ is $Q_W$-symmetric, thus $(W,Q_W,\varphi)$ defines a $K^2$-twisted $\sGL(2,\R)$-Higgs bundle. 

The $\sO(2,\C)$ bundle $(W,Q_W)$ has a first and second Stieffel-Whitney class
\[\xymatrix{w_1(W,Q_W)\in H^1(\Sigma,\Z/2\Z)=\left(\Z/2\Z\right)^{2g}&\text{and}&w_2(W,Q_W)\in H^2(\Sigma,\Z/2\Z)=\Z/2\Z}.\]
There are $2\cdot 2^{2g}-2$ possible values for $(w_1(W,Q_W),w_2(W,Q_W))$ with $w_1(W,Q_W)\neq 0.$ When $w_1(W,Q_W)=0$, the structure group of the $\sO(2,\C)$-bundle lifts to $\sSO(2,\C),$ 
in this case, we have a \em Chern class, \em and \em Proposition 3.20 \em of \cite{MaximalSP4}.
 \begin{Proposition}\label{BGGMdProp}
 	Let $(V,\beta,\gamma)$ be a maximal $\sSp(4,\R)$-Higgs bundle with Cayley partner $(W,Q_W)$ and $w_1(W,Q_W)=0,$ then there is a line bundle $N\ra\Sigma$ so that $V=N\oplus N^{-1}K.$ With respect to this decomposition, 
 	\[\beta=\mtrx{\nu&q_2\\q_2&\mu}:N^{-1}\oplus NK^{-1}\ra NK\oplus N^{-1}K^2\ \ \ \ \text{and}\ \ \ \ \gamma=\mtrx{0&1\\1&0}:N\oplus N^{-1}K\ra N^{-1}K\oplus N.\]
The line bundle $N$ satisfies a degree bound, $g-1\leq deg(N)\leq 3g-3;$ for $g-1<deg(N),$ the line bundle $N$ is \em unique \em and when $deg(N)=g-1,$ the line bundle $N$ is unique up to a multiple of a square root of $\Oo.$ 
 \end{Proposition}
 The proof of this proposition makes extensive use of Mumford's classification of rank $2$ holomorphic orthogonal bundles \cite{MumO2Bun}. The degree of $N$ provides $2g-1$ extra invariants; set $d=deg(N),$ and denote the corresponding moduli space by $\Mm_J^d(\sSp(4,\R)).$ For $deg(N)=3g-3,$ stability forces $N^2=K^3,$ and there are at least $2^{2g}$ connected components corresponding to choices of square roots of $K$. Thus, there are
 \[2\cdot2^{2g}-2+2g-2+2^g\ =\ 3\cdot 2^{2g}+2g-4\]
 invariants for $\sSp(4,\R)$-Higgs bundles with $deg(V)=2g-2,$ and 
 we have 
\begin{equation}\label{componentsEq}
	\Mm^{2g-2}_J(\sSp(4,\R))=\bigsqcup\limits_{w_1\neq 0}\Mm^{w_1,w_2}_J(\sSp(4,\R))\bigsqcup\limits_{g-1\leq d< 3g-3}\Mm_J^d(\sSp(4,\R))\bigsqcup\limits_{L^2=K} \Mm_J^L(\sSp(4,\R)).
\end{equation}
 In \cite{sp4GothenConnComp}, it is shown that each of the above moduli space is nonempty and connected. For $deg(N)=3g-3,$ the $2^{2g}$ connected components are the Hitchin components of Example \ref{HitEx}. When $g-1\leq d<3g-3,$ we call the components $\Mm^d_J(\sSp(4,\R))$ the \em Gothen components\em. 

\begin{Remark}\label{StabilityofMdRemark}
We will restrict ourselves to describing the Higgs bundles in the Hitchin and Gothen components. By Proposition \ref{BGGMdProp}, the $\sSL(4,\C)$ Higgs bundle associated to a Higgs bundle $(V,\beta,\gamma)\in \Mm_J^d(\sSp(4,\R))$ is of the form 
\begin{equation}\label{Sl4d}
	(E,\phi)=\left(N\oplus N^{-1}K\oplus N^{-1} \oplus NK^{1},\mtrx{&&\nu&q_2\\&&q_2&\mu\\0&1&&\\1&0&&}\right)
\end{equation}
If $\mu=0,$ then $N\oplus NK^{-1}$ is an invariant subbundle of $E$ with of degree $2d-2g+2$. 
Thus, for $g-1<d$ stability forces $\mu\neq0.$ 
Furthermore, by Proposition 3.24 of \cite{MaximalSP4}, for $g-1<d$, all isomorphism classes in $\Mm_J^d(\sSp(4,\R))$ are stable and simple. 
When $d=g-1,$ the Higgs bundle is stable if and only if $\mu\neq0.$ 
By Remark \ref{simpleRepsDef}, it follows that $\Mm_J^d(\sSp(4,\R))$, and hence $\Rr_d(\sSp(4,\R)),$ is \em smooth \em if and only if $g-1<d\leq 3g-3.$

\end{Remark}

Let $\Rr_d(\sSp(4,\R))$ be the component the representation variety which corresponds to $\Mm_J^d(\sSp(4,\R)).$
Using the description of the possible Zariski closures of maximal $\sSp(4,\R)$ representations of \cite{BIWmaximalToledoAnnals}, Bradlow, Garcia-Prada, and Gothen showed \cite{MaximalSP4}, if $g-1<d<3g-3$ and $\rho\in\Rr_d(\sSp(4,\R)),$ then $\rho$ is Zariski dense. 
Furthermore, by Remark \ref{StabilityofMdRemark}, the Gothen components $\Rr_d(\sSp(4,\R))$ for $g-1<d<3g-3$ are smooth.

\subsection{Metric splitting and minimal immersions}
Recall that $g-1\leq deg(N)\leq 3g-3,$ the $\sSL(4,\C)$-Higgs bundle \eqref{Sl4d} is determined by the quadruple $(N,\mu,\nu,q_2).$
 The bundle $E=V\oplus V^*$ has symplectic structure $\Omega=\mtrx{0&I\\-I&0},$ with respect to which $\phi^T\Omega+\Omega\phi=0.$ Thus, \eqref{Sl4d} is in fact an $\sSp(4,\C)$-Higgs bundle.
For the group $\sSp(4,\C),$ the polynomial ring $\C[\fsp(4,\C)]^{\sSp(4,\C)}$ has two homogeneous generators $(p_1,p_2)$ of degree two and four. One choice of generators is 
\[\xymatrix{p_1=Tr(X^2)&\text{and}& p_2=Tr(X^4)}.\]
For any other basis $(p_1',p_2'),$ there are constants $A,B,C$ so that 
\[\xymatrix{p_1'=Ap_1&\text{and}& p_2'=Bp_1^2+Cp_2}.\]
Thus, for any choice of basis of the invariant polynomials, the holomorphic quadratic and quartic differentials associate the $(E,\phi)$ via the Hitchin fibration are
\[\xymatrix{ATr(\phi^2)=4Aq_2&\text{and} &BTr(\phi^2)+CTr(\phi^4)=16Bq_2\otimes q_2+C4\mu\otimes\nu}.\]
\begin{Lemma}\label{HigBunForMinSurfLemma}
	Let $\rho\in\Rr_d(\sSp(4,\R))$ and fix a conformal structure $\Sigma=(S,J).$ If the harmonic $\rho$-equivariant map $h_\rho$ is a branched minimal immersion, then the corresponding Higgs bundle $(V,\beta,\gamma)$ is given by 
	\begin{equation}\label{minHiggsbundled}
		\left(N\oplus N^{-1}K,\mtrx{\nu&0\\0&\mu},\mtrx{0&1\\1&0}\right).
	\end{equation}
	Furthermore, up to a constant, the associated holomoprhic quartic differential in the Hitchin base is given by $q_4=\mu\otimes\nu.$
\end{Lemma}

\begin{proof}
Let $\rho\in\Rr_d(\sSp(4,\R))$ and fix a conformal structure $(S,J)=\Sigma.$ By Proposition \ref{BGGMdProp}, the $\sSp(4,\R)$-Higgs bundle corresponding to $\rho$ is given by
\begin{equation*}
	\left(N\oplus N^{-1}K,\mtrx{\nu&q_2\\q_2&\mu},\mtrx{0&1\\1&0}\right).
\end{equation*}
By Remark \ref{HiggsMinimalSurfaceRem}, $h_\rho$ is a branched minimal immersion if and only if 
\[Tr(\phi^2)=Tr\left(\mtrx{0&\beta\\\gamma&0}^2\right)=4q_2=0.\]
In this case, any choice of basis for $\C[\fsp(4,\C)]^{\sSp(4,\C)}$ gives $q_4=p_2(\phi)=CTr(\phi^4)=4C\mu\otimes\nu.$
\end{proof}

Fix a representation $\rho\in\Rr(\sSp(4,\R))$ (or more generally $\Rr(\sG)$) and for each conformal structure denote the corresponding harmonic metric by $h_\rho.$ Consider the following function $\Ee_\rho$ on the Teichm\"uller space $\Tt(S)$
\begin{equation}\label{TeichEnergy}
	\Ee_\rho(J)=\Ee_J(h_\rho)=\haf\int\limits_{S}|dh_\rho|^2 dvol\ :\ \Tt(S)\lra\R
\end{equation}
\begin{Remark}\label{MinimumMinimalSurf}
By \cite{MinImmofRiemannSurf,SchoenYauMinimalSurfEnergy}, critical points of $\Ee_J(h_\rho)$ are a branched minimal immersions, or equivalently, weakly conformal maps. Note that the harmonic map $h_\rho,$ the norm $|dh_\rho|^2$ and the volume element all depend of $J$. 
\end{Remark}
In \cite{CrossRatioAnosoveProperEnergy}, Labourie proved the following theorem:
\begin{Theorem}\label{ProperThm}(\cite{CrossRatioAnosoveProperEnergy}) If $\rho$ is a maximal $\sSp(2n,\R)$ or Hitchin representation, then the energy function $\Ee_\rho:\Tt(S)\ra\R$ is smooth and proper. 
\end{Theorem}
Since $\Ee_\rho$ is proper and bounded below by zero, $\Ee_\rho$ attains a minimum. 
Applying this to our situations, the following corollary follows from remark \ref{MinimumMinimalSurf}.
\begin{Corollary}\label{existenceCor}
	Let $\rho$ be a maximal $\sSp(4,\R)$ representation, then there exists a conformal structure $J$ so that the corresponding $\rho$-equivariant harmonic map is a branched minimal immersion.
\end{Corollary}
\begin{Lemma}\label{RdMinimalImmersion}
	Let $\rho\in\Rr_d(\sSp(4,\R))$ and choose a conformal structure $J$ so that the corresponding $\rho$-equivariant harmonic map $h_\rho$ is a branched minimal immersion, then $h_\rho$ is a \em minimal immersion\em. 
\end{Lemma}

 \begin{proof}
 By Lemma \ref{HigBunForMinSurfLemma}, in the conformal structure $J,$ the $\sSp(4,\R)$-Higgs bundle $(V,\beta,\gamma)$ associated to $\rho$ is given by \eqref{minHiggsbundled}. 
By Remark \ref{HiggsMinimalSurfaceRem} Higgs field is $\phi=\mtrx{0&\beta\\\gamma&0}$ represents the $(1,0)$ part of $dh_\rho.$ Since $\phi$ is nowhere vanishing, by Remark \ref{BranchedMinSurfRem}, the branched minimal immersion $h_\rho$ is branch point free.
 \end{proof}

Our main theorem is the uniqueness of a conformal structure in which the equivariant harmonic map is a minimal surface.

\begin{Theorem}\label{uniqueminThm}
If $S$ be a closed surface of genus at least 2 and $\rho\in\Rr_d(\sSp(4,\R))$ for $g-1<d\leq 3g-3,$ then there exists a unique conformal structure $(S,J_\rho)=\Sigma$ so that the $\rho$-equivariant harmonic map 
 	\[h_\rho:\widetilde\Sigma\ra\sSp(4,\R)/\sU(2)\] is a minimal immersion. 
\end{Theorem}
 
The existence part of Theorem \ref{uniqueminThm} is covered by Corollary \ref{existenceCor} of Labourie's theorem and Lemma \ref{RdMinimalImmersion}. To prove uniqueness, we study special properties of such minimal immersions by studying the Higgs bundles of the form \eqref{minHiggsbundled}. 
The main fact which allows such an analysis is that, for the Higgs bundles of the form \eqref{minHiggsbundled}, the metric solving the Hitchin equations \eqref{Hit} has extra symmetries.
\begin{Proposition}\label{MetricSplittingProposition}
	Let $g-1\leq d\leq 3g-3,$ and $(V,\beta,\gamma)$ be a maximal $\sSp(4,\R)$-Higgs bundle in $\Mm_J^d.$ If  
	\[(V,\beta,\gamma)=\left(N\oplus N^{-1}K,\mtrx{\nu&0\\0&\mu},\mtrx{0&1\\1&0}\right),\]
	has $\mu\neq0,$ then the metric $h$ on $V\oplus V^*$ which solves the Hitchin equations splits as a direct sum \[h=H_1\oplus H_2\oplus H_1^{-1}\oplus H_2^{-1}\] on 
	$V\oplus V^*=N\oplus N^{-1}K\oplus N^{-1}\oplus NK^{-1}.$
\end{Proposition}
\begin{proof}
The metric solving the Hitchin equations for an $\sSp(4,\R)$-Higgs bundle is a reduction of structure of a $\sGL(2,\C)$ bundle to its maximal compact. 
Thus, given a stable and simple $\sSp(4,\R)$ Higgs bundle $(V,\beta,\gamma),$ the unique metric solving the Hitchin equations is a metric $H$ on $V.$ 
On the corresponding $\sSL(4,\C)$ Higgs bundle $(E,\phi)=\left(V\oplus V^*,\mtrx{0&\beta\\\gamma&0}\right),$ the metric solving the Hitchin equations is $h=H\oplus H^{-1}.$

For $\lambda\in\C^*,$ consider the holomoprhic $\sGL(2,\C)$-gauge transformation $g_\lambda$ defined by 
\[g_\lambda=\mtrx{\lambda^{-3}&\\&\lambda^1}:N\oplus N^{-1}K\lra N\oplus N^{-1}K.
 \]
 It acts on the Higgs field $\phi=(\beta,\gamma)=\left(\mtrx{\nu&0\\0&\mu},\mtrx{0&1\\1&0}\right)$ by 
 \[g_\lambda\cdot \left(\mtrx{\nu&0\\0&\mu},\mtrx{0&1\\1&0}\right)=\left(\mtrx{\lambda^{-6}\nu&0\\0&\lambda^2\mu},\mtrx{0&\lambda^21\\\lambda^21&0}\right).\]
When $\lambda=\zeta\dwn{8}=e^\frac{2\pi i}{8},$ the gauge transformation $g_{\zeta\dwn{8}}$ multiplies the Higgs field $(\beta,\gamma)$ by $\zeta_4=e^\frac{2\pi i}{4}.$ 
Thus, the $\sSp(4,\R)$-Higgs bundle $\left(N\oplus N^{-1}K, \mtrx{\nu&0\\0&\mu},\mtrx{0&1\\1&0}\right)$ is a fixed point of the fourth roots of unity $\gen{\zeta\dwn{4}}\subset\sU(1).$ Furthermore, the holomorphic splitting, $V=N\oplus N^{-1}K$ is an eigenbundle splitting for the gauge transformation $g_{\zeta\dwn{8}}.$

Let $g=g_{\zeta\dwn{8}}\oplus g_{\zeta\dwn{8}}^{-1}$ be the gauge transformation of $E=V\oplus V^*.$
Since the triple $(\bar\p_E,\phi,h)$ solves the Hitchin equations (\ref{Hit}), by Remark \ref{MetricRemark}, the triple 
\[(g^{-1}\bar \p_E g,g^{-1}\phi g,hg^{*_h}g)\]
 also solves (\ref{Hit}). 
We have computed 
$(g^{-1}\bar \p_E g,g^{-1}\phi g)=(\bar \p_E,\zeta^{-1}\dwn{4}\phi),$
thus $(\bar \p_E,\zeta\dwn{k}\phi,hg^{*_h}g)$ solves (\ref{Hit}) as well. 
Now, using the $\sU(1)$ action and Remark \ref{U1Remark}, the triple $(\bar \p_E,\phi,hg^{*_h}g)$ solves (\ref{Hit}). By uniqueness, $g$ is unitary: 
\[h=hg^{*_h}g.\]
Since $g_{\zeta\dwn{8}}$ is unitary and preserves the eigenbundle splitting $N\oplus N^{-1}K$, the metric $H$ splits as $H_1\oplus H_2.$ 
\end{proof}
\begin{Remark}\label{metricliftingRemark}
The metric solving \eqref{Hit} gives an equivariant harmonic map $h:\widetilde\Sigma\ra\sSp(4,\R)/\sU(2)\subset\sSp(4,\C)/\sSp(4).$ The holomorphic splitting $V=N\oplus N^{-1}K$ into line bundles gives an equivariant map $s:\widetilde\Sigma\ra\sSp(4,\C)/(\C^*\times\C^*).$ 
By Proposition \ref{MetricSplittingProposition}, the metric splits on the line bundles, and gives an equivariant map $f:\widetilde\Sigma\ra\sSp(4,\C)/(\sU(1)\times\sU(1))$ which makes the following diagram commute
\[\xymatrix{\sSp(4,\C)/(\C^*\times\C^*)&\sSp(4,\C)/(\sU(1)\times\sU(1))\ar[d]\ar[l]\\
\widetilde\Sigma\ar[u]^s\ar[ur]^{f}\ar[r]_h&\sSp(4,\C)/\sSp(4)}\]
The space $\sSp(4,\R)/(\sU(1)\times\sU(1))$ is the space of Hitchin triples (see Definition \ref{HitchinTripleDef}); we will show that the map $f$ is cyclic in the sense of Definition \ref{GG0cyclicsurfaceDef}.
\end{Remark}

\subsection{Parameterizations of $\Mm_J^d(\sSp(4,\R))$ and $\Rr_d(\sSp(4,\R))$} \label{ParaMd}

In this subsection we will only consider components $\Mm_J^d(\sSp(4,\R))$ for $g-1<d\leq 3g-3.$ 
By Proposition \ref{BGGMdProp}, Higgs bundles is such components are given a quadruple $(N,\mu,\nu,q_2)$ with $deg(N)=d$ and $\mu\neq0.$ 
This only describes representatives of the isomorphism classes of Higgs bundles in $\Mm_J^d(\sSp(4,\R)).$ In \em Proposition 3.28 \em of \cite{MaximalSP4}, it is shown that there is only a $1$-parameter family of gauge symmetries to account for: 
\[g_\lambda=\mtrx{\lambda&0\\0&\lambda^{-1}}:N\oplus N^{-1}K\ra N\oplus N^{-1}K.\] 
The gauge transformation $g_\lambda$ acts on the data $(N,\mu,q_2,\nu)$ by
\begin{equation}
	\label{fiberwiseaction}\xymatrix{(N,\mu,q_2,\nu)\ar[r]^{g_\lambda\ \ \ \ \ }&(N,\lambda^2\mu,q_2,\lambda^{-2}\nu)}.
\end{equation}

Since the $\C^*$-action of \eqref{fiberwiseaction} acts trivially on the holomorphic quadratic differential, we have the following parameterization of $\Mm_J^d(\sSp(4,\R)).$
\begin{Proposition}\label{FdJ}
	Let $Pic^d(\Sigma)$ be the variety of degree $d$ line bundles on $\Sigma.$  If
	\[\Ff_J^d=\{(N,\mu,\nu)\ |\ N\in Pic^d(\Sigma)\ \text{with}\ h^0(N^{-2}K^3)>0,\,\ \ 0\neq\mu\in H^0(N^{-2}K^3),\ \ \nu\in H^0(N^2K)\}/\C^*,\]
then the moduli space $\Mm_J^d(\sSp(4,\R))$ is then given by 
\[\Mm_J^d(\sSp(4,\R))=\Ff_J^d\times H^0(K^2).\]
\end{Proposition}

Let $\Pp_d\ra Pic^d(\Sigma)\times \Sigma$ be the Poincar\'e line bundle; it is defined so that the fiber over a point $(N,x)$ is the line $N_x.$ 
Define line bundles
\[\xymatrix{\Ll_1=\Pp_d^{-2}\otimes\pi^*_\Sigma K^3&\text{and}&\Ll_2=\Pp_d^{2}\otimes\pi^*_\Sigma K}\]
on $Pic^d(\Sigma)\times\Sigma.$  
Note that, if $\pi_{Pic}$ and $\pi_\Sigma$ denote the projections onto $Pic^d(\Sigma)$ and $\Sigma,$ then $(\pi_{Pic})_*\Ll_1$ and $(\pi_{Pic})_*\Ll_2$ are coherent sheaves on $Pic^d(\Sigma)$, and the stalks over a point $N\in Pic^d(\Sigma)$ are $ H^0(\Sigma, N^{-2}K^3)$ and $H^0(\Sigma, N^2K)$ respectively. The support of $(\pi_{Pic})_*\Ll_1$ is the subvariety of $Pic^d(\Sigma)$ consisting of linebundles $N$ with $h^0(N^{-2}K^3)>0.$ By Riemann-Roch, 
\[h^0(N^2K)=2d+g-1\]
is independent of $N\in Pic^d(\Sigma),$ and
hence, $(\pi_{Pic})_*\Ll_2$ is a vector bundle supported on all of $Pic^d(\Sigma)$. 
For $g-1<d<2g-2,$ a Riemann Roch computation shows $h^0(N^{-2}K^3)=-2d+5g-5.$ Thus, $(\pi_{Pic})_*\Ll_1$ is also a vector bundle supported on all of $Pic^d(\Sigma)$ for $g-1<d<2g-2.$

The space $\Ff_d^J$ is given by removing the zero section $\widetilde Z$ of the total space of $(\pi_{Pic})_*\Ll_1,$ taking the direct sum with the vector bundle $(\pi_{Pic})_*\Ll_2$ and quotienting by the fiberwise $\C^*$ action given by \eqref{fiberwiseaction}.
When $g-1<d<2g-2,$ the sheaf $(\pi_{Pic})_*\Ll_1$ is a vector bundle and this construction gives \em Proposition 3.29 \em of \cite{MaximalSP4}.
\begin{Proposition}\label{Prop3.29BGG}
	With the above set up, for $g-1<d<2g-2$, the space $\Ff_d^J$ is given by 
	\[\Ff_d^J=(((\pi_{Pic})_*\Ll_1-\widetilde Z)\oplus (\pi_{Pic})_*\Ll_2)/\C^*.\]
	Furthermore, $\pi:\Ff_d^J\ra Pic^d(\Sigma)$ is a fiber bundle, with fiber 
	$\pi^{-1}(N)\cong\Oo_{\P^{-2d+5g-6}}(1)^{\oplus 2d+g-1}.$
\end{Proposition}
For $2g-2\leq d\leq 3g-3$ the dimension of $H^0(N^{-2}K^3)$ depends on $N\in Pic^d(\Sigma)$ and not just on the degree. For $g-1<d\leq3g-3$ we have:
\begin{Proposition}\label{ParamMdProposition}
	There is a surjection from the smooth manifold $\Ff_J^d$ to the support of $(\pi_{Pic})_*(\Pp^{-2}_d\otimes\pi_\Sigma^*)$
	\[\pi:\xymatrix{\Ff_J^d\ar@{->>}[r]&supp((\pi_{Pic})_*(\Pp^{-2}_d\otimes\pi_\Sigma^*))} \subset\ Pic^d(\Sigma).\]
	If $N\in supp(\pi_{Pic})_*(\Pp^{-2}_d\otimes\pi_\Sigma^*))$ with $h^0(N^{-2}K^3)>1$ then the fiber is $\pi^{-1}(N)\cong\Oo_{\P^{a-1}}(1)^{\oplus b},$ where
  	$a=h^0(N^{-2}K^3) $ and  $b=h^0(N^2K)=2d+g-1$. 
  	If $N\in supp(\pi_{Pic})_*(\Pp^{-2}_d\otimes\pi_\Sigma^*))$ with $h^0(N^{-2}K^3)=1$  then $\pi^{-1}(N)\cong\C^{2d+g-1}.$
\end{Proposition}
\begin{proof}
If $[N,\mu,\nu]\in\Ff_J^d$ then $\pi([N,\mu,\nu])=N.$ For $N\in supp((\pi_{Pic})_*(\Pp^{-2}_d\otimes\pi_\Sigma^*)),$ set $a=h^0(N^{-2}K^3) $ and  $b=h^0(N^2K)=2d+g-1$; the fiber of $\Ff_J^d$ is given by $((\C^*)^a\oplus\C^b)/\C^*$ where the action is given by 
 	\[\lambda\cdot(x,y)=(\lambda^2x,\lambda^{-2}y).\]
 	For $a>1$, this is the total space of $\Oo_{\P^{a-1}}(1)^{\oplus b},$ and when $a=1$
 	$(\C^*\oplus\C^{2d+g-1})/\C^*\cong\C^{2d+g-1}.$
\end{proof}
 \begin{Remark}
  For $d=3g-3,$ the degree of $N^{-2}K^3$ is 0, hence, the support of $(\pi_{Pic})_*(\Pp^{-2}_d\otimes\pi_\Sigma^*)$ consists of the $2^{2g}$ points with $N^{-2}K^3=\Oo.$ 
  Thus, $\Ff_J^{3g-3}=\bigsqcup\limits_{j=1}^{2^{2g}} H^0(K^4)$ and $\Mm^{3g-3}_J$ consists of $2^{2g}$ copies of $H^0(K^4)\times H^0(K^2).$ These are the $2^{2g}$ Hitchin components.
 \end{Remark}
The moduli space $\Mm_J^d$ is homotopic to the space obtained by removing the zero section from the ``total space'' of the sheaf $(\pi_{Pic})_*(\Pp_d^{-2}\otimes\pi_\Sigma^* K^3)$ and quotienting by the $\C^*$ action. 
For $g-1<d<2g-2,$ the sheaf $(\pi_{Pic})_*(\Pp_d^{-2}\otimes\pi_\Sigma^* K^3)$ is a vector bundle supported on all of $Pic^d(\Sigma),$ and Proposition \ref{Prop3.29BGG} (Proposition 3.29 of \cite{MaximalSP4}) describes the homotopy type of $\Mm_J^d$ as a projective bundle over the torus $Pic^d(\Sigma).$
For $2g-2\leq d\leq 3g-3$ we can interpret the total space of the sheaf $(\pi_{Pic})_*\Ll_1$ as the kernel of a map of vector bundles over $Pic^d(\Sigma);$ this was suggested to the author by Tom Nevins. 

  \begin{Proposition} \label{kerhomotopyProp}
 If $D$ is an effective divisor on $\Sigma$ and $\Ll_1=\Pp_d^{-2}\otimes\pi^*_\Sigma K^3$ then 
  $(\pi_{Pic})_*\Ll_1$
  is the kernel of a map of vector bundles over $Pic^d(\Sigma)$
  \[\xymatrix{0\ar[r]&(\pi_{Pic})_*\Ll_1\ar[r]&(\pi_{Pic})_*(\Ll_1\otimes\pi_\Sigma^*(\Oo_\Sigma(D)))\ar[r]^{\psi\ \ \ }&(\pi_{Pic})_*(\Ll_1\otimes\pi_\Sigma^*(\Oo_\Sigma(D)/\Oo_\Sigma))}.\]
  The moduli space $\Mm_J^d(\sSp(4,\R))$ is homotopic to the space $(\psi^{-1}(Z)-\widetilde Z)/\C^*$ where $Z$ is the zero section of $(\pi_{Pic})_*(\Ll_1\otimes\pi_\Sigma^*(\Oo_\Sigma(D)/\Oo_\Sigma)),$ $\widetilde Z$ is the zero section of $(\pi_{Pic})_*(\Ll_1\otimes\pi_\Sigma^*(\Oo_\Sigma(D)))$ and $\C^*$ acts by fiberwise multiplication.
  \end{Proposition}
\begin{proof}
Fix an effective divisor $D$ on $\Sigma,$ and consider the exact sequence
\[\xymatrix{0\ar[r]&\Oo_\Sigma\ar[r]&\Oo_\Sigma(D)\ar[r]&\Oo_\Sigma(D)/\Oo_\Sigma\ar[r]&0}.\]
The sheaf $\Oo_\Sigma(D)/\Oo_\Sigma$ is a torision sheaf supported on $D.$ 
Since $\pi_\Sigma$ is flat, we have
\[\xymatrix{0\ar[r]&\pi_\Sigma^*(\Oo_\Sigma)\ar[r]&\pi_\Sigma^*(\Oo_\Sigma(D))\ar[r]&\pi_\Sigma^*(\Oo_\Sigma(D)/\Oo_\Sigma)\ar[r]&0}.\]
Note that $\left.\pi_\Sigma^*(\Oo_\Sigma(D)/\Oo_\Sigma)\right|_{\{N\}\times\Sigma}=coker(N\ra N\otimes\Oo(D)).$ 

Tensoring with a line bundle is exact and pushforward is left exact, thus, twisting by $\Ll_1$ and pushing forward by $\pi_{Pic}$ gives the following sequence of sheaves over $Pic^d(\Sigma)$
\[\xymatrix{0\ar[r]&(\pi_{Pic})_*\Ll_1\ar[r]&(\pi_{Pic})_*(\Ll_1\otimes\pi_\Sigma^*(\Oo_\Sigma(D)))\ar[r]^{\psi\ \ \ }&(\pi_{Pic})_*(\Ll_1\otimes\pi_\Sigma^*(\Oo_\Sigma(D)/\Oo_\Sigma))}.\]
Over a line bundle $N\in Pic^d(\Sigma),$ the stalk of the middle sheaf is $H^0(N^{-2}K^3\Oo(D)).$ Since we can choose the degree of $D$ arbitrarily large, $h^0(N^{-2}K^3\Oo(D))=-2d+5g-5+deg(D).$ 
Thus, the sheaf $(\pi_{Pic})_*(\Ll_1\otimes\pi_\Sigma^*(\Oo_\Sigma(D)))$ is a vector bundle.  
 The last sheaf in the sequence is a vector bundle on $Pic^d(\Sigma)$ of rank $deg(\Oo(D))$. If $Z\subset(\pi_{Pic})_*(\Ll_1\otimes\pi_\Sigma^*(\Oo_\Sigma(D)/\Oo_\Sigma))$ is the zero section, then
$(\pi_{Pic})_*\Ll_1=\psi^{-1}(Z).$

For the homotopy type of $\Mm_J^d,$ if $\widetilde Z$ is the zero section of the vector bundle $(\pi_{Pic})_*(\Ll_1\otimes\pi_\Sigma^*(\Oo_\Sigma(D))),$ then the total space of the sheaf $(\pi_{Pic})_*\Ll_1$ is given by $\psi^{-1}(Z)$ and
 	\[\widetilde Z\subset \psi^{-1}(Z)\subset(\pi_{Pic})_*(\Ll_1\otimes\pi_\Sigma^*(\Oo_\Sigma(D))).\]
 The surjective map from Proposition \ref{ParamMdProposition} is given by
 	\[\pi:(\psi^{-1}(Z)-\widetilde Z)\ra supp((\pi_{Pic})_*\Ll_1),\]
 and has fiber $\pi^{-1}(N)=H^0(N^{-2}K^3)-\{0\}.$ 
By Proposition \ref{ParamMdProposition} the moduli space $\Mm^d_J$ is homotopic to $(\psi^{-1}(Z)-\widetilde Z)/\C^*$ where $\C^*$ acts by fiberwise multiplication.

\end{proof}

Using Theorem \ref{uniqueminThm} and Proposition \ref{ParamMdProposition}, we obtain a mapping class group invariant parameterization of $\Rr_d(\sSp(4,\R))$ for $g-1<d\leq3g-3.$

\begin{Theorem}\label{MCGinvparamTheorem}
 	For $g-1<d\leq3g-3,$ let $\Rr_d(\sSp(4,\R))$ be the component of the maximal $\sSp(4,\R)$ representation variety corresponding to the Higgs bundle component $\Mm_J^d(\sSp(4,\R)).$ If $\pi:\Ff^d\ra\Tt(S)$ is the fiber bundle over Teichm\"uller space with $\pi^{-1}([J])= \Ff_J^d$ from Proposition \ref{FdJ}, then there is a mapping class group equivariant diffeomorphism 
 	\[\Psi:\Ff^d\lra\Rr_d(\sSp(4,\R)).\] 
 \end{Theorem}
The fact that $\Psi$ is a bijection follows as a corollary of Theorem \ref{uniqueminThm}, and the fact that $\Psi$ is a diffeomorphism follows as a corollary of the proof of Theorem \ref{uniqueminThm} (see Corollary \ref{MCGinvParamCor}).
\begin{Remark}
	For the Hitchin component, $\Ff_J^{3g-3}=\bigsqcup\limits_{j=1}^{2^{2g}}H^0(K^4)$ and we recover Labourie's mapping class group invariant parameterization of the Hithin component as a vector bundle over Teichm\"uller space.
\end{Remark}
 \begin{proof}(of bijection)
 Let $\rho_{J,N,\mu,\nu}\in\Rr_d(\sSp(4,\R))$ be the representation associated to the $\sSp(4,\R))$-Higgs bundle \[\left(N\oplus N^{-1}K,\mtrx{\nu&0\\0&\mu},\mtrx{0&1\\1&0}\right)\] over the Riemann surface $(S,J).$ 
 The map $\Psi$ is defined by 
 \[\xymatrix@R=.2em{\Ff^d\ar[rr]^{\Psi\ \ \ \ \ }&&\Rr_d(\sSp(4,\R))\\
 			(J,[N,\mu,\nu])\ar@{|->}[rr] &&\rho_{J,N,\mu,\nu}}\]
The inverse of $\Psi$ is defined by Theorem \ref{uniqueminThm},
\[\xymatrix@R=.2em{\Rr_d(\sSp(4,\R))\ar[rr]^{\ \ \ \ \ \Psi^{-1}}&&\Ff^d\\
 			\rho\ar@{|->}[rr] &&(J_\rho,[N,\mu,\nu])}\]
 	
 \end{proof}

\bigskip
\section{Lie theory background}\label{LieTheory}
In this section some Lie theory preliminaries are introduced in general and then in the specific case of $\fsp(4,\C)$. For the purpose of maximal $\sSp(4,\R)$ representations, we do not need the full generality of this section.
However, we wish to define cyclic surfaces in a general setting, prove a rigidity result and apply the result to the maximal $\sSp(4,\R)$ representations in the Gothen components. 

\subsection{Lie theory preliminaries in general} 
For all Lie theory notions, we follow \cite{knappbeyondintro} and \cite{realssliealgreps}. Let $\sG$ be a complex simple Lie group with Lie algebra $\fg$ and Killing form $B_\fg.$ 
Recall that a Cartan involution is a conjugate linear involution 
$\theta:\fg\ra\fg$ so that $-B_\fg(\cdot,\theta\cdot)$ is positive definite. The fixed point set $\fk$ of $\theta$ is the Lie algebra of a maximal compact subgroup $\sK_\theta\subset\sG.$ Cartan involutions exist and are unique up to conjugation. 
Furthermore, under the conjugation action, the stabilizer of a Cartan involution $\theta$ is the group $\sK_\theta.$ Thus we obtain:
\begin{Proposition}\label{symSpaceCartanInv}
	Let $\sG$ be a complex simple Lie group with maximal compact $\sK$ then
	\[\sG/\sK\cong\{\theta:\fg\ra\fg\ |\ \theta\ \text{a Cartan involution}\}.\]
\end{Proposition}

Let $\fc\subset\fg$ be a Cartan subalgebra, that is, a maximally abelian subalgebra consisting of semisimple elements. Cartan subalgebras exist, are unique up to conjugation. The dimension of $\fc$ is called the \em rank \em of $\fg$ and the Killing form $B_\fg|_{\fc\times\fc}$ is nondegenerate. 
An element $\alpha\in\fc^*$ is called a \em root \em if $\alpha\neq0$ and
\[\fg_\alpha=\{X\in\fg\ |\ [H,X]=\alpha(H)X\ \text{for\ all\ } H\in\fc\}\neq \{0\}.\]
 Denote the set of roots by $\Delta(\fg,\fc)\subset\fc^*.$ 
If $\alpha$ is a root, the space $\fg_\alpha$ is called the root space of $\alpha;$ the dimension of a root space $\fg_\alpha$ is always 1.  
Given two roots $\alpha,\beta\in\Delta(\fg,\fc),$ a simple, but fundamental calculation shows $[\fg_\alpha,\fg_\beta]\subset\fg_{\alpha+\beta}.$ 

Note that if $\alpha$ is a root, then $-\alpha$ is also a root. 
This allows us to choose a subset $\Delta^+\subset\Delta(\fg,\fc)$ of \em positive roots\em; $\alpha\in\Delta^+$ if and only if $-\alpha\notin\Delta^+.$
 A choice of positivity defines a set of \em simple roots \em 
 \[\Pi=\{\alpha_1,\cdots,\alpha_\ell\}\subset\Delta(\fg,\fc)\subset \fc^*,\]
where $\alpha\in\Delta^+$ implies $\alpha=\sum\limits_{j=1}^\ell n_i\alpha_i$ with $n_i\in\N$ and $\alpha_i\in\Pi.$ 
The integer $l(\alpha)=\sum_in_i$ is called the \em height \em or length of the root $\alpha.$ Let $m_\ell$ be the maximum height, then there is a unique root $\mu$ with $l(\mu)=m_\ell$ called the \em highest root \em.
\begin{Proposition}\label{GCcartanRem}
	The group $\sG$ acts transitively on the space of a Cartan subalgebra with a choice of simple roots, and the stabalizer of a point is the corresponding Lie group $\sC\subset \sG$ with Lie algebra $\fc;$ thus 
\begin{equation}\label{GCcartan}
	\sG/\sC\cong \{(\fc,\Delta^+)\ |\ \fc\subset\fg\ \text{a\ Cartan\ subalgebra\ ,\ } \Delta^+\subset\fc^* \text{a\ simple\ root\ system}\}.
\end{equation}
\end{Proposition}

Define $\fc(\R)=\{H\in\fc\ |\ \alpha(H)\in\R\ \text{for\ all\ }\alpha\in\Delta(\fg,\fc)\},$
then $\fc(\R)$ is a real form of $\fc.$ The Killing form $B_\fg$ is real and positive definite on $\fc(\R)$ and $\fc(\R)^*=Span\{\Delta(\fg,\fc)\}.$ Furthermore, the Killing form satisfies $B_\fg(X,Y)=0$ for $X\in\fg_\alpha,$ $Y\in\fg_\beta$ and $\alpha+\beta\neq 0.$ Thus $\fc$ and the vector subspaces $(\fg_\alpha\oplus\fg_{-\alpha})$ are pairwise orthogonal. 

Since the Killing form restricted to $\fc$ is nondegenerate, we can define the \em coroot \em $H_\alpha\in\fc$ of a root $\alpha$ by duality 
\[\beta(H_\alpha)=\frac{2B_{\fg^*}(\beta,\alpha)}{B_{\fg^*}(\alpha,\alpha)}.\]
By construction, $\alpha(H_\alpha)=2,$ thus $H_\alpha\in\fc(\R)$ and $\{H_{\alpha_i}\}_{i=1}^\ell$ forms a basis for $\fc(\R).$
A collection $\{X_{\alpha}\}_{\alpha\in\Delta}$ satisfying 
\begin{itemize}
	\item $[X_\alpha,X_{-\alpha}]=H_\alpha$
	\item $[X_\alpha,X_\beta]=N_{\alpha,\beta}\ X_{\alpha+\beta}$
	with $N_{\alpha,\beta}=-N_{-\alpha,-\beta}\in\N$ and $N_{\alpha,\beta}=0$ if $\alpha+\beta$ is not a root.
\end{itemize}  is called a \em Chevalley basis\em; Chevalley bases exist (Theorem 6.6 \cite{knappbeyondintro}). 
\begin{Definition}
	A Cartan involution which globally preserves a Cartan subalgebra $\fc$ is called a \em $\fc$-Cartan involution\em.
\end{Definition}
\begin{Lemma}\label{fcCartanInvLemma}
A $\fc$-Cartan involution takes a root space $\fg_\alpha$ to $\fg_{-\alpha}$ and $\theta(H_\alpha)=-H_\alpha.$ 
 \end{Lemma}
 \begin{proof}
 	Since $\theta$ is an isomorphism and $\alpha(H_\beta)$ real, for all $X\in\fg_\alpha,$ we have $\theta([H_\beta,X])=\alpha(H_\beta)\theta(X).$
 	 So $\theta$ takes root spaces to roots spaces. 
 	 Recall that for $\alpha+\beta\neq0,$ the root spaces $\fg_\alpha$ and $\fg_\beta$ are orthogonal. 
 	 By definitinion of a Cartan involution, $-B_\fg(\cdot,\theta\cdot)$ is positive definite. 
 	 Thus $\theta$ takes $\fg_\alpha$ to $\fg_{-\alpha}.$ 
 	 Let $X_{\pm\alpha}\in\fg_{\pm\alpha}$ with $[X_\alpha, X_{-\alpha}]=H_\alpha,$ then 
 	 \[\theta(H_\alpha)=[\theta(X_\alpha),\theta(X_{-\alpha})]=[\lambda_1 X_{-\alpha},\lambda_2 X_{\alpha}]=-\lambda_1\lambda_2H_\alpha.\] 
 	 Since, $\theta$ is an involution, $B_\fg(H_\alpha,H_\alpha)>0$ and $-B_\fg(H_\alpha,\theta(H_\alpha)>0$ we conclude $\theta(H_\alpha)=-H_\alpha.$
  \end{proof}

The existence of a Chevally basis gives the existence of two real forms, the split real form and the compact real form. 
The Lie subalgebra 
\begin{equation}\label{splitNonPTDS}
	\fg'=\bigoplus\limits_{i=1}^\ell\R H_{\alpha_i}\oplus \bigoplus\limits_{\alpha\in \Delta}\R X_\alpha
\end{equation}
is a split real form (Corollary 6.10 \cite{knappbeyondintro}). In terms of the Chevalley basis, $\fg'$ is the fixed point set of the conjugate linear involution $\lambda$ defined by $\lambda(X_\alpha)=X_{-\alpha}$ and $ \lambda(H_{\alpha_i})=H_{\alpha_i}.$
The subalgebra 
\[\fk=\bigoplus\limits_{i=1}^\ell \R iH_{\alpha_i}\oplus\bigoplus\limits_{\alpha\in\Delta(\fg,\fc)}\R(X_\alpha-X_{-\alpha})\oplus \R i(X_\alpha+X_{-\alpha})\]
is a compact real form of $\fg$ (Theorem 6.11 \cite{knappbeyondintro}). In terms of the Chevalley basis, $\fk$ is the fixed point set the conjugate linear Cartan involution defined by
\begin{equation}
	\label{ChevalleyCartanInv}\xymatrix{\theta(X_\alpha)=-X_{-\alpha}& \text{and}&\theta(H_\alpha)=-H_\alpha}.
\end{equation}

Following Kostant \cite{ptds}, we define the principal three dimensional subalgebra (PTDS) with respect to the Chevally basis. If $\{\epsilon_1,\cdots,\epsilon_\ell\}$ is the basis of $\fc$ dual to the simple roots, set  
\begin{equation}\label{gradingElementofPTDS}
	x=\sum\limits_{i=1}^\ell\epsilon_i\ =\ \haf\sum\limits_{\alpha\in\Delta^+}H_\alpha\ =\ \haf\sum\limits_{\i=1}^\ell r_{\alpha_i}H_{\alpha_i}\ .
\end{equation}
By construction of $x,$ if $X\in\fg_\alpha,$ then $[x,X]=l(\alpha)X.$ 
The eigenspace decomposition of $\fg$ with respect to $ad_x$ gives a $\Z$-grading on $\fg$ called the height decomposition:
\begin{equation}\label{Zgrading}
	\fg=\fg_{-m_\ell}\oplus\cdots\oplus\fg_{-1}\oplus\oplus\fc\oplus\fg_1\oplus\cdots\oplus\fg_{m_\ell}
\end{equation}
where $\fg_j=\bigoplus\limits_{l(\alpha)=j}\fg_\alpha.$ 
Define
\[\xymatrix{e_1=\sum\limits_{i=1}^\ell \sqrt{r_\alpha}X_{\alpha_i}&\text{and}&\tilde e_1=\sum\limits_{i=1}^\ell \sqrt{r_\alpha}X_{-\alpha_i}}.\]
By construction $\fs=\gen{\tilde e_1 , x , e_1}$ satisfies the bracket relations
\[\xymatrix{[x,e_1]=e_1\ ,&[x,\tilde e_1]=-\tilde e_1\ ,&[e_1,\tilde e_1]=x},\]
and thus $\fs\cong\fsl(2,\C).$ 

The adjoint action of $\fs$ on $\fg$ decomposes into a direct sum of irreducible $\fsl(2,\C)$-representations $\fg=\bigoplus V_j.$ Kostant \cite{ptds} showed that there are exactly $\ell=rank(\fg)$ irreducible summands
\begin{equation}\label{irrRepPTDS}
	\fg=\bigoplus\limits_{j=1}^\ell V_j.
\end{equation}
Furthermore, $dim(V_j)=2m_j+1$ and the integers $\{m_j\}$ are independent of all the choices. 
The numbers $\{m_1,\cdots,m_{\ell}\}$ are called the exponents of $\fg$ and always satisfy $m_1=1$ and $m_\ell=l(\mu).$ A three dimensional subalgebra with this property is unique up to conjugation \cite{ptds}. 
\begin{Definition}
Any subalgebra $\fs'$ conjugate to $\fs$ is called a \em principal three dimensional subalgebra (PTDS)\em, if $\fs'\cap\fc\neq \{0\}$ the PTDS is called a $\fc$-PTDS.
\end{Definition}

\begin{Theorem}\label{KostantThm42}
	(Theorem 4.2 \cite{ptds}) Let $\fs\subset \fg$ be PTDS and $x\in\fs$ be a semisimple element with centeralizer $\fg_x$. Then any other PTDS $\fs'\subset\fg$ containing $x$ is conjugate to $\fs$ by an element in Lie group $\sG_x$ with Lie algebra $\fg_x=Ker(ad_x).$
 \end{Theorem}

Let $e_j\in V_j$ be the highest weight vector, by definition, $[e_1,e_j]=0.$ 
Since $[x,e_\ell]=m_\ell e_\ell,$ one can always take $e_\ell=X_\mu,$ where $\mu$ is the highest root. 
The decomposition \eqref{irrRepPTDS} allows us to define the an involution $\sigma:\fg\ra\fg$ by 
\[\xymatrix{\sigma(e_j)=-e_j&\sigma(\tilde e_1)=-\tilde e_1}\]
and extended by the bracket relations.
\begin{Proposition}(\cite{ptds})
 	 The involution $\sigma$ commutes with the $\fc$-Cartan involution $\theta$ defined by $\theta(X_\alpha)=-X_{-\alpha}.$ Furthermore, the resulting real form $\lambda=\theta\circ\sigma$ is a split real form.
 \end{Proposition} 
\begin{Remark}
The involution $\sigma$ can be represented pictorial using the theory of irreducible representations of $\fsl(2,\C).$ For instance, when $\fg=\fsp(4,\C)$ the exponents $(m_1,m_2)=(1,3).$ The irreducible representations of equation \eqref{irrRepPTDS} have dimensions $(3,7)$ and the involution $\sigma$ is defined by:  
\begin{equation}\label{PTDSpictureSp4}
	\xymatrix@R=-.5em{&-3&-2&-1&0&1&2&3\\\\V_1&&&\underset{-1}{\overset{\tilde e_1}{\bullet}}&\underset{1}{\overset{x}{\bullet}}&\underset{-1}{\overset{e_1}{\bullet}}&&\\
V_2&\underset{-1}{\overset{ad_{\tilde e_1}^{6}(e_2)}{\bullet}}&\underset{1}{\bullet}&\underset{-1}{\bullet}&\underset{1}{\bullet}&\underset{-1}{\bullet}&\underset{1}{\bullet}&\underset{-1}{\overset{e_2}{\bullet}}}
\end{equation}
where the $\pm 1$ below each bullet is the value of the involution $\sigma,$ and the top row represents the height grading of \eqref{Zgrading}. 
By construction, the involution $\sigma$ is complex linear and preserves the height grading of \eqref{Zgrading}. In particular, it preserves the middle column which is the Cartan subalgebra $\fc.$ Thus $\{ad_{\tilde e_1}^{m_j}e_j\}$ generate the $\fc,$ and whenever $m_j$ is odd, $\sigma(ad_{\tilde e_1}^{m_j}e_j)=1.$   
\end{Remark}
Recall that we may take $e_{m_\ell}=X_\mu$ where $\mu$ is the highest root.
Since $\sigma$ commutes with $\theta,$ 
\[\sigma(X_{-\mu})=\theta(\sigma(\theta (X_{-\mu})))=-X_{-\mu},\]
and thus, $\sigma(H_\mu)=\sigma([X_\mu, X_{-\mu}])=H_\mu.$
Following Labourie \cite{cyclicSurfacesRank2}, we note that the involution $\sigma$ is unique.
\begin{Proposition}(Proposition 2.5.6 \cite{cyclicSurfacesRank2})
	Let $\fc$ be a Cartan subalgebra with a positive root system and $\mu$ the highest root. If $\sigma$ be an involution which preserves globally preserves $\fc$ and a $\fc$-PTDS $\fs$ with $\sigma(H_\mu)=H_{\mu},$ then $\sigma$ is unique. 
\end{Proposition}

The involutions $\theta$ and $\sigma$ give eigenspace decompositions 
\[\xymatrix{\fg=\fg^\theta\oplus i\fg^\theta&\fg=\fg^\sigma\oplus\fg^{-\sigma}.}\]
Since the compact form $\theta$ and the involution $\sigma$ commute, the restriction of $\sigma$ to the split real form $\fg_0=\fg^{\lambda}$ is a Cartan involution for $\fg_0:$
\[\fg_0=\fh\oplus\fm=(\fg^\theta\cap\fg^\sigma)\oplus(i\fg^\theta\cap\fg^{-\sigma}).\]
 Since both $\theta$ and $\sigma$ globally preserve $\fc,$ we may write $\fc^\lambda=\fc_0=\ft\oplus\fa$ where $\ft\subset\fh$ and $\fa\subset\fm,$
\[\fc=\ft_\C\oplus\fa_\C.\]
Recall that the coroots $\{H_\alpha\}\subset\fc$ are in the $(-1)$-eigenspace of the compact real form $\theta.$ 
In terms of the decomposition $\fc=\ft_\C\oplus\fa_\C,$ the $(-1)$-eigenspace of $\theta$ is $i\ft\oplus\fa.$ 
This leads to the notion of real, imaginary and complex roots. 
A root $\alpha\in(i\ft^*\oplus\fa^*)$ is called \em real \em if $\alpha|_{i\ft}=0,$ \em imaginary \em if $\alpha|_{\fa}=0$ and \em complex \em otherwise. 
By construction if $\alpha$ is real then $\theta(\alpha)=-\alpha,$ if $\alpha$ is imaginary then $\theta(\alpha)=\alpha$ and if $\alpha$ is complex then $\theta(\alpha)$ is a root different than $\alpha.$ An imaginary root $\alpha$ is called \em compact \em if $\fg_\alpha\subset\fh_\C,$ and \em noncompact \em if $\fg_\alpha\subset \fm_\C.$

By definition, the Cartan involution $\sigma|_{\fg_0}$ preserves the set of positive roots, so there are no real roots. 
Thus, the Cartan subalgebra $\fc_0$ is a maximally compact Cartan subalgebra (see Proposition 6.70 \cite{knappbeyondintro}). 
Furthermore, since $\sigma(e_1)=-e_1,$ by definition of $e_1,$ it follows that there are no imaginary compact simple roots. Thus, we have proven: 
\begin{Proposition}\label{Tmaxcompact}
The Cartan subalgebra $\fc_0\subset\fg_0$ is a maximally compact Cartan subalgebra and,
with respect to the Cartan involution $\sigma$ on $\fg_0,$ all simple roots are noncompact imaginary or complex. 
Furthermore, the subgroup $\sT\subset \sG$ with Lie algebra $\ft$ is a maximal compact torus of $\sG_0.$ 
\end{Proposition}
\begin{Remark}
	It is important to note that the split real form $\fg_0=\fg^{\sigma\circ\theta}$ is very different than the split real form $\fg'$ of equation \eqref{splitNonPTDS}. For $\fg_0$, the Cartan subalgebra $\fc$ is maximally compact, and for $\fg',$ the Cartan subalgebra $\fc$ is maximally noncompact. Thus
	\[\xymatrix{\fc\cap\fk\cap\fg_0=\ft\neq\emptyset&\fc\cap\fk\cap\fg'=\emptyset}.\]
\end{Remark}


\subsection{Lie theory preliminaries for $\fsp(4,\C)$}\label{sp4subsection} We now make all the above notions explicit for the Lie algebra $\fsp(4,\C)$.
Recall that if $\Omega=\mtrx{0&I\\-I&0},$ the Lie algebra $\fsp(4,\C)$ is the set of $4\times 4$ matrices $M$ so that $M^T\Omega+\Omega M=0.$ A simple computation shows
\begin{equation}
	\fsp(4,\C)=\left\{\left. \mtrx{X&Y\\Z&- X^T}\right| Y^T=Y,\ Z^T=Z\right\}.\label{matrixfsp4}
\end{equation}

Let $\fc\subset\fsp(4,\C)$ be a Cartan subalgebra consiting of the diagonal matrices
\[\fc=\gen{diag(1,0,-1,0), diag(0,1,0,-1)}=\gen{H_1, H_2}.\] 
If $L_1,L_2\in\fc^*$ is the basis dual to $\{H_1, H_2\},$ the set of roots is given by  
\[\Delta(\fsp(4,\C),\fc)=\{\pm L_i\pm L_j\}.\]
In terms of the matrix description \eqref{matrixfsp4}, a computation shows that the root spaces are
\begin{equation}\label{mtrxrootspacesp4}
	\mtrx{0&\fg_{L_1-L_2}&\fg_{2L_2}&\fg_{L_1+L_2}\\
			\fg_{-L_1-L_2}&0&\fg_{L_1+L_2}&\fg_{2L_1}\\
			\fg_{-2L_2}&\fg_{-L_1-L_2}&0&\fg_{-L_1+L_2}\\
			\fg_{-L_1-L_2}&\fg_{-2L_1}&\fg_{L_1-L_2}&0}.
\end{equation}
For any choice of positivity, there are two simple roots $\alpha_1,\alpha_2$ and the set of positive roots is given by
\[\{\alpha_1,\alpha_2, \alpha_1+\alpha_2,2\alpha_1+\alpha_2\}.\]
The choice of simple roots we will work with is $\{\alpha_1,\alpha_2\}=\{L_1+L_2,-2L_1\},$ the set of positive roots is then
\[\Delta^+=\{\alpha_1,\alpha_2, \alpha_1+\alpha_2,2\alpha_1+\alpha_2\}=\{L_1+L_2, -2L_1, -L_1+L_2, 2L_2\}.\]  
The coroots $\{H_\alpha\}$ are given by
\[\xymatrix{H_{L_1+L_2}=H_1+H_2& H_{-2L_1}=-H_1& H_{-L_1+L_2}=-H_1+H_2& H_{2L_2}=H_2.}\]
The height space decomposition \eqref{Zgrading} is given by 
\[\fsp(4,\C)=\bigoplus\limits_{j=-3}^3\fg_j=(\fg_{-2L_2})\oplus(\fg_{L_1-L_2})\oplus (\fg_{-L_1-L_2}\oplus\fg_{2L_1})\oplus \fc\oplus (\fg_{L_1L_2}\oplus\fg_{-2L_1})\oplus (\fg_{-L_1+L_2})\oplus (\fg_{2L_2})\]
For the Chevalley basis, set $X_\alpha\in\fg_\alpha$ to be the matrix with a 1 in the root space $\fg_\alpha.$ 
Then the compact real form $\theta$ given by $\theta(X_\alpha)=-X_{-\alpha}$ is  
 \[\theta(Y)=-\overline Y ^T.\]
The PTDS $\fs$ is defined by 
\[x=\mtrx{-\haf&&&\\&\frac{3}{2}&&\\&&\haf&\\&&&-\frac{3}{2}}\ \ \ \ \ \
e_1=\mtrx{&&0&\sqrt{\frac{3}{2}}\\
&&\sqrt{\frac{3}{2}}&0 \\
\sqrt{2}&0&&&\\
0&0&&
}\ \ \ \ \ \ \tilde e_1=\mtrx{&&\sqrt{2}&0\\
&&0&0 \\
0&\sqrt{\frac{3}{2}}&&&\\
\sqrt{\frac{3}{2}}&0&&
}.\]
The Lie algebra $\fsp(4,\C)$ has rank $2$ and the exponents $(m_1,m_2)=(1,3),$ thus the adjoint action of $\fs$ on $\fsp(4,\C)$ splits as a direct sum of two irreducible representations $\fsp(4,\C)=\fs\oplus V_2$ of dimensions $3$ and $7.$ 
The highest weight vector of $V_2$ can be taken to be $e_2=X_{2L_2}.$ 
The involution $\sigma$ is given by \eqref{PTDSpictureSp4}.
From the root space description \eqref{mtrxrootspacesp4}, the eigenspace decomposition of $\fsp(4,\C)$ with respect to the involution $\sigma$ is 
\[\fsp(4,\C)=\mtrx{X&0\\0&- X^T}\oplus\mtrx{0&Y\\Z&0}\cong\fgl(2,\C)\oplus\left(sym(\C^2)\oplus sym(\C^2)\right)=\fh_\C\oplus\fm_\C.\]
One checks that $\sigma\circ\theta=\theta\circ \sigma$ and that the fixed point set of $\lambda=\sigma\circ\theta$ is isomorphic to $\fsp(4,\R).$ 
Since $\sigma$ acts as the identity of the Cartan subalgebra $\fc,$ the fixed point subalgebra $\fc_0$ of the involution $\lambda$ is compact,
\[\fc=\fc_0\oplus i\fc_0=\ft\oplus i\ft.\] 
The Lie group $\sT\subset\sSp(4,\R)$ with Lie algebra $\ft$ is isomorphic to $\sU(1)\times\sU(1).$ 
Furthermore, since $\sigma$ acts as $(-1)$ on the simple root spaces, the simple roots $L_1+L_2$ and $-2L_1$ are noncompact imaginary roots, i.e. $\fg_{L_1+L_2},\fg_{-2L_1}\subset \fm_\C$. 

\begin{Remark}\label{MaxSp4dLiealgbundleRemark}
	Recall that we are interested in studying $\sSp(4,\R)$-Higgs bundles of the form $(V,\beta,\gamma)=\left(N\oplus N^{-1}K,\mtrx{\nu&q_2\\q_2&\mu},\mtrx{0&1\\1&0}\right).$ In terms of the root space decomposition \eqref{mtrxrootspacesp4}, the $\fsp(4,\C)$-Lie algebra subbundle of $End(N\oplus N^{-1}K\oplus N^{-1}\oplus NK^{-1})$ is given by
	\[(N^{-2}\otimes\fg_{-2L_2})\oplus (N^2K^{-1}\otimes\fg_{L_1-L_2})\oplus (K^{-1}\otimes\fg_{-L_1-L_2}\oplus N^{-2}K^2\otimes\fg_{2L_1})\oplus(\Oo\otimes\fc)\]\[\oplus (N^{2}K^{-2}\otimes\fg_{-2L_1}\oplus K^{1}\otimes\fg_{L_1+L_2})\oplus (N^{-2}K^{1}\otimes\fg_{-L_1+L_2})\oplus (N^{2}\otimes\fg_{2L_2}).\]
	The Higgs field $\phi=\left(\mtrx{\nu&q_2\\q_2&\mu},\mtrx{0&1\\1&0}\right)$ is given by
	\[\phi=1\otimes \fg_{-L_1-L_2}+\mu\otimes\fg_{2L_1}+q_2\otimes\fg_{L_1+L_2}+\nu\otimes\fg_{2L_2}. \]
\end{Remark}
\bigskip

 \section{The spaces of Cartan triples and Hitchin triples} \label{Homspaces}\smallskip We now define the main reductive homogeneous spaces we will study. The spaces we will be interested in are $\sG/\sT$ and $\sG/\sT_0$ where $\sG$ is a complex simple Lie group and $\sT$ is a maximal compact torus of $\sG$ and $\sT_0$ is the maximal compact torus of a split real form of $\sG_0\subset\sG.$ We start by considering a more geometric set of objects.
\begin{Definition}\label{CartanTripleDef}
 	A \em Cartan triple is a triple \em $(\fc,\Delta^+,\theta)$ where
 	\begin{itemize}
 	\item $\fc\subset \fg$ is a Cartan subalgebra
 	\item $\Delta^+\subset\fc^*$ is a choice of positive roots
 	\item $\theta$ is a $\fc$-Cartan involution 
 	\end{itemize}
 \end{Definition} 

 Let $\sT\subset\sG$ be a maximal compact torus, Proposition \ref{GCcartan} and Lemma \ref{fcCartanInvLemma} imply the following proposition. 

\begin{Proposition}\label{GCKcartanProp}
The space of Cartan triples is isomorphic to $\sG/\sT$
\end{Proposition}
 Note that we could equivalently define $(\fc,\Delta^+,\theta)$ to be a Cartan triple where $\theta$ is Cartan involution, and $(\fc,\Delta^+)$ a Cartan subalgebra with positive root system and $\fc$ is preserved by $\theta.$ There are natural projection maps 
 \[\xymatrix@=1.3em{\sG/\sC&\sG/\sT\ar[l]_{\pi_1}\ar[d]^{\pi_2}\\&\sG/\sK}\]
 where $\pi_1(\fc,\Delta^+,\theta)=(\fc,\Delta^+)$ and $\pi_2(\fc,\Delta^+,\theta)=\theta.$
\begin{Definition}\label{HitchinTripleDef}
	A \em Hitchin triple \em  is a triple $(\Delta^+\subset\fc^*,\theta,\lambda)$ where 
	\begin{itemize}
		\item $\fc$ is a Cartan subalgebra
		\item   $\Delta^+\subset\Delta(\fg,\fc)\subset\fc^*$ is a choice of positive roots
		\item $\theta$ is a $\fc$-Cartan involution which globally preserves a PTDS $\fs$ which contains $x=\haf\sum\limits_{\alpha\in\Delta^+}H_\alpha.$ 
		\item $\lambda$ is a split real form which commutes with $\theta,$ globally preserves $\fc$, globally preserves a PTDS $\fs$ which contains $x=\haf\sum\limits_{\alpha\in\Delta^+}H_\alpha$ and satisfies $\lambda(H_\mu)=-H_\mu.$  
	\end{itemize}
\end{Definition}
\begin{Proposition}Let $\sG$ be a complex simple Lie group, and $\sG_0$ be a split real form of $\sG.$ The space of Hitchin triple is diffeomorphic to $\sG/\sT_0$ where $\sT_0$ is the maximal compact torus of $\sG_0.$
\end{Proposition}
\begin{proof}
We first show that the $\sG$ acts transitively on the space of Hitchin triple. Let $(\Delta_1^+\subset\fc_1^*,\theta_1,\lambda_1)$ and $(\Delta_2^+\subset\fc_2^*,\theta_2,\lambda_2)$ be two such Hitchin triples. By Remark \ref{GCcartanRem}, we can conjugate $(\Delta_2^+\subset\fc_2^*)$ to $(\Delta_1^+\subset\fc_1^*).$ 
Thus we may assume $(\Delta_1^+\subset\fc_1^*)=(\Delta_2^+\subset\fc_2^*).$ 
 Let $x=\haf\sum\limits_{\alpha\in\Delta_1^+}H_\alpha,$ and suppose $\theta_1$ stabilizes an $\fc$-PTDS $\fs_1$ and $\theta_2$ stabilizes an $\fc$-PTDS $\fs_2$ with $x\in\fs_1$ and $x\in\fs_2.$ 
 By Theorem 4.2 of \cite{ptds} (q.f. Theorem \ref{KostantThm42}), the PTDSs $\fs_1$ and $\fs_2$ are conjugate via an element of $\sC.$ 
Thus we may assume $\fs_1=\fs_2.$ 
Since $\theta_1$ and $\theta_2$ are both $\fc$-Cartan involutions, $\theta_1|_\fc=\theta_2|_\fc.$ 
Furthermore, $\theta_1$ and $\theta_2$ are both $\fc\cap\fs$-Cartan involutions of $\fs,$ by Proposition \ref{GCKcartanProp}, $\theta_2|_\fs$ can be conjugated to $\theta_1|_\fs$ by an element of the subgroup $\sC'\subset\sC$ with Lie algebra $\fc\cap\fs.$ 
Observe that conjugating by $\sC'$ preserves $(\Delta_1^+\subset\fc_1^*,\fs_1).$ 
Furthermore, $\fg$ is generated by $\fc+\fs,$ thus after conjugating by such an element of $\sC'$, we obtain $\theta_1=\theta_2.$ 
Since $\theta_1=\theta_2$ and $\fs_1=\fs_2,$ by uniqueness of the involution $\sigma,$ the splits real forms $\lambda_1$ and $\lambda_2$ are equal.   

The stabilizer of $(\Delta^+\subset\fc,)$ is a maximal torus $\sC,$ and the stabilizer of a $\fc$-Cartan involution is $\sC\cap\sK.$
 The stabilizer of the split real form $\lambda$ is the corresponding split real group $\sG_0\subset \sG.$
Thus the stabilizer of a Hitchin triple $(\Delta^+\subset\fc,\theta,\lambda)$ is $\sT_0=\sG_0\cap\sK\cap\sC.$
\end{proof}
\begin{Remark}
	A real form $\sG_0$ is called a \em group of Hodge type \em if the maximal compact torus $T_0\subset\sG_0$ is a maximal compact torus of the complex group $\sG.$ 
	For split real forms, only $\sSL(n,\R),\ \sSO(2n+1,2n+1),$ and the split real form of $\sE_6$ are \em not \em of Hodge type.
	When a split real form $\sG_0$ is of Hodge type, the space of Cartan triples and the space of Hitchin triples are the same. 
	In this case, the involution $\sigma$ determined by a $\fc$-PTDS containing $x=\haf\sum\limits_{\alpha\in\Delta^+}H_\alpha$ acts as $+Id$ on $\fc,$ and $\fc=\ft\oplus i\ft.$
\end{Remark}

	

\subsection{Reductive homogeneous spaces}
We now recall some important geometry of reductive homogeneous spaces, the main reference is the first chapter of \cite{TwistorTheoryHarmonicMapsBOOK}.
Let $M$ be a manifold with a smooth transitive action of $\sG.$ If we fix a base point $x_0\in M$ and define $\sH=Stab_\sG(x_0),$ then, since the action is transitive, we have a prinicpal $\sH$-bundle
\[\xymatrix@R=0em@C=1em{\sH\ar@{-}[r]&\sG\ar[rr]^\pi&&M\\
     			&g\ar@{|->}[rr]&& g\cdot x_0}
\]
Thus, the tangent bundle is given by $TM=\sG\times_H\fg/\fh.$ A homogeneous space $M$ is called \em reductive \em if the Lie algebra $\fg$ has a decomposition $\fg=\fh\oplus\fm$ as $Ad_\sH$-modules. All homogeneous spaces will be assumed to be reductive. 
If $W$ is a linear representation of $\sH$ we will denote the associated bundle by
$\sG\times_\sH W=[W].$ Thus, \[[\fm]\cong TM.\]


 Since $\fm$ is an $Ad_\sH$-invariant subspace of $\fg$, we have $[\fm]\subset [\fg].$ 
The action of $\sH$ on $\fg$ is the restriction of the $\sG$ action, hence $[\fg]$ is trivializable
\[\xymatrix@R=.2em{\sG\times_\sH \fg\ar@{<->}[rr]^\cong &&M\times \fg\\
[g,\xi]\ar@{|->}[rr]&&(\pi(g),Ad_g\xi)}
\]
\begin{Example}
	When $\sG$ is a complex simple Lie group with maximal compact $\sK$, the symmetric space $\sG/\sK$ is a reductive homogeneous space. Furthermore, since $\fk\otimes\C=\fg,$ we have 
\[T(\sG/\sK)\otimes \C\cong [\fk]\oplus[i\fk]=[\fg]\cong \sG/\sH\times \fg.\] 
\end{Example}
Using $[\fm]\cong TM,$ we have $TM\subset M\times \fg.$ This inclusions can be interpreted as an equivariant 1-form on $M$ valued in $\fg,$ $\omega\in\Omega^1(M,\fg).$ 
\begin{Definition}\label{MCFormDef}
	The equivariant $\fg$-valued 1-form $\omega\in\Omega^1(M,\fg)$ is called the \em Maurer Cartan form of the homogeneous space $M.$\em
\end{Definition}

 We will view a reductive homogeneous space as coming equipped with a fixed summand $\fm\subset \fg.$ 
 Let $\omega_\sG\in\Omega^1(\sG,\fg)^\sG$ be the left Maurer-Cartan form of $\sG,$ it is $\sG$-equivariant. Since $\fg=\fh\oplus\fm,$ we may split 
 \[\omega_\sG=P_\fh\omega_\sG\oplus P_\fm\omega_\sG.\]
 This an $Ad_\sH$-invariant splitting since $\fg=\fh\oplus\fm$ is $Ad_\sH$-invariant, thus 
\[P_\fh\omega_\sG\in\Omega^1(\sG,\fh)^\sH\ \ \ \ \text{and}\ \ \ \ \  P_\fm\omega_\sG\in\Omega^1(\sG,\fm)^\sH.\]
The form $P_\fh\omega_\sG$ is a connection on the principal $\sH$-bundle $\sG\ra M$ which we call the \em canonical connection\em. 
For any $\H$-representation $V$, the canonical connection induces a covariant derivative $\nabla^c$ on any associated bundle $[V].$ 
By construction, if $s\in C^\infty(M,[V])$ is $\sG$-equivariant, then $\nabla^cs=0.$ 
The form $P_\fm\omega_\sG$ is an equivariant horizontal $1$-form, i.e. it vanishes along vector fields induced by the action. 
Thus, $P_\fm\omega_\sG$ descends to a 1-form on $M$ valued in $[\fm]$ which is the Maurer Cartan form $\omega.$


When $V$ is the restriction of a representation of $\sG,$ $[V]$ is trivializable,
in which case, there is a simple relationship between flat differentiation on $M\times V$ and covariant differentiation by the canonical connection. This will be important for our later considerations of cyclic surfaces and the Hitchin equations.
\begin{Lemma}(see chapter 1 \cite{TwistorTheoryHarmonicMapsBOOK})\label{CanonicalConnectionandTorsion}
 	Let $f:M\ra M\times V$ be a smooth section, then 
 	\[df=\nabla^c f+\omega\cdot f.\]
 	If $V=\fg$ is the adjoint representation, then $\nabla^c=d-ad_\omega$ and the torision of the canonical connection on $TM=[\fm]$ is given by
 	\[T_{\nabla^c}=-\haf[\omega,\omega]^\fm.\]
 \end{Lemma} 
\begin{Remark}
	\label{parallelsubbundles}
	If $[\ ,\ ]^\fm$ and $[\ ,\ ]^\fh$ denote the projections onto $[\fm]$ and $[\fh],$ then the flatness of $d$ can be written in terms of $\nabla^c$ and $\omega$ as
\[\begin{dcases}
	F_{\nabla^c}+\haf[\omega,\omega]^\fh=0 & \fh-part\\
	d^{\nabla^c}\omega+\haf[\omega,\omega]^\fm=0& \fm-part
\end{dcases}\]
Moreover, if we decompose $\fm=\bigoplus\limits_{j}\fm_j$ into irreducible $\sH$-representations then 
\[TM\cong \bigoplus\limits_{j}\sG\times_\sH\fm_j= \bigoplus\limits_{j}[\fm_j].\]
This gives a decomposition of the trivial bundle $\underline \fg\ra M$ as a direct sum of $\nabla^c$-parallel vector bundles 
\[[\fg]=\underline{\fg}=[\fh]\oplus\bigoplus\limits_j[\fm_j].\]
Furthermore, the Maurer Cartan form decomposes $\omega=\sum\limits_j \omega_j,$ and the zero curvature equations are
\[\begin{dcases}
	F_{\nabla^c}+\haf\sum\limits_{j,k}[\omega_j,\omega_k]^\fh=0 & \fh-part\\
	d^{\nabla^c}\omega_j+\haf\sum\limits_{k,\ell}[\omega_k,\omega_\ell]^{\fm_j}=0& \fm_j-part
\end{dcases}\]
\end{Remark}

\begin{Example}
When $\sG$ is a semisimple Lie group, with $\sK\subset\sG$ a maximal compact, any $G$-invariant metric on $\sG/\sK$ is a $G$-equivariant section of an assoicated bundle. 
Thus, the canonical connection $\nabla^c$ is a metric connection. 
Since $[\fm,\fm]\subset\fk,$ the torsion of the canonical connection vanishes. 
Hence, for a symmetric space $\sG/\sK,$ the canonical connection is the Levi Civita connection of any $G$-invariant metric.
Furthermore, the flatness equations decompose as
 \[\begin{dcases}
	F_{\nabla^c}+\haf[\omega,\omega]=0 & \fk-part\\
	d^{\nabla^c}\omega=0& \fm-part
\end{dcases}\]
 \end{Example}
 Recall from Proposition \ref{symSpaceCartanInv}, the symmetric space $\sG/\sK$ is the space of Cartan involutions. The following lemma will be important for our defintion of cyclic surfaces. 
 \begin{Lemma}
 	\label{canonicalAuts} Let $\underline{\fg}=[\fk]\oplus[\fm]$ denote the trivializable Lie algebra bundle over  $M=\sG/\sK$ the symmetric space of Cartan involutions of $\fg.$ There is a canonical automorphism $\Theta:\underline\fg\lra\underline\fg$ given by \[\Theta(\theta,X)=(\theta, \theta(X)).\]
 	Furthermore, the invariant metric on $\underline{\fg}$ induced by the Killing form is given by 
 	$B_\Theta(X,Y)=-B_\fg(X,\Theta(Y)),$ and is parallel with respect to the canonical connection, $\nabla^cB_\Theta=0.$
 \end{Lemma}
\begin{Remark}\label{Thetaforms}
The automorphism $\Theta$ has a natural extension to complex forms valued in $\underline\fg.$ If $\alpha\in\Omega^*(\sG/\sK, \underline\fg)$ is of the form $\alpha=A\cdot a$ where $A\in\Omega^*(\sG/\sK)$ and $a$ is a section of $\underline\fg,$ then 
$\Theta(\alpha)=\overline A\cdot \Theta(a).$
 \end{Remark}
\begin{Proposition}\label{MapstoSymmSpaceProp}
	Let $N$ be a simply connected manifold and $(\widetilde\fg,\widetilde D)$ be a flat $\fg$-bundle. Suppose
	\begin{itemize}
		\item $\widetilde\Theta:\widetilde\fg\ra\widetilde\fg$ be a smoothly varying Cartan involution with $\widetilde\fg=\widetilde\fk\oplus\widetilde\fm$ the corresponding eigenbundle decomposition.  
		\item  $\widetilde\nabla$ a connection with $\widetilde\nabla\widetilde\Theta=0$
		\item $\widetilde\omega\in\Omega^1(N,\widetilde\fm)$ with $\widetilde D=\widetilde\nabla+ad_{\widetilde\omega}.$ 
	\end{itemize} 
	Then there exists a map $f:N\ra\sG/\sK,$ unique up to postcomposition by an element of $\sG$ so that 
	\[f^*(\ \underline\fg\ ,\ \Theta\ ,\ \nabla^c\ ,\ \omega\ )\ =\ (\ \widetilde\fg\ ,\ \widetilde\Theta\ ,\ \widetilde\nabla\ ,\ \widetilde\omega\ ).\]
\end{Proposition}
\begin{proof}
	Since $N$ is simply connected, choose a trivialization $(\widetilde\fg,\widetilde D)=(N\times\fg, d).$ 
	In this trivialization, the gauge transformation $\widetilde\Theta$ defines the map $f:N\ra\sG/\sK$ with $(f^*\underline\fg,f^*\Theta)=(\widetilde\fg,\widetilde \Theta).$ 
	Another trivialization produces a map which differs from $f$ by postcomposition by an element of $\sG.$

	Thus, $\widetilde\Theta$ is parallel with respect to $f^*\nabla^c$ and $\widetilde\nabla.$ Since the stabilizer of a Cartan involution is $\sK,$ we have $f^*\nabla^c-\widetilde\nabla\in\Omega(N,\widetilde\fk),$ and thus 
	$f^*\omega-\widetilde\omega\in\Omega^1(N,\widetilde\fk).$
	But $f^*\omega-\widetilde\omega\in\Omega^1(N,\widetilde\fm),$ thus $f^*\nabla^c=\widetilde\nabla$ and $f^*\omega=\widetilde\omega.$ 
  \end{proof}
    \begin{Remark}
  	\label{Higgsfieldanddh}
  	Given a map $f:N\ra\sG/\sK,$ by definition of the Maurer Cartan form $\omega,$ the form $f^*\omega$ is  derivative of the map $df.$  
  \end{Remark}
Since the canonical connection $\nabla^c$ is the Levi Civita connection of any $\sG$-invariant metric on $\sG/\sK,$ the harmonic map equations for a map $f:(N,g)\ra\sG/\sK$ are given by $(D^{f^*\nabla^c})^*(df)=0.$ The following corollary is immediate.
  \begin{Corollary}\label{UniqueMaptoSymmSpaceCor}
  	Let $N$ be a smooth manifold and $\widetilde\fg\ra N$ be a flat $\fg$-Lie algebra bundle equipped with the structure of Proposition \ref{MapstoSymmSpaceProp}, then there exists
 	\begin{enumerate}
 	\item A representation $\rho:\pi_1(N)\ra \sG$ unique up to conjugation 
 	\item A $\rho$-equivariant map $f$ from the universal cover $\widetilde N$ of $N$ to the space of Cartan involutions $\sG/\sK$ satisfying the conclusion of Proposition \ref{MapstoSymmSpaceProp}. 
  	\end{enumerate}
  	Furthermore, if $\widetilde\nabla^*\widetilde\omega=0,$ then the map $f$ is harmonic.
  \end{Corollary}

\subsection{Maurer Cartan form for $\sG/\sT$ and $\sG/\sT_0$} 
Let $M$ be the space of Cartan triples of Definition \ref{CartanTripleDef}, then $M\cong\sG/\sT$ where $\sT$ is the maximal compact torus of $\sG.$ 
If $(\fc,\Delta^+,\theta)$ is a Cartan triple, let $\ft=\fc^\theta,$ then $\ft$ is the Lie algebra of $\sT.$ We have the following $Ad_\sT$ invariant decompositions
\[\xymatrix{\fg=\ft\oplus i\ft \oplus  \bigoplus\limits_{\alpha\in\Delta(\fg,\fc)} \fg_\alpha&\text{and}&
\fg=\fk\oplus i\fk.}\]
Thus, the Lie algebra bundle $[\fg]=\underline\fg\ra \sG/\sT$ has corresponding compatible $\nabla^c$-parallel decompositions 
\[\xymatrix{\underline\fg=[\ft]\oplus [i\ft]\oplus  \bigoplus\limits_{\alpha\in\Delta(\fg,\fc)} [\fg_\alpha]
&\text{and}&
\underline\fg=[\fk]\oplus [i\fk].}\]
Recall that $T\sG/\sT\cong[i\ft]\oplus  \bigoplus\limits_{\alpha\in\Delta(\fg,\fc)} [\fg_\alpha],$ thus the Maurer-Cartan form vanishes of $[\ft],$ i.e. $\omega|_{[\ft]}\equiv 0.$

If $\ell=rank(\fg),$ then a set simple roots gives a $\Z^\ell$-grading of $\fg$ called the root space decomposition 
\[\fg=\fc\oplus\bigoplus\limits_{\alpha\in \Delta(\fg,\fc)}\fg_\alpha.\]
 Since this decomposition is $Ad_\sT$-invariant and $[\fg_\alpha,\fg_\beta]\subset \fg_{\alpha+\beta},$ the zero curvature equations decomposes as 
 \begin{equation}\label{ZellgradEq}
 	\begin{dcases}
	F_{\nabla^c}+\sum\limits_{\alpha\in\Delta^+(\fg,\fc)}[\omega_\alpha,\omega_{-\alpha}]^{\ft}=0 & \ft-part\\
	d^{\nabla^c}\omega_{i\ft}+\sum\limits_{\alpha\in\Delta^+(\fg,\fc)}[\omega_\alpha,\omega_{-\alpha}]^{i\ft}=0 & i\ft-part\\
	d^{\nabla^c}\omega_\alpha+[\omega_{i\ft},\omega_\alpha]+\sum\limits_{\substack{\beta,\gamma\in\Delta(\fg,\fc)\\ \alpha=\beta+\gamma}}[\omega_\beta,\omega_\gamma]=0& \fm_\alpha-part
\end{dcases}
 \end{equation} 

Recall that if $\{\alpha_i\}$ is the collection of simple roots, then every root $\alpha$ can be written uniquely as $\alpha=\sum\limits n_i\alpha_i,$ and the integer $\ell(\alpha)=\sum n_i$ is called the height of $\alpha$. 
From equation \eqref{Zgrading}, the grading element $x$ from the PTDS $\fs$ gives a $\Z$-grading on $\fg$
 	\[\fg=\fg_{-m_\ell}\oplus\cdots\oplus\fg_{-1}\oplus\oplus\fc\oplus\fg_1\oplus\cdots\oplus\fg_{m_\ell}\]
 	where $\fg_j=\bigoplus\limits_{l(\alpha)=j}\fg_\alpha.$ 
Since, $[\fg_j,\fg_k]\subset\fg_{j+k}$, in terms of the height decomposition, the flatness equations decompose as
 \begin{equation}\label{ZgradEq}
 \begin{dcases}
	F_{\nabla^c}+d^{\nabla^c}\omega_0+\sum\limits_{j>0}[\omega_j,\omega_{-j}]=0 & \fc-part\\
	d^{\nabla^c}\omega_j+\haf\sum\limits_{k}[\omega_k,\omega_{j-k}]=0& \fg_j-part
\end{dcases}
 \end{equation} 
 Set $g_+=exp(\frac{2\pi i\cdot x}{m_\ell+1}),$ and consider the autormorphism $Ad_{g_+}:\fg\ra\fg.$ 
 Since $ad(x)$ acts on $\fg_j$ with eigenvalue $j,$ the automorphism $Ad_{g_+}$ acts on $\fg_j$ with eigenvalue $\zeta\dwn{m_\ell+1}^j=e^\frac{2\pi i\cdot j}{m_\ell+1}.$ Note that, by construction, $Ad_{g_+}(X)=X$ if and only if $X\in \fc.$
An eigenspace decomposition of $Ad_g$ gives a $\Z/(m_\ell+1)\Z$-grading on $\fg$:
\[\fg=\bigoplus\limits_{j\in\Z/(m_\ell+1)\Z}\widehat\fg_j\]
where $\widehat\fg_j=\bigoplus\limits_{k=j\ mod\ m_\ell+1}\fg_k.$ 
The Maurer Cartan form decomposes as
\begin{equation}
 	\label{ZmellMCdecom}\omega=\sum\limits_{j\in\Z/(m_\ell+1)\Z}\widehat\omega_j
 \end{equation} and the flatness equations decompose as
\begin{equation}\label{ZmellgradingEq}
\begin{dcases}
	F_{\nabla^c}+d^{\nabla^c}\widehat\omega_0+\sum\limits_{j>0}[\widehat\omega_j,\widehat\omega_{-j}]=0 & \widehat\fg_0=\fc-part\\
	d^{\nabla^c}\widehat\omega_j+\haf\sum\limits_{k}[\widehat\omega_k,\widehat\omega_{j-k}]=0& \widehat\fg_j-part
\end{dcases}
\end{equation} 
\begin{Remark}\label{ZmZgradingRemark}
This grading will be essential for our definition of cyclic surfaces.
The automorphism $Ad_{g_+}$ makes the space $\sG/\sT$ into a $(m_\ell+1)$-symmetric space. It will be important that the subspaces $\widehat\fg_{\pm 1}$ are 
\begin{equation}\label{gradingequivs}
	\xymatrix@C=1em{\widehat\fg_1=\fg_1\oplus\fg_{-m_\ell}=\fg_{\alpha_1}\oplus\cdots\oplus\fg_{\alpha_\ell}\oplus\fg_{-\mu}&\text{and}& \widehat{\fg}_{-1}=\fg_{-1}\oplus\fg_{m_\ell}=\fg_{-\alpha_1}\oplus\cdots\oplus\fg_{-\alpha_\ell}\oplus\fg_{\mu}}
\end{equation}
where $\{\alpha_i\}$ is the set of simple roots and $\mu$ is the highest root. Furthermore, the compact involution $\theta$ maps $\widehat\fg_1$ to $\widehat\fg_{-1}.$ 
\end{Remark}
For the space of Hitchin triples $\sG/\sT_0,$ the Cartan subalgebra decomposes as $\fc=\ft_0\oplus i\ft\oplus\fa\oplus i\fa.$ The tangent bundle is given by 
\[T\sG/\sT_0=[i\ft_0]\oplus[\fa]\oplus[i\fa]\oplus\bigoplus\limits_{\alpha\in\Delta(\fg,\fc)}[\fg_\alpha]\]
and Maurer Cartan form vanishes on $[\ft_0].$ The decompositions \eqref{ZellgradEq}, \eqref{ZgradEq}, and \eqref{ZmellgradingEq} of the flatness equations still hold. 
\begin{Lemma}\label{propertiesOfTriplesLem}
	Let $\ft_0\oplus\fm=\fg$ be the reductive decomposition corresponding to a Hitchin triple. The trivial Lie algebra bundle $\underline\fg\ra\sG/\sT$ has the following data 
	\begin{itemize}
	\item $\omega\in\Omega^1(\sG/\sT,[\fm]\subset\underline\fg)$ the Maurer Cartan form 
	\item the canonical connection $\nabla^c$ with flat differentiation given by $d=\nabla^c+ad_\omega$
	\item $[\fc]\subset\underline\fg$ which decomposes as $[\fc]=[\ft_0]\oplus [i\ft_0]\oplus[\fa]$
	\item $\nabla^c$-parallel subbundles $[\fn^+]\subset\underline\fg$ and $[\fn^-]\subset\underline\fg$ with $[\fn^-]\oplus[\fc]\oplus [\fn^+]=\underline\fg.$
 	\item $\nabla^c$-parallel conjugate linear involution $\Theta:\underline\fg\ra\underline\fg$ and $\Lambda:\underline\fg\ra\underline\fg$ with fixed point set $[\fk]$ and $[\fg_0].$
 	\item A $\nabla^c$-parallel complex linear involution $\sigma=\Theta\circ\Lambda$ with eigenbundle decomposition $\underline\fg=[\fh_\C]\oplus[\fm_\C],$ where $[\fh]\subset[\fg_0]$ is the fixed point set of $\Theta|_{[\fg_0]}.$
 	\item A $\nabla^c$-parallel order $(m_\ell+1)$ automorphism $\Xx_+:\underline\fg\ra\underline\fg$ with eigenbundles $[\widehat\fg_j]$ and $[\fc]$ the identity eigenbundle.  
 \end{itemize} 
\end{Lemma}
\begin{proof}
The splitting of $\fg$ into root space is $Ad_{\sT_0}$ invariant, thus we have $\nabla^c$-parallel subbundles 
	\[\xymatrix{[\fn^+]=\bigoplus\limits_{\alpha\in\Delta^+}[\fg_\alpha]&\text{and}&[\fn^-]=\bigoplus\limits_{\alpha\in\Delta^-}[\fg_\alpha]}.\]
	The fiber of $[\fn^+]$ over a Hitchin triple $(\Delta^+\subset\fc^*,\theta,\lambda)$ is $\bigoplus\limits_{\alpha\in\Delta^+}\fg_\alpha.$
	For $X\in\fg,$ the conjugate linear involutions $\Theta$ and $\Lambda$ are defined by 
	\[\xymatrix@C=1.2em{\Theta((\Delta^+\subset\fc^*,\theta,\lambda),X))=((\Delta^+\subset\fc^*,\theta,\lambda),\theta(X))&\text{and}&\Lambda((\Delta^+\subset\fc^*,\theta,\lambda),X))=((\Delta^+\subset\fc^*,\theta,\lambda),\lambda(X))}.\]
The subbundle $[\ft_0]$ is defined by
\[[\ft_0]=\{X\in[\fc]\ |\ \Lambda(X)=X=\Theta(X)\}.\]
By definition, the conjugate linear involutions $\Theta$ and $\Lambda$ commute. Thus, we also obtain a complex linear involution $\sigma$ which is the complex linear extension of a Cartan involution of the split real form $\fg_0.$ 
If $\fg_0=\fh\oplus\fm$ is the corresponding Cartan decomposition, then the eigenbundle splitting of $\underline\fg$ is given by 
	\[\sigma=\Theta\circ\Lambda:[\fh_\C]\oplus[\fm_\C].\]
Recall that for $x=\haf\sum\limits_{\alpha\in\Delta^+}H_\alpha,$ and if the highest root has height $m_\ell$ then we defined $g_{+}=exp(\frac{2\pi ix}{m_\ell+1}).$ The $\nabla^c$-parallel automorphism $\Xx^+$ is defined by
\[\Xx^+((\Delta^+\subset\fc^*,\theta,\lambda),X))=((\Delta^+\subset\fc^*,\theta,\lambda),g_{+}(X)).\]
\end{proof}

The following proposition and corollary are proven in section 4 of \cite{cyclicSurfacesRank2}, the proofs are analogous to Proposition \ref{MapstoSymmSpaceProp}.
 \begin{Proposition}\label{MapstoTriplesProp}
 	Let $N$ be a smooth simply connected manifold and $\fg$ be a complex simple Lie algebra. Let $(\widetilde\fg,\widetilde D)\ra N$ be a flat $\fg$-Lie algebra bundle with the following
 	\begin{itemize}
 		\item  A smoothly varying Hitchin triple $(\widetilde\fc,\widetilde{\fn^+},\widetilde\Theta,\widetilde\Lambda)$ with corresponding decompositions \[\widetilde\fg=\widetilde\ft_0\oplus i\widetilde\ft\oplus\widetilde\fa\oplus\widetilde\fn^+\oplus\widetilde\fn^-\ =\ \widetilde\ft_0\oplus\widetilde\fm\]
 		\item $\widetilde\nabla$ a connection so that $(\widetilde\fc,\widetilde{\fn^+},\widetilde\Theta,\widetilde\Lambda)$ is parallel. 
 		\item $\widetilde\omega\in\Omega^1(N,\widetilde\fm)$ with $\widetilde\nabla+ad_{\widetilde\omega}=\widetilde D.$
 	\end{itemize}
 	 Then there is a map $f:N\ra\sG/\sT_0,$ unique up to post composition by an element of $\sG,$ so that 
 	\[(\widetilde\fg,\ \widetilde\nabla,\ \widetilde\omega,\ \widetilde\fc,\ \widetilde\fn^+,\ \widetilde\Theta,\ \widetilde\Lambda)=f^*(\underline\fg,\ \nabla^c,\ \omega,\ [\fc],\ [\fn^+],\ \Theta,\ \Lambda).\]
 \end{Proposition}
 \begin{Corollary}\label{MapstoTriplesCor}
 	Let $N$ be a smooth manifold and $\widetilde\fg\ra N$ be a flat $\fg$-Lie algebra bundle equipped with the structure of Proposition \ref{MapstoTriplesProp}, then there exists
 	\begin{enumerate}
 	\item A representation $\rho:\pi_1(N)\ra \sG$ unique up to conjugation 
 	\item A $\rho$-equivariant map $f$ from the universal cover $\widetilde N$ of $N$ to the space of Hitchin triples $\sG/\sT_0$ satisfying the conclusion of Proposition \ref{MapstoTriplesProp}.
  	\end{enumerate}
 \end{Corollary}

\begin{Remark}
	Lemma \ref{propertiesOfTriplesLem}, Proposition \ref{MapstoTriplesProp} and Corollary \ref{MapstoTriplesCor} all have analogous versions for the space of Cartan triples. 
\end{Remark}
\begin{Proposition}\label{HiggsMaptoTrip}
	The $\fsp(4,\C)$-Lie algebra bundle associated to a $\sSp(4,\R)$-Higgs bundle  \[\left(N\oplus N^{-1}K,\mtrx{\nu&0\\0&\mu},\mtrx{0&1\\1&0}\right)\] satisfies the hypotheses of Corollary \ref{MapstoTriplesCor}. 
\end{Proposition}
\begin{proof}
For any $\sSp(4,\R)$-Higgs bundle $(V,\beta,\gamma)$, there is an involution $\widetilde\sigma,$ defined by \eqref{PTDSpictureSp4}, which gives a decomposition
 \[\widetilde{\fsp(4,\C)}=\fgl(V)\oplus ( Sym^2(V)\oplus Sym^2(V^*))\subset End(V\oplus V^*).\] 
Since the Higgs bundle $(V,\beta,\gamma)$ is a polystable, there exists a metric $H$ solving the Hitchin equations which is compatible with the involutions $\sigma$. 
This metric gives a decomposition of the bundle
\[\widetilde{\fsp(4,\C)}=\widetilde{\fsp(4)}\oplus i\widetilde{\fsp(4)}\]
which is induced by decomposition of the $End(V\oplus V^*)$ into self adjoint and skew adjoint endomorphisms.
Such a splitting defines a smoothly varying Cartan involution $\widetilde\Theta=\mtrx{1&0\\0&-1}$. 
Since $(V,\beta,\gamma)$ is a $\sSp(4,\R)$-Higgs bundle, the metric $H$ and the involution $\widetilde\sigma$ are compatible with a smoothly varying split real form
\[\widetilde\Lambda=\widetilde\Theta\circ\widetilde\sigma=\widetilde\sigma\circ\widetilde\Theta\]
 Furthermore, the involutions $\widetilde\Theta$ and $\widetilde\Lambda$ are covariantly constant with respect to the Chern connection $\nabla_H.$

By Remark \ref{MaxSp4dLiealgbundleRemark}, the $\fsp(4,\C)$-Lie algebra bundle $\widetilde{\fsp(4,\C)}$ associated to $\left(N\oplus N^{-1}K,\mtrx{\nu&q_2\\q_2&\mu},\mtrx{0&1\\1&0}\right)$ has 
\[\widetilde\fc\oplus\widetilde{\fn^+}=(\Oo\otimes\fc)\oplus (N^{2}K^{-2}\otimes\fg_{-2L_1}\oplus K^{1}\otimes\fg_{L_1+L_2})\oplus (N^{-2}K^{1}\otimes\fg_{-L_1+L_2})\oplus (N^{2}\otimes\fg_{2L_2}).\]
In general, these subbundles are not parallel with respect to $\nabla_H$. 
However, if the quadratic differential in the Higgs field vanishes, by Proposition \ref{MetricSplittingProposition} the metric is diagonal, and thus the bundles $\widetilde\fc$ and $\widetilde{\fn^+}$ are parallel with respect to the Chern connection. 
Furthermore, the adjoint $\phi^{*_H}=H^{-1}\overline\phi^TH=-\Theta(\phi)$ does not contain any components in $\widetilde\fc$. Thus,
\[\xymatrix{\widetilde\omega=\phi+\phi^{*_H}=\Omega^1(\Sigma,\widetilde\fm)}\] 
with $\nabla_H(\widetilde\omega)=0$ and $\nabla_H+ad_{\widetilde\omega}$ flat.
\end{proof}
 \begin{Corollary}\label{liftingtoTripCor}
 	Let $\rho\in\Rr_d(\sSp(4,\R)$ for $g-1\leq d\leq 3g-3.$ For each choice of conformal structure $(S,J)=\Sigma$ so that the Higgs bundle corresponding to $\rho$ is given by $\left(N\oplus N^{-1}K,\mtrx{\nu&0\\0&\mu},\mtrx{0&1\\1&0}\right)$ with $\mu\neq 0,$ there is unique a $\rho$-equivariant map $f:\widetilde\Sigma\ra\sG/\sT_0$ lifting the harmonic metric. 
 \end{Corollary} 
 \begin{proof}
 	 By Proposition \ref{HiggsMaptoTrip}, to each metric solving the Hitchin equations, there is a unique $\rho$-equivariant map $f:\widetilde\Sigma\ra\sG/\sT_0$ which lifts the $\rho$-equivariant harmonic metric $h_\rho:\widetilde\Sigma\ra\sG/\sK$. 
 	 If $\mu\neq0,$ by Remark \ref{StabilityofMdRemark}, the Higgs bundle $\left(N\oplus N^{-1}K,\mtrx{\nu&0\\0&\mu},\mtrx{0&1\\1&0}\right)$ is stable.
 	  Thus, there is a unique metric solving the Hitchin equations.
 \end{proof}
\bigskip
\section{Cyclic surfaces}\label{cyclicSurf}\smallskip
The surfaces we will be interersted are solutions to certain Pfaffian systems in the spaces of Cartan triples and Hitchin triples. The cyclic surfaces defined below are more general than \cite{cyclicSurfacesRank2}, yet, we show that deformations of the below cyclic surfaces have many similarities with deformations of Labourie's cyclic surfaces. 

\subsection{Cyclic Pfaffian systems and cyclic surfaces} 
The general Pfaffian system defintions in this section come from section 7 of \cite{cyclicSurfacesRank2}. 
\begin{Definition}\label{PfaffDef}
	Let $E\ra N$ be a vector bundle over a smooth manifold $N,$ and $(\eta_1,\cdots,\eta_n)$ be a collection of differential forms on $N$ valued in $E.$ A submanifold $L\subset N$ is called a solution to the \em Pfaffian system defined by $(\eta_1,\cdots,\eta_n)$ \em if $\eta_j|_L\equiv0$ for all $j.$
\end{Definition}
The Pfaffian systems we will be interested are defined as follows:
\begin{Definition}
  	\label{PfaffGcyclicDef} Let $\omega\in\Omega^1(\sG/\sT,\underline\fg)$ be the Maurer Cartan form of the space of Cartan triples $\sG/\sT.$ A \em $\sG$-cyclic Pfaffian system \em is defined by the vanishing of the following $\underline\fg$-valued forms
  	\[((\widehat\omega_0,\widehat\omega_2,\cdots,\widehat\omega_{m_\ell-1})\ ,\ [\widehat\omega_{-1},\widehat\omega_{-1}]\ ,\ \widehat\omega_{-1}+\Theta(\widehat\omega_1))\]
  	where $\omega=\sum \widehat\omega_j$ is the decomposition of \eqref{ZmellMCdecom}.
  \end{Definition} 
 For the space of Hitchin triples, we define a $\sG_0$-cyclic Pfaffian system as follows. 
   \begin{Definition}
  	\label{PfaffG0cyclicDef} Let $\omega\in\Omega^1(\sG/\sT_0,\underline\fg)$ be the Maurer Cartan form of the space of Hitchin triples $\sG/\sT_0.$ The \em $\sG_0$-cyclic Pfaffian system \em is defined by the vanishing of the following $\underline\fg$-valued forms
  	\[((\widehat\omega_0,\widehat\omega_2,\cdots,\widehat\omega_{m_\ell-1})\ ,\ [\widehat\omega_{-1},\widehat\omega_{-1}]\ ,\ \widehat\omega+\Theta(\widehat\omega)\ ,\ \Lambda(\omega)-\omega)\]
  	where $\omega=\sum \widehat\omega_j$ is the decomposition of \eqref{ZmellMCdecom}.
  \end{Definition}
The above definition are related to the $\tau$-primitive maps consider for compact groups $\sG$ in \cite{minimalsurfacesaffinetoda}. 
In the context of representations of surface groups, we are interested in maps from a Riemann surface $\Sigma$ to the spaces of Cartan triples and Hitchin triples. 

\begin{Definition}\label{GG0cyclicsurfaceDef}
	Let $\Sigma$ be a Riemann surface (not necessarily compact), a map $f:\Sigma\ra\sG/\sT$ is a \em $\sG$-cyclic surface \em if it is a $\sG$-cyclic Pfaffian system and $f^*\widehat\omega_{-1}$ is a $(1,0)$-form. Similarly, a map $f:\Sigma\ra\sG/\sT_0$ is a \em$\sG_0$-cyclic surface \em if it is a $\sG_0$-cyclic Pfaffian system and $f^*\widehat\omega_{-1}$ is a $(1,0)$-form.
\end{Definition}
\begin{Remark}
	The reality condition $f^*(\Lambda(\omega))=f^*(\omega)$ for a $\sG_0$-cyclic surface implies $f(\Sigma)$ lies in a $\sG_0$ orbit.
If $\sG_0$ is a split real form of Hodge type, then $T_0=T,$ and the $\sG_0$-cyclic condition is just an extra symmetry a $\sG$-cyclic map must satisfy. 
\end{Remark}

Maximal $\sSp(4,\R)$-Higgs bundles of the form $\left(N\oplus N^{-1}K,\mtrx{\nu&0\\0&\mu},\mtrx{0&1\\1&0}\right)$ with $\mu\neq 0$ give rise to special equivariant $\sSp(4,\R)$-cyclic surfaces:
\begin{Proposition}\label{cyclicSurfacesforRd} Let $\rho\in\Rr_d(\sSp(4,\R))$ for $g-1\leq d\leq 3g-3.$ If $(S,J)=\Sigma$ is a conformal structure so that the Higgs bundle corresponding to $\rho$ is given by $\left(N\oplus N^{-1}K,\mtrx{\nu&0\\0&\mu},\mtrx{0&1\\1&0}\right)$ with $\mu\neq 0,$ then the unique $\rho$-equivariant map \[f:\widetilde\Sigma\ra\sSp(4,\R)/(\sU(1)\times\sU(1))\]
 from Corollary \ref{liftingtoTripCor} is a $\sSp(4,\R)$-cyclic surface. Furthermore, $f^*\omega_{-\alpha_i}\neq0$ for each simple root $\alpha_i$ and there exists a simple root $\beta$ such that $f^*\omega_{-\beta}$ nowhere vanishing. 
\end{Proposition}
\begin{proof}
By Corollary \ref{liftingtoTripCor}, there is a unique equivariant lift of the harmonic metric $f:\widetilde\Sigma\ra\sG/\sT_0,$ it remains to show that this map is cyclic. 
The pullback of the Maurer-Cartan form is given by $\phi+\phi^{*_H}$ where $\phi=\mtrx{&\beta\\\gamma&}$ is the Higgs field and $H$ is the metric solving the Hitchin equations. 
By construction, see Remark \ref{MaxSp4dLiealgbundleRemark}, the Higgs field $\phi$ is valued in $\Omega^{10}(\Sigma,f^*[\widehat{\fsp(4,\C)}_{-1}]),$ and thus $\phi=f^*\widehat\omega_{-1}.$
Also, $\phi^{*_H}=-\widetilde\Theta(\phi)\in\Omega^{01}(\Sigma,f^*[\widehat{\fsp(4,\C)}_{1}]).$ 
It follows that $f^*\widehat\omega_j=0$ for $j\neq \pm1$ and $f^*(\Lambda\omega)=f^*\omega.$ 
The Higgs field is given by
\[f^*\widehat\omega_{-1}=\phi=1\otimes \fg_{-L_1-L_2}+\mu\otimes\fg_{2L_1}+q_2\otimes\fg_{L_1+L_2}+\nu\otimes\fg_{2L_2}\]
where $L_1+L_2$ and $-2L_1$ are the two simple roots (see Remark \ref{MaxSp4dLiealgbundleRemark}). 
Since $f^*\omega_{-L_1-L_2}$ is nowhere vanishing and $\mu\neq0,$ we have $f^*\omega_{-\alpha_i}\neq0$ for each simple root $\alpha_i$ and $f^*\omega_{-\beta}$ is nowhere vanishing for a simple root $\beta.$
\end{proof}

The following theorem relates equivariant cyclic surfaces and Higgs bundles that are fixed points of $\gen{\zeta\dwn{m_\ell}}\subset\sU(1).$  

\begin{Theorem}\label{GCyclicSurfacesAndHiggsBundlesTheorem}
	Let $\sG$ be a complex simple Lie group of rank at least 2, and $\rho\in\Rr(\sG).$ If $\underline{\fg}\ra\sG/\sT$ is the associated Lie algebra bundle and $f:\widetilde \Sigma\ra\sG/\sT$ be a $\rho$-equivariant $\sG$-cyclic surface, then  $(f^*\underline{\fg},(f^*\nabla^c)^{01},f^*\widehat\omega_{-1})$ is a $\sG$-Higgs bundle that is a fixed point of the $\gen{\zeta\dwn{m_\ell+1}}$-action. 
	Furthermore, $f^*B_\Theta$ solves the Hitchin equations \eqref{Hit} which simplify to 
	\[F_{f^*\nabla^c}+\sum\limits_{i=1}^\ell [f^*\omega_{\alpha_i},f^*\omega_{-\alpha_i}]+[f^*\omega_{\mu},f^*\omega_{-\mu}]=0.\] 
\end{Theorem}
\begin{proof}
To prove that $(f^*\underline\fg,(f^*\nabla^c)^{01},f^*\widehat\omega_{-1})$ is a $\sG$-Higgs bundle we just need to show $f^*\widehat\omega_{-1}$ is holomorphic. 
By equations \eqref{ZmellMCdecom}, the flatness equations for $\nabla^c+\omega$ we have 
\[d^{\nabla}\widehat\omega_{-1}+\sum\limits_{j=0}^{m_\ell}[\widehat\omega_j,\widehat\omega_{-j-1}]=0.\]
By the cyclic assumption, $f^*\widehat\omega_j=0$ for $j\neq \pm 1,$ thus, pulling back the flatness equations, we have
\[d^{f^*\nabla}(f^*\widehat\omega_{-1})=0.\]
Since $f^*\widehat\omega_{-1}$ is a $(1,0)$-form, we conclude that $(d^{f^*\nabla})^{01}f^*\omega_{-1}=0.$ 

To see that it is a fixed point of $\gen{\zeta\dwn{m_\ell+1}},$ recall from Lemma \ref{propertiesOfTriplesLem} that there is an automorphism $\Xx_+:\underline\fg\ra\underline\fg,$ of order $(m_\ell+1),$ which acts as $\zeta_{m_\ell+1}^{-1}$ on $[\widehat\fg_{-1}].$
Thus $f^*(\Xx_+)^{-1}$ is a gauge transformation of $f^*\underline\fg$ which acts as $\zeta_{m_\ell+1}$ on the Higgs field $f^*\widehat\omega_{-1}.$

Recall that by definition of a $\sG$-cyclic surface, we have $f^*(-\Theta(\widehat\omega_{-1}))=f^*\widehat\omega_1,$ thus the adjoint of the Higgs field $f^*\widehat\omega_{-1}$ is given by $f^*(-\Theta(\widehat\omega_{-1}))=f^*\widehat\omega_{1}.$
Using the decompositions of \eqref{gradingequivs} we have
\[\xymatrix{\widehat\omega_{1}=\sum\limits_{i=1}^\ell \omega_{\alpha_i}+\omega_{-\mu}&\widehat\omega_{-1}=\sum\limits_{i=1}^\ell \omega_{-\alpha_i}+\omega_{\mu}}.\] 
The assumption that $f^*\omega_0=0,$ and the flatness equations of \eqref{ZmellMCdecom} imply
\[F_{f^*\nabla^c}+\sum\limits_{i=1}^\ell [f^*\omega_{\alpha_i},f^*\omega_{-\alpha_i}]+[f^*\omega_{\mu},f^*\omega_{-\mu}]=0.\]
Since $f^*\nabla^c$ is a metric connection for the hermitian metric $f^*B_\Theta,$ and the holomorphic structure is on $f^*\fg$ is defined to be $(f^*\nabla^c)^{01},$ we conclude that $f^*B_\Theta$ solves the Hitchin equations. 
\end{proof}
\begin{Corollary} \label{GCyclicSurfacesAndHiggsBundlesCor}
Let $\sG$ be a complex simple Lie group with rank at least 2, $\rho\in\Rr(\sG),$ and $f:\widetilde \Sigma\ra\sG/\sT$ be an $\rho$-equivariant $\sG$-cyclic surface. If $\pi:\sG/\sT\ra\sG/\sK$ is the natural projection, then the associated equivariant harmonic map $h_{\rho,J}=f\circ\pi:\widetilde\Sigma\ra\sG/\sH$ is a minimal surface.
\end{Corollary}
\begin{proof}
	Since the Higgs bundle admits a solution to the Hitchin equations, it is polystable. Since it is a fixed point of $\gen{\zeta_{m_\ell+1}}$ and $rank(\fg)\geq2,$ the quadratic differential is the image of the Hitchin fibration vanishes, thus the Hopf differential of the harmonic map is zero and we conclude the harmonic map is a branched minimal immersion. 
\end{proof}
Similarly, for $\sG_0$-cyclic surfaces we have the following theorem.
\begin{Theorem}\label{G0CyclicSurfacesAndHiggsBundlesTheorem}
	Let $\sG$ be a complex simple Lie group of rank at least 2, and $\rho\in\Rr(\sG).$ If $\underline{\fg}\ra\sG/\sT_0$ is the associated Lie algebra bundle and $f:\widetilde \Sigma\ra\sG/\sT_0$ be a $\rho$-equivariant $\sG_0$-cyclic surface. 
	If $\underline{\fg}\ra\sG/\sT_0$ is the associated Lie algebra bundle, then $(f^*[\fh_\C],(f^*\nabla^c)^{01},f^*\widehat\omega_{-1})$ is a $\sG_0$-Higgs bundle that is a fixed point of the $\gen{\zeta\dwn{m_\ell+1}}$-action. Furthermore, $f^*B_\Theta$ solves the Hitchin equations \eqref{Hit} which simplify to 
	\[F_{f^*\nabla^c}+\sum\limits_{i=1}^\ell [f^*\omega_{-\alpha_i},f^*\omega_{\alpha_i}]+[f^*\omega_{\mu},f^*\omega_{-\mu}]=0.\] 
\end{Theorem}
\begin{Remark}
	In this case, the representation $\rho\in\Rr(\sG)$ is actually in $\Rr(\sG_0).$
\end{Remark}
\begin{proof}
	Recall from Lemma \ref{propertiesOfTriplesLem}, that the Lie algebra bundle $\underline\fg\ra\sG/\sT_0$ has a complex linear involution $\sigma=\Theta\circ\Lambda$ which has eigenbundle decomposition $\underline\fg=[\fh_\C]\oplus[\fm_\C]$ where $\fg_0=\fh\oplus\fm$ is the corresponding Cartan decomposition. 
	To show that $(f^*[\fh_\C],(f^*\nabla^c)^{01},f^*\widehat\omega_{-1})$ is a $\sG_0$-Higgs bundle, we must show that $f^*\widehat\omega_{-1}\in\Omega^{10}(\Sigma, f^*[\fm_\C]).$ 
	Recall from Remark \ref{PTDSpictureSp4}, the involution $\sigma$ preserves the height grading
	\[\fg_{-m_\ell}\oplus\cdots\fg_{-1}\oplus\fc\oplus\fg_1\oplus\cdots\oplus\fg_{m_\ell},\]
	and thus, $\sigma$ preserves both $\widehat\fg_1$ and $\widehat\fg_{-1}.$
	By the definition of a $\sG_0$-cyclic surface \ref{GG0cyclicsurfaceDef}, $f^*\omega=f^*\widehat\omega_{-1}+f^*\widehat\omega_1$ and
	 \[\xymatrix{f^*\Theta(f^*\widehat\omega_{-1}+f^*\widehat\omega_1)=-f^*\widehat\omega_{1}-f^*\widehat\omega_{-1}&f^*\Lambda(f^*\widehat\omega_{-1}+f^*\widehat\omega_1)=f^*\widehat\omega_{-1}+f^*\widehat\omega_1}.\]
	Hence, $\sigma(f^*\widehat\omega_{-1}+f^*\widehat\omega_1)=-(f^*\widehat\omega_{-1}+f^*\widehat\omega_1),$ and, furthermore, since $\sigma$ preserves $[\widehat\fg_{-1}],$
	 \[\sigma(\widehat\omega_{-1})=-\widehat\omega_{-1}.\] 
	This proves that $f^*\widehat\omega_{-1}\in\Omega^{10}(\Sigma,f^*[\fm_\C]).$ 

	We also need to check that the gauge transformation $f^*\Xx_+$ is an $f^*[\fh_\C]$-gauge transformation.
	 Recall that the grading element $x$ of the PTDS is in the $+1$-eigenspace of $\sigma$ (see Remark \ref{PTDSpictureSp4}).
	Since $\Xx_+$ is obtained from exponentating $x,$ it follows that $f^*\Xx_+$ is an $[\fh_\C]$-gauge transformation.
	The proof of the rest of the theorem is identical to the proof of Theorem \ref{GCyclicSurfacesAndHiggsBundlesTheorem}.
	 \end{proof}


\subsection{Deformations of cyclic Pfaffian systems and cyclic surfaces}
\begin{Definition}\label{1storderDefPfaff}
Let $F=(f_t):L\ra N$ be a one parameter family with $f_0$ being the inclusion and set 
\[\xi=\left.\frac{d}{dt}\right|_{t=0}f_t.\]
Then $\xi\in\Omega^0(L,f_0^*TN)$ is a vector field along $L$ in $N$ called the \em tangent vector field to the family \em $F.$ 
A family $F=(f_t)$ is a \em first order deformation of the Pfaffian system $L$ defined by $(\eta_1,\cdots,\eta_n)$ \em if, for all $j,$ 
	\[\left.\frac{d}{dt}\right|_{t=0}f^*_t\eta_i=0.\] 
\end{Definition}
In the above definition, we have chosen a connection to identify $f^*_tE$ and $f^*_0E,$ this choice does not effect the definition. 
\begin{Definition}\label{VectFieldalongPfaffDef}
A vector field $\xi$ along a solution $L$ of a Pfaffian system given by $\eta=(\eta_1,\cdots,\eta_n)$ is an \em infinitesimal variation of the Pfaffian system \em if, for any connection $\nabla,$ and all $j,$
\[\left.\iota_\xi d^{\nabla}\eta_j\right|_{L}=\left.-d^\nabla(\iota_\xi \eta_j)\right|_L.\]
\end{Definition}
The relation between first order deformations and variations is given by Proposition 7.1.4 of \cite{cyclicSurfacesRank2}:
\begin{Proposition}\label{DefsofPfaffProp}
	Let $\xi$ be a tangent vector to a family of first order deformations of a Pfaffian system $\eta=(\eta_1,\cdots,\eta_n)$, then $\xi$ is an infinitesimal variation of the Pfaffian system. 
 \end{Proposition}
 \begin{proof}
 	See Proposition 7.1.4 of \cite{cyclicSurfacesRank2}.
 \end{proof}

\begin{Definition}\label{G0cyclicdeformation}
An \em infinitesimal variation of a $\sG$-cyclic surface \em is an infinitesimal variation of a $\sG$-cyclic Pfaffian system. An \em infinitesimal variation of a $\sG_0$-cyclic surface \em is an infinitesimal variation $\xi$ of a $\sG_0$-cyclic Pfaffian system such that $\Lambda(\xi)=\xi.$
\end{Definition}

\begin{Definition}
Let $\rho:\pi_1(S)\ra\sG$ be a representation and $f:\widetilde \Sigma\ra \sG/\sT$ be a $\rho$-equivariant $\sG$-cyclic surface. If $\xi$ is an infinitesimal variation of a $\sG$-cyclic surface, then $\xi$ is an \em infinitesimal variation of the equivariant $\sG$-cyclic surface \em if it is $\rho$-equivariant. 
Similarly for an equivariant $\sG_0$-cyclic surface. 
\end{Definition} 

The signs in the following lemma will be crucial.
\begin{Lemma}\label{signLemma}
	Let $\Sigma$ be a compact Riemann surface and $f$ a $\sG$-cyclic surface or a $\sG/\sT_0$-cyclic surface. Let $\alpha\in\Omega^{10}(\Sigma,f^*\underline\fg)$ and $\beta\in\Omega^{01}(\Sigma,f^*\underline\fg)$ then 
	\[\xymatrix{-i\int\limits_\Sigma B_\fg(\alpha,\Theta\alpha)\geq 0&\text{and}&i\int\limits_\Sigma B_\fg(\beta,\Theta\beta)\geq 0}.\]
	Also, if $\alpha,\beta\in\Omega^1(\Sigma,f^*\underline\fg)$ and $\gamma\in\Omega^0(\Sigma,f^*\underline\fg),$ then 
	\begin{equation}\label{invariancesign}
		B_\fg(\gamma,[\beta,\alpha])=B_\fg([\gamma,\alpha],\beta).
	\end{equation}
\end{Lemma}
\begin{proof}
It is suffices to check the sign on a form $\alpha=A\cdot a$ where $a$ is a section of $f^*\underline\fg$ and $A\in\Omega^{10}(\Sigma).$ By Remark \ref{Thetaforms},
$\Theta(\alpha)=\overline A\cdot\Theta(a)$ and hence, since $-B_\fg(\cdot,\Theta\cdot)$ is positive definite,
\[-i\int\limits_\Sigma B_\fg(\alpha,\Theta\alpha)=-i\int\limits_\Sigma A\wedge\overline A\cdot B_\fg(a,\Theta a)\geq0.\]
Equation \eqref{invariancesign} follows from a calculation using invariance of the Killing form. 
\end{proof}
Let $f:\Sigma\ra\sG/\sT$ be a $\sG$-cyclic surface, we will use the following notation
\[\xymatrix{\left.\widehat\omega_{-1}\right|_{f(\Sigma)}=\Phi=\Phi_{-1}+\Phi_{m_\ell}&\text{and}&\left.\widehat\omega_{1}\right|_{f(\Sigma)}=\Phi^*=\Phi^*_{1}+\Phi^*_{-m_\ell}.}\]
Let $\xi$ is an infinitesimal variation of $f$, and denote the contraction with $\omega$ by
\[\zeta=\iota_\xi(\omega).\]
Using the various decompositions of the Maurer Cartan form $\omega,$ we have the following decompositions 
\begin{equation}\label{zetadecomps}
\xymatrix{\zeta=\zeta_0+\sum\limits_{\alpha\in\Delta(\fc,\fg)}\zeta_\alpha\ ,&\zeta=\sum\limits_{j=-m_\ell}^{m_\ell} \zeta_j\ ,& \zeta=\sum\limits_{j\in\Z/(m_\ell+1)\Z}\widehat\zeta_j}.\end{equation}
The following notation will also be useful
\begin{equation}
	\label{Ynotation}\zeta=\widehat\zeta_{-1}+\widehat\zeta_{0}+\widehat\zeta_{1}+\widehat\zeta_{Y}, 
\end{equation}
where $\widehat\zeta_Y=\sum\limits_{j\neq 0,1,-1}\widehat\zeta_j.$ 

Using the decomposition of the flatness equations \eqref{ZgradEq} we have
\[d^{\nabla^c}\omega_j+\sum\limits_{k=-m_\ell}^{m_\ell}[\omega_k,\omega_{j-k}]=0 .\]
By Definition \ref{VectFieldalongPfaffDef}, on the surface $f(\Sigma),$ we have $\iota_\xi(d^{\nabla^c}\omega_j)=-d^{\nabla^c}(\zeta_j)$ for $j\neq -m_\ell,-1,1,m_\ell.$ Contracting the wedge product is given by 
\[\iota_\xi[\omega_j,\omega_{j-k}]\ =\ [\zeta_j,\omega_{j-k}]-[\omega_j,\zeta_{j-k}]\ =\ [\zeta_j,\omega_{j-k}]+[\zeta_{j-k},\omega_j].\]
 Thus, contracting the flatness equations with $\xi$ yeilds
\begin{eqnarray}\label{jneq1-1}
	d^{\nabla^c}(\zeta_j)=\sum\limits_{k=-m_\ell}^{m_\ell}([\zeta_{k}, \omega_{j-k}]+[\zeta_{j-k},\omega_k])&\ \ \ \ \ \ j\neq -m_\ell,-1,1,m_\ell
\end{eqnarray}  
The assumption on a cyclic surface that $f^*\widehat\omega_j=0$ for $j\neq\pm1$ and the fact that $\Phi$ is a $(1,0)$-form and $\Phi^*$ is a $(0,1)$-form allows us to simplify the equations. For $1<j<m_\ell$ we have
\begin{equation}
	\label{PosJparital}
	\xymatrix{\p^{\nabla^c}(\zeta_j)=2([\zeta_{j+1},\Phi_{-1}]+[\zeta_{j-m_\ell},\Phi_{m_\ell})]&\text{and}&\bar\p^{\nabla^c}(\zeta_j)=2[\zeta_{j-1},\Phi_{1}^*]}
\end{equation}
and for $-m_\ell<j<-1$ we have
\begin{equation}
\label{NegJpartial}
	\xymatrix{\p^{\nabla^c}(\zeta_{j})=2[\zeta_{j+1},\Phi_{-1}]&\text{and}&\bar\p^{\nabla^c}(\zeta_j)=2([\zeta_{j-1},\Phi^*_1]+[\zeta_{j+m_\ell},\Phi^*_{-m_\ell}])}
\end{equation}
Let $\pi_Y$ denote the projection onto the $Y$ component of equation \eqref{Ynotation}, then equations \eqref{PosJparital} and \eqref{NegJpartial} can be written compactly as:
\begin{equation}
	\label{pYand1}
	\xymatrix@R=0em{
	\p^{\nabla^c}\widehat\zeta_Y=2\pi_Y\left([\Phi,(\widehat\zeta_Y+\widehat\zeta_{-1})]\right)&\text{and}&
 \bar\p^{\nabla^c}\widehat\zeta_Y=2\pi_Y\left([\Phi^*,(\widehat\zeta_Y+\widehat\zeta_{1})]\right)}
\end{equation}
For $j=-1,1$ even though $\iota_\xi \left.(d^{\nabla^c}\widehat\omega_j)\right|_\Sigma\neq-d^{\nabla^c}\widehat\zeta_j,$ by equations \eqref{ZmellgradingEq} we have
\begin{equation}
	\label{bp1,p-1}
	\xymatrix{\left.\iota_\xi(\p^{\nabla^c}\widehat\omega_{-1})\right|_\Sigma= -2[\zeta_0,\Phi]&\text{and}& \left.\iota_\xi(\bar\p^{\nabla^c}\widehat\omega_{1})\right|_\Sigma= -2[\zeta_0,\Phi^*] }
\end{equation}
Similarly, 
\begin{equation}
	\label{p1,-bp1}
	\xymatrix@R=0em{\left.\iota_\xi(\bar\p^{\nabla^c}\omega_{-1})\right|_\Sigma=-2([\zeta_{-2},\Phi^*_1]+[\zeta_{-1+m_\ell},\Phi^*_{-m_\ell}])\ , & \left.\iota_\xi(\p^{\nabla^c}\omega_{1})\right|_\Sigma= -2([\zeta_2,\Phi_{-1}]+[\zeta_{-1-m_\ell},\Phi_{m_\ell}])\ , \\
	\left.\iota_\xi(\bar\p^{\nabla^c}\omega_{m_\ell})\right|_\Sigma=-2[\zeta_{m_\ell-1},\Phi^*_1]\ , & \left.\iota_\xi(\p^{\nabla^c}\omega_{-m_\ell})\right|_\Sigma= -2[\zeta_{-m_\ell+1},\Phi_{-1}]\ .}
\end{equation}

\begin{Proposition}\label{secondDerivativeProp}
	The second derivatives are given by 

\begin{equation}
	\label{typeppzetaY}
	\xymatrix{ \bar\p^{\nabla^c}(\p^{\nabla^c}\widehat\zeta_Y)=4\pi_Y\left(\left[\left[\widehat\zeta_Y,\Phi^*\right],\Phi\right]\right)\ , &
	\p^{\nabla^c}(\bar\p^{\nabla^c}\widehat\zeta_Y)=4\pi_Y\left(\left[\left[\widehat\zeta_Y,\Phi\right],\Phi^*\right]\right)},
\end{equation}
\begin{equation}
	\label{typeppzeta0}
	\xymatrix{ \bar\p^{\nabla^c}(\p^{\nabla^c}\widehat\zeta_0)=4\pi_{i\ft}\left(\left[\left[\widehat\zeta_0,\Phi^*\right],\Phi\right]\right)\ , &
	\p^{\nabla^c}(\bar\p^{\nabla^c}\widehat\zeta_0)=4\pi_{i\ft}\left(\left[\left[\widehat\zeta_0,\Phi\right],\Phi^*\right]\right)}.
\end{equation}
\end{Proposition}
\begin{proof}
Recall that on a $\sG$-cyclic surface we have $\p^{\nabla^c}\Phi=\bar \p^{\nabla^c}\Phi=\p^{\nabla^c}\Phi^*=\bar\p^{\nabla^c}\Phi^*=0$. We will first show equation \eqref{typeppzetaY}. Using equations \eqref{PosJparital} and \eqref{NegJpartial}, a direct computation shows
	\begin{equation*}
\begin{dcases}
	\bar\p^{\nabla^c}(\p^{\nabla^c}\zeta_j)=4\left(\left[\left[\zeta_j,\Phi_1^*\right],\Phi_{-1}\right]+\left[\left[\zeta_{j-m_\ell-1},\Phi_1^*\right],\Phi_{m_\ell}\right]+\left[\left[\zeta_j,\Phi^*_{-m_\ell}\right],\Phi_{m_\ell}\right]\right)&1<j<m_\ell-1 \\
	 \bar\p^{\nabla^c}(\p^{\nabla^c}\zeta_j)=4\left(\left[\left[\zeta_j,\Phi^*_{1}\right],\Phi_{-1}\right]+\left[\left[\zeta_{j+1+m_\ell},\Phi^*_{-m_\ell}\right],\Phi_{-1}\right]\right)&-m_\ell<j<-2\\
	\p^{\nabla^c}(\bar\p^{\nabla^c}\zeta_j)=4\left(\left[\left[\zeta_j,\Phi_{-1}\right],\Phi^*_1\right]+\left[\left[\zeta_{j-1-m_\ell},\Phi_{m_\ell}\right],\Phi^*_{1}\right]\right)& 2<j<m_\ell\\
	\p^{\nabla^c}(\bar\p^{\nabla^c}\zeta_j)=4\left(\left[\left[\zeta_j,\Phi_{-1}\right],\Phi_{1}^*\right]+\left[\left[\zeta_{j+m_\ell+1},\Phi_{-1}\right],\Phi^*_{-m_\ell}\right]+\left[\left[\zeta_j,\Phi{m_\ell}\right],\Phi_{-m_\ell}^*\right]\right)& -m_\ell+1<j<-1
	\end{dcases}
\end{equation*}
The remaining cases are given by 
\[\xymatrix@R=0em{\bar\p^{\nabla^c}(\p^{\nabla^c}\zeta_{m_\ell-1})=2\left(\left[\bar\p^{\nabla^c}\zeta_{m_\ell},\Phi_{-1}\right]+\left[\bar\p^{\nabla^c}\zeta_{-1},\Phi_{m_\ell}\right]\right)\ ,& \bar\p^{\nabla^c}(\p^{\nabla^c}\zeta_{-2})=2[\bar\p^{\nabla^c}\zeta_{-1},\Phi_{-1}]\ ,
\\\p^{\nabla^c}(\bar\p^{\nabla^c}\zeta_{-m_\ell+1})=2\left(\left[\p^{\nabla^c}\zeta_{-m_\ell},\Phi_{1}^*\right]+\left[\p^{\nabla^c}\zeta_{1},\Phi_{-m_\ell}^*\right]\right)\ ,&\p^{\nabla^c}(\bar\p^{\nabla^c}\zeta_2)=2[\p^{\nabla^c}\zeta_{1},\Phi_{1}^*]\ .
}
	\]
We will compute the first two cases, the remaining two cases follow by a symmetric argument.
Recall that $\p^{\nabla^c}\Phi=\bar\p^{\nabla^c}\Phi=\p^{\nabla^c}\Phi^*=\bar\p^{\nabla^c}\Phi^*=0$ on a cyclic surface, thus 
\begin{equation}\label{omegawedgeomegaeq}
\xymatrix@R=0em{\bar\p^{\nabla^c}\left(\iota_\xi[\omega_{m_\ell},\omega_{-1}]\right)= [\bar\p^{\nabla^c}\zeta_{m_\ell},\Phi_{-1}]+ [\bar\p^{\nabla^c}\zeta_{-1},\Phi_{m_\ell}]\ ,&\bar\p^{\nabla^c}\left(\iota_\xi[\omega_{-1},\omega_{-1}]\right)=2[\bar\p^{\nabla^c}\zeta_{-1},\Phi_{-1}]}
\end{equation}
However, since $[\omega_1,\omega_1]=[\omega_1,\omega_{-m_\ell}]=[\omega_{-1},\omega_{-1}]=[\omega_{-1},\omega_{m_\ell}]=0$ on a cyclic surface, we have
\eqtns{
	\bar\p^{\nabla^c}\left(\iota_\xi\left[\omega_{m_\ell},\omega_{-1}\right]\right)=-\left[\left.\iota_\xi(\bar\p^{\nabla^c}\omega_{m_\ell})\right|_\Sigma,\Phi_{-1}\right]-\left[\left.\iota_\xi(\bar\p^{\nabla^c}\omega_{-1})\right|_\Sigma,\Phi_{m_\ell}\right]\\
	\bar\p^{\nabla^c}\left(\iota_\xi\left[\omega_{-1},\omega_{-1}\right]\right)=-2\left[\left.\iota_\xi(\bar\p^{\nabla^c}\omega_{-1})\right|_\Sigma,\Phi_{-1}\right]
}{phiwedgephivanisheq}
Using equations \eqref{p1,-bp1} and \eqref{omegawedgeomegaeq}, we have the desired result:
\[2([\bar\p^{\nabla^c}\zeta_{m_\ell},\Phi_{-1}]+ [\bar\p^{\nabla^c}\zeta_{-1},\Phi_{m_\ell}])= 4\left(\left[\left[\zeta_{m_\ell-1},\Phi_1^*\right],\Phi_{-1}\right]+\left[\left[\zeta_{-2},\Phi_2^*\right],\Phi_{m_\ell}\right]+\left[\left[\zeta_{-1+m_\ell},\Phi^*_{-m_\ell}\right],\Phi_{m_\ell}\right]\right) \]
and 
\[2\left[\bar\p^{\nabla^c}\zeta_{-1},\Phi_{-1}\right]=4\left(\left[\left[\zeta_{-2},\Phi_1^*\right],\Phi_{-1}\right]+\left[\left[\zeta_{-1+m_\ell},\Phi_{-m_\ell}^*\right],\Phi_{-1}\right]\right).\]
Thus, we obtain the desired formula:
\[\xymatrix{ \bar\p^{\nabla^c}(\p^{\nabla^c}\widehat\zeta_Y)=4\pi_Y\left(\left[\left[\widehat\zeta_Y,\Phi^*\right],\Phi\right]\right)&\text{and} &
	\p^{\nabla^c}(\bar\p^{\nabla^c}\widehat\zeta_Y)=4\pi_Y\left(\left[\left[\widehat\zeta_Y,\Phi\right],\Phi^*\right]\right)}\]

We now prove formula \eqref{typeppzeta0}, for $\bar\p^{\nabla^c}\p^{\nabla^c}\widehat\zeta_0$ and $\p^{\nabla^c}\bar\p^{\nabla^c}\widehat\zeta_0.$ Since $\widehat\omega_0$ vanishes along a $\sG$-cyclic surface, by the flatness equations \eqref{jneq1-1}, we have
\[\xymatrix{\p^{\nabla^c}\widehat\zeta_0=2[\widehat\zeta_1,\Phi]&\text{and}&\bar\p^{\nabla^c}\widehat\zeta_0=2[\widehat\zeta_{-1},\Phi^*].}\]
Recall that $\widehat\zeta_0$ vanishes along the subbundle $[\ft],$ that is, $\pi_{i\ft}\widehat\zeta_0=\widehat\zeta_0.$
Thus, the second derivatives are 
\begin{equation}\label{pp0eq}
	\xymatrix{\bar\p^{\nabla^c}\p^{\nabla^c}\widehat\zeta_0=2[\bar\p^{\nabla^c}\widehat\zeta_1,\Phi]=2\pi_{i\ft}\left([\bar\p^{\nabla^c}\widehat\zeta_1,\Phi]\right)\ ,& \p^{\nabla^c}\bar\p^{\nabla^c}\widehat\zeta_0=2[\p^{\nabla^c}\widehat\zeta_{-1},\Phi^*]=2\pi_{i\ft}\left([\p^{\nabla^c}\widehat\zeta_{-1},\Phi^*]\right) }.
\end{equation}
Since $\widehat\omega_1=-\Theta(\widehat\omega_{-1})$ on a $\sG$-cyclic surface, it follows that $\pi_{i\ft}([\widehat\omega_1,\widehat\omega_{-1}])=0$ along a $\sG$-cyclic surface. Thus,
\[\left.\iota_\xi d^{\nabla^c}\pi_{i\ft}([\widehat\omega_1,\widehat\omega_{-1}])\right|_\Sigma=-d^{\nabla^c}\left.(\iota_\xi \pi_{i\ft}([\widehat\omega_1,\widehat\omega_{-1}])\right|_\Sigma).\]
The subbundle $[i\ft]$ is parallel with respect to $\nabla^c,$ thus
\[\pi_{i\ft}\left(\left.\iota_\xi d^{\nabla^c}([\widehat\omega_1,\widehat\omega_{-1}])\right|_\Sigma\right)=-\pi_{i\ft}\left(d^{\nabla^c}\left.(\iota_\xi ([\widehat\omega_1,\widehat\omega_{-1}])\right|_\Sigma)\right).\] 
For the (1,0) part, we have 
\[\pi_{i\ft}\left(\left.\iota_\xi \p^{\nabla^c}([\widehat\omega_1,\widehat\omega_{-1}])\right|_\Sigma\right)=\pi_{i\ft}\left([\left.(\p^{\nabla^c}\widehat\omega_1)\right|_\Sigma,\Phi]+[\left.(\p^{\nabla^c}\widehat\omega_{-1})\right|_{\Sigma},\Phi^*]\right)=\pi_{i\ft}\left(\left.(\p^{\nabla^c}\widehat\omega_{-1})\right|_{\Sigma},\Phi^*])\right).\]
The term $[\p^{\nabla^c}\widehat\zeta_1,\Phi]$ vanishes since it is a $(2,0)$-form. A similar calculations for the (0,1) part gives 
\[\xymatrix{
	\pi_{i\ft}\left(\left.\iota_\xi\left(\bar\p^{\nabla^c}([\widehat\omega_1,\widehat\omega_{-1}])\right)\right|_\Sigma\right)=
	\pi_{i\ft}\left([\bar\p^{\nabla^c}\widehat\zeta_{1},\Phi] \right).
}\]
Thus, by equations \eqref{pp0eq},
\[\xymatrix{\bar\p^{\nabla^c}\p^{\nabla^c}\widehat\zeta_0=2\pi_{i\ft}\left(\left.\iota_\xi\left(\bar\p^{\nabla^c}([\widehat\omega_1,\widehat\omega_{-1}])\right)\right|_\Sigma\right)&\text{and}&\p^{\nabla^c}\bar\p^{\nabla^c}\widehat\zeta_0}= 2\pi_{i\ft}\left(\left.\iota_\xi \p^{\nabla^c}([\widehat\omega_1,\widehat\omega_{-1}])\right|_\Sigma\right) \]
The term $-\pi_{i\ft}\left(d^{\nabla^c}\left.(\iota_\xi ([\widehat\omega_1,\widehat\omega_{-1}])\right|_\Sigma)\right)$ is computed using equation \eqref{bp1,p-1}:
\[-\pi_{i\ft}\left(\p^{\nabla^c}\left.(\iota_\xi ([\widehat\omega_1,\widehat\omega_{-1}])\right|_\Sigma)\right)=-\pi_{i\ft}\left([\left.\left(\p^{\nabla^c}\iota_\xi\widehat\omega_1\right)\right|_\Sigma,\Phi]+[\left.\left(\p^{\nabla^c}\iota_\xi\widehat\omega_{-1}\right)\right|_\Sigma,\Phi^*]
\right)=2\pi_{i\ft}\left([[\widehat\zeta_0,\Phi],\Phi^*]\right). \]
A similar computation shows 
\[-\pi_{i\ft}\left(\bar\p^{\nabla^c}\left.(\iota_\xi ([\widehat\omega_1,\widehat\omega_{-1}])\right|_\Sigma)\right)=2\pi_{i\ft}\left([[\widehat\zeta_0,\Phi^*],\Phi]\right).\]
Thus, on a $\sG$-cyclic surface,

\[	\xymatrix{\p^{\nabla^c}\bar\p^{\nabla^c}\widehat\zeta_0=4\pi_{i\ft}\left([[\widehat\zeta_0,\Phi],\Phi^*]\right)&\text{and}&\bar\p^{\nabla^c}\p^{\nabla^c}\widehat\zeta_0}= 4\pi_{i\ft}\left([[\widehat\zeta_0,\Phi^*],\Phi]\right). \]
\end{proof}

\begin{Proposition}\label{zetaYzeta0flatProp}
		Let $\rho:\pi_1(S)\ra\sG$ and $f:\widetilde\Sigma\ra\sG/\sT$ be a $\rho$-equivariant $\sG$-cyclic surface. Let $\widehat\zeta_Y,$ $\widehat\zeta_0,$ $\Phi$ and $\Phi^*$ be as above, then 
		\[\xymatrix{\p^{\nabla^c}\widehat\zeta_Y=0\ ,&\bar\p^{\nabla^c}\widehat\zeta_Y=0\ ,&[\Phi,\widehat\zeta_Y]=0\ ,&[\Phi^*,\widehat\zeta_Y]=0}\]
		and
		\[\xymatrix{\p^{\nabla^c}\widehat\zeta_0=0\ ,&\bar\p^{\nabla^c}\widehat\zeta_0=0\ ,&[\Phi,\widehat\zeta_0]=0\ ,&[\Phi^*,\widehat\zeta_0]=0}.\]
\end{Proposition}
\begin{proof} Recall from Lemma \ref{signLemma} that
	\[0\leq-i\int\limits_{\Sigma}B_\fg\left(\p^{\nabla^c}\widehat\zeta_Y,\Theta\left(\p^{\nabla^c}\widehat\zeta_Y\right)\right) .\]
Since the canonical connection is a metric connection, we have 
\[d\left(-B_\fg\left(\widehat\zeta_Y,\Theta\left(\p^{\nabla^c}\widehat\zeta_Y\right)\right)\right)=-B_\fg\left(\p^{\nabla^c}\widehat\zeta_Y,\Theta\left(\p^{\nabla^c}\widehat\zeta_Y\right)\right)-B_\fg\left(\widehat\zeta_Y,\Theta\left(\bar\p^{\nabla^c}\p^{\nabla^c}\widehat\zeta_Y\right)\right).\]
Integrating over $\Sigma$ gives 
\[0\leq -i\int\limits_{\Sigma}B_\fg\left(\p^{\nabla^c}\widehat\zeta_Y,\Theta\left(\p^{\nabla^c}\widehat\zeta_Y\right)\right)=i\int\limits_\Sigma B_\fg\left(\widehat\zeta_Y,\Theta\left(\bar\p^{\nabla^c}\p^{\nabla^c}\widehat\zeta_Y\right)\right).\]
Recall that $[\widehat\fg_Y]=\bigoplus\limits_{j\neq -1,0,1}[\widehat\fg_j]$ and $\Theta\left([\widehat\fg_j]\right)\subset[\widehat\fg_{-j}].$ Also, if $i+j\neq0\ mod\ \left(m_\ell+1\right)$ then $\widehat\fg_j$ and $\widehat\fg_i$ are orthogonal with respect to $B_\Theta.$ Thus, the bundles $[\widehat\fg_j]$ and $[\widehat\fg_i]$ are orthogonal. 
Thus, using equations \eqref{typeppzetaY} we have 
\[0\leq 4i\int\limits_\Sigma B_\fg\left(\widehat\zeta_Y,\Theta\left(\pi_Y\left(\left[\left[\widehat\zeta_Y,\Phi^*\right],\Phi\right]\right)\right)\right)=4i\int\limits_\Sigma B_\fg\left(\widehat\zeta_Y,\Theta\left(\left[\left[\widehat\zeta_Y,\Phi^*\right],\Phi\right]\right)\right).\]
 Lemma \ref{signLemma} and the cyclic surface assumption $\Phi=-\Theta\left(\Phi^*\right)$ yield 
\[0\leq-i\int\limits_{\Sigma}B_\fg\left(\p^{\nabla^c}\widehat\zeta_Y,\Theta\left(\p^{\nabla^c}\widehat\zeta_Y\right)\right)=-4i\int\limits_\Sigma B_\fg\left(\left[\widehat\zeta_Y,\Phi^*\right],\Theta\left(\left[\widehat\zeta_Y,\Phi^*\right]\right)\right)\leq 0.\]
Thus 
\begin{equation}\label{flatzetaY1}
	\xymatrix{\p^{\nabla^c}\widehat\zeta_Y=0&\text{and}&
	[\Phi^*,\widehat\zeta_Y]=0}
\end{equation}
By a symmetric argument, we obtain
\begin{equation}\label{flatzetaY2}
	\xymatrix{\bar\p^{\nabla^c}\widehat\zeta_Y=0&\text{and}&
	[\Phi,\widehat\zeta_Y]=0}
\end{equation}
For $\widehat\zeta_0,$ consider the following integral:
\[0\leq -i\int_\Sigma B_\fg\left(\p^{\nabla^c}\widehat\zeta_0,\Theta(\p^{\nabla^c}\widehat\zeta_0)\right)=i\int_\Sigma B_\fg\left(\widehat\zeta_0,\Theta(\bar\p^{\nabla^c}\p^{\nabla^c}\widehat\zeta_0)\right)\]
Using equations \eqref{typeppzeta0}, the fact that $i\ft\oplus\ft$ is is orthogonal, Lemma \ref{signLemma} and $\Theta(\Phi^*)=-\Phi$ we have 
\[0\leq4i\int_\Sigma B_\fg\left(\widehat\zeta_0,\Theta([[\widehat\zeta_0,\Phi^*],\Phi])\right)=-4i\int_\Sigma B_\fg\left([\widehat\zeta_0,\Phi^*],\Theta\left([\widehat\zeta_0,\Phi^*]\right)\right)\leq0.\]
Thus, 
\[\xymatrix{\p^{\nabla^c}\widehat\zeta_0=0&\text{and}&[\Phi^*,\widehat\zeta_0]=0}.\]
A symmetric argument shows 
\[\xymatrix{\bar\p^{\nabla^c}\widehat\zeta_0=0&\text{and}&[\Phi,\widehat\zeta_0]=0}.\]
\end{proof}
The same calculations show that the analogous proposition for equivariant $\sG_0$-cyclic surfaces is also true.
\begin{Corollary}\label{G0zeta10covconstCor}
		Let $\rho:\pi_1(S)\ra\sG$ and $f:\widetilde\Sigma\ra\sG/\sT_0$ be a $\rho$-equivariant $\sG_0$-cyclic surface. Let $\widehat\zeta_Y,$ $\Phi$ and $\Phi^*$ be as above, then  
		\[\xymatrix{\p^{\nabla^c}\widehat\zeta_Y=0\ ,&\bar\p^{\nabla^c}\widehat\zeta_Y=0\ ,&[\Phi,\widehat\zeta_Y]=0\ ,&[\Phi^*,\widehat\zeta_Y]=0\ ,}\]
		\[\xymatrix{\p^{\nabla^c}\widehat\zeta_0=0\ ,&\bar\p^{\nabla^c}\widehat\zeta_0=0\ ,&[\Phi,\widehat\zeta_0]=0\ ,&[\Phi^*,\widehat\zeta_0]=0}.\]	
		Furthermore, if $\fc=\ft\oplus i\ft\oplus\fa\oplus i\fa$ is the decomposition of the Cartan subalgebra, then $\widehat\zeta_0$ vanishes along $i\ft\oplus i\fa.$ In particular, if $\sG_0$ is of Hodge type then $\widehat\zeta_0=0$
\end{Corollary}
\begin{proof}
	The first part is an immediate corollary of the proof of Proposition \ref{zetaYzeta0flatProp}. The variation $\widehat\zeta_0$ is along $i\ft\oplus\fa\oplus i\fa,$ where $\Lambda$ acts as $+1$ on $\fa$ and $-1$ on $i\ft\oplus i\fa.$  
	But, by the reality condition of variations of $\sG_0$-cyclic surfaces, $\Lambda(\widehat\zeta_0)=\widehat\zeta_0;$ thus, $\widehat\zeta_0$ vanishes along $i\ft\oplus i\fa.$ Recall that a $\sG_0$ is of Hodge type then the $\fa=\{0\},$ thus, in this case $\widehat\zeta_0=0.$
\end{proof}

If $\rho:\pi_1(S)\ra\sG$ is representation and $f:\widetilde\Sigma\ra\sG/\sT$ is a $\sG$-cyclic surface, then Proposition \ref{zetaYzeta0flatProp} says that $\widehat\zeta_0$ and $\widehat\zeta_Y$ are covariantly constant with respect to the \em flat connection \em $f^*\nabla^c+\Phi+\Phi^*.$ 
Thus, if either $\widehat\zeta_0$ or $\widehat\zeta_Y$ is non zero, then they are in the centralizer of the representation $\rho.$ However, if the centralizer of $\rho$ is discrete then the centralizing subalgebra is zero, thus we have the following proposition. 

\begin{Proposition}\label{irrRepRigid}
	Let $\sG$ be a complex simple Lie group, and $\rho:\pi_1(S)\ra\sG$ be an irreducible representation. If $f:\widetilde\Sigma\ra\sG/\sT$ be a $\rho$-equivariant cyclic surface, then for any variation $\xi,$ we have 
	\[\xymatrix{\iota_\xi(f^*\widehat\omega_0)=\widehat\zeta_0=0&\text{and}&\iota_\xi(f^*\widehat\omega_Y)=\widehat\zeta_Y=0}.\]
\end{Proposition}


	

\bigskip
\section{Special cyclic surfaces and the proof of Theorem \ref{uniqueminThm}}\label{proofofTHM}\smallskip
In this section we consider equivariant cyclic surfaces with extra conditions on $f^*\omega_{-1}$ and show that for these special equivariant cyclic surfaces are rigid. We then relate these special equivariant cyclic surfaces to Higgs bundles and specifically maximal $\sSp(4,\R)$-Higgs bundles. 

\begin{Proposition}\label{RigidProposition}
	Let $(S,J)=\Sigma$ be a compact Riemann surface, $\sG$ be a complex simple Lie group of rank at least 2 and not $\sSL(3,\C).$  Let $\rho:\pi_1(S)\ra\sG$ be an irreducible representation and $f:\widetilde\Sigma\ra\sG/\sT_0$ be a $\rho$-equivariant $\sG_0$-cyclic surface so that $f^*\omega_{-\alpha_i}\neq0$ for all simple roots $\alpha_i$. If $\xi$ is an infinitesimal variation with the property that there exists a simple root $\alpha$ with $\iota_\xi\omega_{-\alpha}\equiv 0,$ then \[\iota_\xi\omega\equiv0.\] 
\end{Proposition}
\begin{Remark}
	The analogous statement follows for $\sG$-cyclic surfaces if one assumes that there are simple roots $\alpha$ and $\beta$ so that $\iota_\xi\omega_{-\alpha}\equiv 0\equiv\iota_\xi\omega_{+\beta}.$ 
	For $\sG_0$-cyclic surfaces, if $\iota_\xi\omega_{-\alpha}\equiv 0,$ the reality condition $\Lambda\xi=\xi$ on an infinitesimal variation implies that $\iota_\xi\omega_{\Lambda(-\alpha)}\equiv 0.$ Furthermore, since $\Theta$ flips positive simple roots and negative simple roots, $\sigma$ preserves the set of positive simple roots and $\Lambda=\Theta\circ\sigma$, it follows that $\Lambda(-\alpha)$ is a positive simple root. If $\sG_0$ is of Hodge type, then $\Lambda(-\alpha)=\alpha.$
 \end{Remark}
\begin{proof}
Let $\xi$ be a variation of the $\rho$-equivariant $\sG_0$-cyclic surface $f:\widetilde\Sigma\ra\sG/\sT_0,$ and $\zeta=\iota_\xi\omega.$
Using the decompositions of \eqref{zetadecomps} and \eqref{Ynotation}, by Corollary \ref{G0zeta10covconstCor},
\[\xymatrix{\widehat\zeta_0=0&\text{and}&\widehat\zeta_Y=0}.\]
It remains to show $\widehat\zeta_1=0=\widehat\zeta_{-1}.$ Recall that $\sG\neq \sSL(3,\C),$ thus, $\widehat\fg_{Y}\neq\{0\},$ in particular $\fg_{\pm2}\neq \{0\}.$ 
A infinitesimal variation $\xi$ of a $\sG_0$-cyclic surface satisfies the reality condition $\Lambda\xi=\xi.$ 
By Lemma \ref{fcCartanInvLemma}, $\Theta(\fg_\alpha)=\fg_{-\alpha}$ for all roots, and by Proposition \ref{Tmaxcompact}, the involution $\sigma$ sends roots simple root spaces to simple root spaces. 
Since there is a simple root $\alpha$ so that $\zeta_{-\alpha}\equiv0$ and $\Lambda\zeta=\zeta,$ it follows that there is a simple root $-\Lambda(\alpha)$ so that $\zeta_{-\Lambda\alpha}\equiv0.$ 

By equation \eqref{pYand1}, we have 
 	\[\xymatrix{0=\p^{\nabla^c}(\zeta_{-2})=2[\zeta_{-1}, \Phi_{-1}]&\text{and}&0=\bar\p^{\nabla^c}(\zeta_{2})=2[\zeta_{1},\Phi_1^*]}.\]
 	Thus for each pair of simple roots $\alpha_i,\alpha_j$ so that $\alpha_i+\alpha_j$ is a root, we have 
 	\[\xymatrix{[\zeta_{-\alpha_i},\Phi_{-\alpha_j}]+[\zeta_{-\alpha_j},\Phi_{-\alpha_i}]=0&\text{and}&[\zeta_{\alpha_i},\Phi^*_{\alpha_j}]+[\zeta_{\alpha_j},\Phi^*_{\alpha_i}]=0}.
 		\]
 	Since $\Phi_{-\alpha_i}=f^*\omega_{-\alpha_i},$ by assumption $\Phi_{-\alpha_i}$ is a nonzero holomorphic section. By the definition of a $\sG_0$-cyclic surface, 
 	\[f^*\Theta(\Phi_{\alpha_i})=f^*(\Theta\omega_{-\alpha_i})=-f^*(\omega_{\alpha_i})=-\Phi^*_{\alpha_i}.\]
 	Thus, $\Phi_{\alpha_i}^*$ is also nonzero for all simple roots. 
 	
 	The group $\sG$ is simple, thus the Dynkin diagram is connected and we conclude that $\zeta_{\pm\alpha_i}=0$ for all simple roots. It remains to show that for the highest root $\mu,$ we have $\zeta_{\pm m_\ell}=\zeta_{\pm\mu}=0.$ 
 	By equations \eqref{pYand1}, we have 
 	\[\xymatrix{0=\bar\p^{\nabla^c}(\zeta_{m_\ell-1})=2[\zeta_{m_\ell},\Phi_{1}^*]&\text{and}&0=\p^{\nabla^c}(\zeta_{-m_\ell+1})=2[\zeta_{m_\ell},\Phi_{-1}]}.\]
 	Since $\sG\neq\sSL(3,\C),$ we have $\fg_{\pm1}\neq\fg_{\pm (m_\ell-1)}\neq \{0\}.$ 
 	Thus, for each roots $\gamma=\mu-\alpha_i\in\fg_{ m_\ell-1}$ we have $0=[\zeta_{\mu},\Phi_{-\alpha_i}].$ 
 	Hence $\zeta_\mu=0,$ and similarly, $\zeta_{-\mu}=0.$
\end{proof}
\begin{Remark}
Proposition \ref{RigidProposition} is also true when $\sG=\sSL(3,\C)$, see Proposition 7.7.4 of \cite{cyclicSurfacesRank2}.
\end{Remark}

In \cite{cyclicSurfacesRank2}, Labourie considers maps $f:S\ra\sG/\sT_0$ from the surface $S,$ without a conformal structure, to the space of Hitchin Triples $\sG/\sT_0$ that satisfy $f^*\widehat\omega_j=0$ for $j\neq\pm1,$  $f^*(\Theta(\widehat\omega_{-1}))=-f^*(\widehat\omega_1),$ $f^*([{\widehat\omega_{-1},\widehat\omega_{-1}}])=0,$ $f^*(\Lambda\omega)=f^*\omega$ and satisfy the extra assumption that \em for all \em simple roots $\alpha_i,$
	 \[f^*\omega_{-\alpha_i}\ \ \ \text{is nowhere vanishing}.\]
It is then proven that there is a unique conformal structure on $S$ so that $f^*\widehat\omega_{-1}$ is a $(1,0)$-form. 

\begin{Proposition}
Let $rank(\sG)\geq 2,$ a map $f:S\ra\sG/\sT$ satisfies: $f^*\widehat\omega_j=0$ for $j\neq\pm1$ and
\[\xymatrix{f^*(\Theta(\widehat\omega_{-1}))=-f^*(\widehat\omega_1)\ ,&f^*([\widehat\omega_{-1},\widehat\omega_{-1}])=0\ ,&f^*(\Lambda\omega)=f^*\omega.}\]
Suppose that $f^*\omega_{\alpha_i}$ has discrete zeros for all simple roots $\alpha_i$ and that there exists a simple root $\beta$ so that  
	\[f^*\omega_{-\beta}\ \ \ \text{is nowhere vanishing}.\]
	 Then there exists a unique conformal structure $(S,J)=\Sigma,$ so that $f:\Sigma\ra\sG/\sT$ is a cyclic surface.
\end{Proposition}
Thus, Definition \ref{GG0cyclicsurfaceDef} and the cyclic surfaces in Proposition \ref{RigidProposition} are generalizations of the cyclic surfaces in \cite{cyclicSurfacesRank2}. 
The cyclic surfaces related to maximal $\sSp(4,\R)$ representations are more special than those considered in Proposition \ref{RigidProposition} and more general than Labourie's. 
Namely, we only require that there exists a simple root $\alpha_i$ so that $f^*\omega_{-\alpha_i}$ is nowhere vanishing.
\begin{proof}
Let $\beta\in\Delta^+(\fg,\fc)$ be a simple roots for which $f^*\omega_{-\beta}$ is nowhere vanishing. 
Since $f^*\omega$ is nowhere vanishing, $df:TS\ra[\fg_{-\beta}]$ is an isomorphism. 
Thus, there is a unique complex structure $(S,J)=\Sigma$ so that $f^*\omega_{-\beta}$ is a $(1,0)$-form.

Since $f^*([\widehat\omega_{-1},\widehat\omega_{-1}])=0,$ decomposing this in terms of root spaces we have for all simple roots $\alpha$ and $\gamma$
\[\xymatrix{[f^*\omega_{-\alpha}, f^*\omega_{-\gamma}]=0&\text{and}& [f^*\omega_{-\alpha}, f^*\omega_{\mu}]=0}.\]
Recall that $\fg$ is simple, so there is a simple root $\alpha$ so that $-\beta-\alpha$ is a root, in particular, 
\[[[\fg_{-\alpha}],[\fg_{-\beta}]]\neq 0.\]
By $[f^*\omega_{-\alpha}, f^*\omega_{-\beta}]=0,$ it follows that $f^*\omega_{-\alpha}$ is a $(1,0)$-form. 
Using the fact that $\fg$ is simple and that $f^*\omega_{-\alpha}$ has discrete zeros, we conclude that for all simple roots $\alpha,$ the form $f^*\omega_{-\alpha}$ is a $(1,0)$-form. 
Similary, there is there is a simple root $\alpha$ so that $\mu-\alpha$ is a root. 
We again conclude that $f^*\omega_{\mu}$ is a $(1,0)$-form, proving $f^*\widehat\omega_{-1}$ is a $(1,0)$-form.
\end{proof}

Putting everything together, we obtain the following theorem which is the analogue to the transversality of the Hitchin map in \cite{cyclicSurfacesRank2}. 

\begin{Theorem}\label{RigidityOfCyclicSurf}
	Let $\sG$ be a complex simple Lie group of rank at least 2, $\rho:\pi_1(S)\ra\sG$ an irreducible representation, and $(S,J)=\Sigma$ be a conformal structure. Suppose $f:\widetilde\Sigma\ra\sG/\sT_0$ a $\rho$-equivariant $\sG_0$-cyclic surface such that there exists a simple root $\alpha$ so that $f^*\omega_{-\alpha}$ is nowhere vanishing and, for all simple roots $\alpha_i,$ the form $f^*\omega_{-\alpha_i}$ is nonzero. Let $(\rho_t,J_t)$ is a one parameter family with $(\rho_0,J_0)=(\rho,J)$ and $f_t:\widetilde{(S,J_t)}\ra\sG/\sT_0$ be a family of $\rho_t$-equivariant $\sG_0$-cyclic surfaces with $f_0=f$. If $[\left.\frac{d}{dt}\right|_{t=0}\rho_t]=0,$ then $\left.\frac{d}{dt}\right|_{t=0}J_t=0.$ 
\end{Theorem}

\begin{proof}
	Let $\rho\in\Rr(\sG)$ be an irreducible representation and let $(S,J)=\Sigma$ be a conformal structure. 
	Let $f:\widetilde\Sigma\ra\sG/\sT_0$ be a $\rho$-equivariant $\sG_0$-cyclic surface so that there is a simple root $\alpha$ with $f^*\omega_{-\alpha}$ nowhere vanishing and, for all simple roots $\alpha_i,$ the form $f^*\omega_{-\alpha_i}\neq0.$ 
	Suppose $(\rho_t,J_t)$ is a one parameter family and $f_t:\widetilde{(S,J_t)}\ra\sG/\sT_0$ is a family of $\rho_t$-equivariant $\sG_0$-cyclic surfaces with $f_0=f$, that is for all $\gamma\in\pi_1(S),$
	\[f_t(\gamma(s))=\rho_t(\gamma)\cdot f_t(s).\]
	If $[\left.\frac{d}{dt}\right|_{t=0}\rho_t]=[\rho],$ then, since $\rho$ is irreducible, the tangent space a $\rho$ is given by 
	\[T_\rho\Rr(\sG)=H_\rho^1(S,\fg).\]
	Thus, after conjugating the family $\rho_t$ by a family of elements of $\sG,$ and preforming a similar transformation for $f_t$, for all $\gamma\in\pi_1(S)$ we have
	\[\left.\frac{d}{dt}\right|_{t=0}f_t(\gamma(s))=\rho(\gamma)\cdot\left.\frac{d}{dt}\right|_{t=0}f_t(s).\]
	In particular, $\xi(s)=\left.\frac{d}{dt}\right|_{t=0}f_t(s)$ is an $\rho$-equivariant infinitesimal deformation of $f.$ 
	Since $f^*\omega_{-\alpha}$ is nowhere vanishing, $f^*\omega_{-\alpha}:T\Sigma\ra[\fg_{-\alpha}]$ is a bijection. 
	Let $X$ be the vector field along $\Sigma$ so that $\iota_{\omega_{-\alpha}}\xi=f^*\omega_{-\alpha(X)},$ then $df(X)$ is an infinitesimal variation of $f.$
By construction, $\xi-df(X)$ is an equivariant infinitesimal variation of $f$ which vanishes along the simple root space $[\fg_\alpha].$ Thus, by Proposition \ref{RigidProposition}, $\xi-df(X)=0.$ 

To see that $\left.\frac{d}{dt}\right|_{t=0}J_t=0,$ we employ an argument of Marco Spinaci \cite{MarcoWorkshopNotes}. We have 
	\[\xi=\left.\frac{\p f_t}{\p t}\right|_{t=0}=df(X),\] thus
	\[\zeta=\omega(df_0(X))=\Phi(X)+\Phi^*(X).\]
In particular, $\zeta$ is self adjoint and hence lives in the subbundle $[i\fk].$ Also, $\widehat\zeta_{-1}=\Phi(X)$ is holomorphic and $\widehat\zeta_1=\Phi^*(X)$ is antiholomorphic. 

Let $\Psi_t=f_t^*\omega=\Phi_t+\Phi_t,$ by definition, for all tangent vectors, we have
\begin{equation}\label{eqPsiJt}
	\Psi_t(J_tv)=i\Phi_t(v)-i\Phi^*_t(v)=\left(i\pi_{\widehat\fg_{-1}}-i\pi_{\widehat\fg_{1}}\right)\Psi_t(v).
\end{equation} 
Recall that, for vector fields $Y$ and $Y'$ on $\sG/\sT_0,$ we have
\begin{equation}
	\label{torsioneq}\omega(d^{\nabla}_YY')=\omega(d^{\nabla}_{Y'}Y)+\omega([Y,Y'])+\omega(T(Y,Y'))
\end{equation}
where $T(Y,Y')$ is the torsion tensor given by Lemma \ref{CanonicalConnectionandTorsion}. Differentiating equation \eqref{eqPsiJt} yields
\[\left.d^{f_t^*\nabla^c}_{\frac{\p}{\p t}}\left(f^*_t\omega(J_tv)\right)\right|_{t=0}=\left.\left(i\pi_{\widehat\fg_{-1}}-i\pi_{\widehat\fg_{1}}\right)d^{f_t^*\nabla^c}_{\frac{\p}{\p t}}\left(f^*_t\omega(v)\right)\right|_{t=0}.\]
Using the pullback of equation \eqref{torsioneq} by $f_t,$ the left hand side of the above equations is given by
\[\left.\left(d^{f_t^*\nabla^c}_{J_tv}\left(f^*_t\omega\left(\frac{\p}{\p t}\right)\right)+f^*_t\omega\left[\frac{\p}{\p t},J_tv\right]+f^*_t\omega \left(T\left(\frac{\p}{\p t},J_tv\right)\right)\right)\right|_{t=0}.\] 
The expression for the torsion in Lemma \ref{CanonicalConnectionandTorsion} and the decomposition $\fg=\ft_0\oplus\fm$ imply the above expression can be rewritten as
\[\left.\left(d^{f_t^*\nabla^c}_{J_tv}\left(\omega\left(\frac{\p f_t}{\p t}\right)\right)+\omega\left(df_t\left(\frac{\p J_t}{\p t}v\right)\right)+\pi_{[\fm]}\left(\left[f^*_t\omega\left(\frac{\p}{\p t}\right) , f^*_t\omega (J_tv)\right]\right)\right)\right|_{t=0}.\]
Using $\Psi_0=f^*_0\omega$ and $\left.\zeta=\frac{\p f_t}{\p t}\right|_0,$ evaluating at $t=0$ yields
\[d^{f_0^*\nabla^c}_{J_0v}\left(\zeta\right)+\Psi_0\left(\left.\frac{\p J_t}{\p t}\right|_0 v \right)+\pi_{[\fm]}\left(\Psi_0(J_0v),\zeta \right).\]
Since $\left[\frac{\p}{\p t},v\right]=0,$ a similar computation shows that the left hand side of equation \eqref{eqPsiJt} is 
\[d^{f_0^*\nabla^c}_{v}\left(\zeta\right)+\pi_{[\fm]}\left(\Psi_0(v),\zeta \right).\]
Thus, we have
\begin{equation}
	\label{simplifiedtderveq}
	d^{f_0^*\nabla^c}_{J_0v}\left(\zeta\right)+\Psi_0\left(\left.\frac{\p J_t}{\p t}\right|_0 v \right)+\pi_{[\fm]}\left(\Psi_0(J_0v),\zeta \right)=\left(i\pi_{\widehat\fg_{-1}}-i\pi_{\widehat\fg_{1}}\right)\left( d^{f_0^*\nabla^c}_{v}\left(\zeta\right)+\pi_{[\fm]}\left(\Psi_0(v),\zeta \right)\right)
\end{equation}
Recall $\zeta=\omega\left(\left.\frac{\p f_t}{\p t}\right|_0\right)$ is in $\omega(df(T\Sigma))\subset[i\fk].$ Also, $\pi_{[\fm]}$ and $\left(i\pi_{\widehat\fg_{-1}}-i\pi_{\widehat\fg_{1}}\right)$ commute with the Cartan involution $\Theta.$ Thus, we can consider the $[i\fk]$ part of equation \eqref{simplifiedtderveq}. 
This yields
\begin{equation}
	d^{f_0^*\nabla^c}_{J_0v}\left(\zeta\right)+\Psi_0\left(\left.\frac{\p J_t}{\p t}\right|_0 v \right)=\left(i\pi_{\widehat\fg_{-1}}-i\pi_{\widehat\fg_{1}}\right)\left( d^{f_0^*\nabla^c}_{v}\left(\zeta\right)\right).
\end{equation}
Rearranging the equations and using the fact that $\widehat\zeta_{-1}$ is holomorphic and $\widehat\zeta_1$ is antiholomorphic gives
\[\Psi_0\left(\left.\dfrac{\p J_t}{\p t}\right|_0 v \right)=2i(\bar\p\widehat\zeta_{-1}-\p\widehat\zeta_1)=0.\]
Since $\Psi$ is injective, it follows that $\left.\dfrac{\p J_t}{\p t}\right|_{t=0}=0,$ as desired.
 \end{proof}

For maximal representations $\rho\in\Rr_d(\sSp(4,\R))$ with $g-1\leq d\leq3g-3$ we obtain local uniqueness of the conformal structures $J_\rho$ in which the $\rho$-equivariant harmonic map is minimal.
\begin{Theorem}\label{LocalUniqueminThm}
 	Let $\rho\in\Rr_d(\sSp(4,\R))$ for $g-1< d\leq 3g-3$ or $\rho\in\Rr_{g-1}(\sSp(4,\R))$ and $\rho$ irreducible. Then the collection of conformal structures $\{J_\rho\}$ so that the $\rho$-equivariant harmonic mapping $\widetilde\Sigma\ra \sSp(4,\R)/\sU(2)$ is a minimal immersion is nonempty and discrete. 
 \end{Theorem}
\begin{proof}
Fix a representation $\rho\in \Rr_d(\sSp(4,\R)),$ by Corollary \ref{existenceCor} and Lemma \ref{RdMinimalImmersion}, there exists a conformal structure $(S,J_\rho)=\Sigma_\rho$ in which the corresponding $\rho$-equivariant harmonic map $f_\rho:\tilde\Sigma_\rho\ra \sSp(4,\R)/\sU(2)$ is a minimal immersion.
Since $\rho$ is irreducible, in any such conformal structure $J_\rho$ the corresponding Higgs bundle is stable and hence  given by 
\[\left(N\oplus N^{-1}K,\mtrx{\nu&0\\0&\mu},\mtrx{0&1\\1&0}\right)\]
with $\mu\neq0.$ 

By Proposition \ref{cyclicSurfacesforRd}, to such a Higgs bundle there is a $\rho$-equivariant $\sSp(4,\R)$-cyclic surface 
\[f_\rho:\widetilde\Sigma_{\rho}\ra\sSp(4,\C)/(\sU(1)\times\sU(1))\]
so that $f_\rho^*\omega_{-\alpha_i}\neq0$ for all simple roots and there exists a simple root $\alpha$ so that $f_\rho^*\omega_{-\alpha}$ is nowhere vanishing.  
Suppose $(\rho_t,J_{\rho_t})$ be a one parameter family with $[\left.\frac{d}{dt}\rho_t\right|_{t=0}]=0.$ Let 
\[f_t:\widetilde{(S,J_t)}\ra\sSp(4,\C)/(\sU(1)\times\sU(1))\]
be the family of $\rho_t$-equivariant $\sSp(4,\R)$-cyclic surfaces associated to the family Higgs bundles determined by $(\rho_t,J_t).$  
Then by Theorem \ref{RigidityOfCyclicSurf}, $\left.\frac{d}{dt}J_t\right|_{t=0}=0,$ proving local uniqueness of $J_\rho.$
\end{proof}

Recall that Theorem \ref{uniqueminThm} says that, for $g-1<d\leq3g-3$ and $\rho\in\Rr_d(\sSp(4,\R)$, there is a unique conformal structure in which the $\rho$-equivariant harmonic map is a minimal immersion. To prove it, we need to go from the local uniqueness of Theorem \ref{LocalUniqueminThm} to global uniqueness. We follow Labourie's general differential geometric arguments in section 8 of \cite{cyclicSurfacesRank2}. 
\begin{Theorem}\label{DiffGeomThmLabourie}
	(Theorem 8.1.1 \cite{cyclicSurfacesRank2})Let $\pi:P\ra M$ be a smooth fiber bundle with connected fibers and $F:P\ra \R$ be a positive smooth function. Define
	\[N=\{x\in P\ |\ \left.d_x (F\right|_{P_{\pi(x)}})=0\}.\]
	and assume for all $m\in M$ the function $\left.F\right|_{P_m}$ is proper and that $N$ is connected and everywhere transverse to the fibers. Then $\pi$ is a diffeomeorphism from $N$ onto $M$ and $\left.F\right|_{P_m}$ has a unique critical point which is an absolute minimum.
\end{Theorem}
 \begin{proof}
	(Of Theorem \ref{uniqueminThm}) By Remark \ref{StabilityofMdRemark}, the space $\Rr_d(\sSp(4,\R))$ is smooth if and only if $g-1<d\leq 3g-3.$ Consider the fiber bundle $\pi:\Tt(S)\times\Rr_d(\sSp(4,\R))\ra\Rr_d(\sSp(4,\R)).$ Define a positive function $F$ by 
	\[F((J,\rho))=\Ee_\rho(J)=\Ee_J(h_\rho)=\haf\int\limits_{S}|dh_\rho|^2 dvol_J.\]
By \cite{CrossRatioAnosoveProperEnergy}, the map $\left.F\right|_{\rho}$ is proper and smooth, furthermore, the critical points of $\left.F\right|_{P_\rho}$ are minimal surfaces. Set 
\[N=\{(J,\rho)\in P\ |\ \left.d_x (F\right|_{P_{\rho}})=0\}.\]
By Theorem \ref{LocalUniqueminThm}, $N$ is everywhere transverse to the fibers. Applying Theorem \ref{DiffGeomThmLabourie}, when $g-1<d\leq 3g-3,$ for each $\rho\in\Rr_d(\sSp(4,\R))$ there is a unique conformal structure $(S,J_\rho)=\Sigma$ in which the $\rho$-equivariant harmonic map 
\[h_\rho:\widetilde\Sigma\ra\sSp(4,\R)/\sU(2)\]
is a minimal immersion. 
\end{proof}
The following corollary is immediate from the above proof.
\begin{Corollary}\label{MCGinvParamCor}
	The map $\Psi$ in Theorem \ref{MCGinvparamTheorem} is a diffeomeorphism. 
\end{Corollary}

\bibliography{../../../mybib}{}
\bibliographystyle{amsalpha}

\end{document}